\documentclass[11pt]{article}
\usepackage[english]{babel}
\usepackage{amsmath,amsthm,amssymb}
\usepackage{mathrsfs}
\usepackage{bbm}
\usepackage{amsfonts}
\usepackage[margin=0.8in]{geometry}
\usepackage{graphicx}
\usepackage{caption}
\usepackage{subcaption}
\usepackage{wasysym}
\usepackage{enumerate}
\usepackage{algorithm2e}
\usepackage{tabularx}
\usepackage{color}

\newtheorem{thm}{Theorem}[section]
\newtheorem{ass}[thm]{Assumption}

\newtheorem{lem}[thm]{Lemma}
\newtheorem{prop}[thm]{Proposition}
\theoremstyle{definition}
\newtheorem{defn}[thm]{Definition}

\newtheorem{example}[thm]{Example}
\theoremstyle{remark}
\newtheorem{rem}[thm]{Remark}
\numberwithin{equation}{section}
\newcommand{\R}{\mathbb R}

\newcommand{\bbF}{\mathbb F}
\newcommand{\bbG}{\mathbb G}

\newcommand{\mcA}{\mathcal{A}}
\newcommand{\mcAabs}{\mathcal{A}^{\rm abs}}
\newcommand{\mcB}{\mathcal{B}}
\newcommand{\mcD}{\mathcal D}
\newcommand{\mcG}{\mathcal{G}}
\newcommand{\mcN}{\mathcal N}
\newcommand{\mcF}{\mathcal F}
\newcommand{\mcI}{\mathcal I}

\newcommand{\mcAI}{\mcJ} % kan behöva ändras
\newcommand{\mcJ}{\mathcal J}

\newcommand{\mcT}{\mathcal T}
\newcommand{\mcP}{\mathcal P}

\newcommand{\mcH}{\mathcal H}
\newcommand{\mcK}{\mathcal K}
\newcommand{\mcU}{\mathcal U}
\newcommand{\mcS}{\mathcal S}

\newcommand{\E}{\mathbb{E}}
\newcommand{\Prob}{\mathbb{P}}

\newcommand{\veca}{\bold{a}}
\newcommand{\vecb}{\bold{b}}

\newcommand{\vece}{\bold{e}}
\newcommand{\esssup}{\mathop{\rm{ess}\,\sup}}
\newcommand{\ett}{\mathbbm{1}}
\newcommand{\bigSET}{\mcD}

\newcommand{\cadlag}{c\`adl\`ag~}
\newcommand{\cadlagp}{c\`adl\`ag.~}

\newcommand{\ie}{\textit{i.e.\ }}
\newcommand{\eg}{\textit{e.g.\ }}

\begin{document}

\title{On the Finite Horizon Optimal Switching Problem with Random Lag\footnote{This work was supported by the Swedish Energy Authorities through grant }}%

\author{Magnus Perninge\footnote{M.\ Perninge is with the Department of Physics and Electrical Engineering, Linnaeus University, V\"axj\"o,
Sweden. e-mail: magnus.perninge@lnu.se.}} %
\maketitle
% ----------------------------------------------------------------
\begin{abstract}
We consider an optimal switching problem with random lag and possibility of component failure. The random lag is modeled by letting the operation mode follow a regime switching Markov-model with transition intensities that depend on the switching mode. The possibility of failures is modeled by having absorbing components. We show existence of an optimal control for the problem by applying a probabilistic technique based on the concept of Snell envelopes.
\end{abstract}

% ----------------------------------------------------------------
\section{Introduction}
The standard optimal switching problem (sometimes referred to as starting and stopping problem) is a stochastic optimal control problem of impulse type that arises when an operator controls a dynamical system by switching between the different members in a set of switching modes $\mcI=\{\vecb_1,\ldots,\vecb_m\}$. In the two-modes case ($m=2$) the modes may represent, for example, ``operating'' and ``closed'' when maximizing the revenue from mineral extraction in a mine as in~\cite{Brennan}. In the multi-modes case the operating modes may represent different levels of power production in a power plant when the owner seeks to maximize her total revenue from producing electricity~\cite{CarmLud} or the states ``operating'' and ``closed'' of single units in a multi-unit production facility as in~\cite{BrekkeOksendal}.

In optimal switching the control takes the form $u=(\tau_1,\ldots,\tau_N;\beta_1,\ldots,\beta_N)$, where $\tau_1\leq\tau_2\leq\cdots\leq\tau_N$ is a sequence of (random) times when the operator intervenes on the system and $\beta_j\in\mcI$ is the switching mode that the operator switches to at time $\tau_j$. The standard multi-modes optimal switching problem in finite horizon ($T<\infty$) can then be formulated as finding the control that maximizes
\begin{equation*}
\E\bigg[\int_0^T \psi_{\alpha_s}(s)ds+\Upsilon_{\alpha_T}-\sum_{j=1}^Nc_{\beta_{j-1},\beta_{j}}(\tau_j)\bigg],
\end{equation*}
where $\alpha_t=\beta_0\ett_{[0,\tau_{1})}(t)+\sum_{j=1}^N \beta_j\ett_{[\tau_{j},\tau_{j+1})}(t)$ is the \emph{operation mode} (when starting in a predefined mode $\beta_0\in\mcI$), $\psi_\vecb$ and $\Upsilon_\vecb$ are the running and terminal revenue in mode $\vecb\in\mcI$, respectively and $c_{\vecb,\vecb'}(t)$ is the cost of switching from mode $\vecb$ to mode $\vecb'$ at time $t\in [0,T]$.

The standard optimal switching problem has been thoroughly investigated in the last decades after being popularised in~\cite{Brennan}. In~\cite{HamJean} a solution to the two-modes problem was found by rewriting the problem as an existence and uniqueness problem for a doubly reflected backward stochastic differential equation. In~\cite{BollanMSwitch1} existence of an optimal control for the multi-modes optimal switching problem was shown by a probabilistic method based on the concept of Snell envelopes. Furthermore, existence and uniqueness of viscosity solutions to the related Bellman equation was shown for the case when the switching costs are constant and the underlying uncertainty is modeled by a stochastic differential equation (SDE) driven by a Brownian motion. In~\cite{ElAsri} the existence and uniqueness results of viscosity solutions was extended to the case when the switching costs depend on the state variable. Since then, results have been extended to Knightian uncertainty~\cite{HuTang,HamZhang,ChassElieKharr} and non-Brownian filtration and signed switching costs~\cite{MartyrSigned}. For the case when the underlying uncertainty can be modeled by a diffusion process, generalization to the case when the control enters the drift and volatility term was treated in \cite{EliKharrCV}. This was further developed to include state constraints in~\cite{KharrSC}. Another important generalization is to the case when the operator only has partial information about the present state of the diffusion process as treated in~\cite{UUganget}.

As many physical systems do not immediately respond to changes in the control variables, including delays is an important aspect when seeking to derive applicable results in optimal control. General impulse control problems with lag have been considered in a variety of different settings including the novel paper~\cite{BarIlan95}, where an explicit solution to an inventory problem with uniform delivery lag is found by taking the current stock plus pending orders as one of the states. Similar approaches are taken in \cite{RAid2} where explicit optimal solutions of impulse control problems with uniform delivery lags are derived for a large set of different problems and in~\cite{Bruder} where an iterative algorithm is proposed. In~\cite{OksenImpulse} the authors propose a solution to general impulse control problems with lag, by defining an operator that circumvents the delay period. The optimal switching problem with non-uniform (but deterministic) lag and ramping was solved in~\cite{LimFEED} by state space augmentation in combination with the probabilistic approach initially developed in~\cite{BollanMSwitch1}.

The aim of the present article is to extend the applicability of optimal switching further by considering the case of random lag and component failure during startup. As in~\cite{LimFEED} we consider the problem of operating $n>0$ different production units, that can be either in operation or turned off, and thus let the switching modes be the set of all $n$-dimensional vectors of zeroes and ones, \ie $\mcI:=\{0,1\}^n$. To model the random lags and failures we let the operation mode, $\alpha^u_t$, be a continuous-time, finite-state, observable Markov-process taking values in $\mcA:=\{-1,0,1\}^n$, where $-1$ represents ``malfunction'', 0 represents ``off'' and 1 represents ``operating''. We assume that the transition intensities of $\alpha^u_t$ depend on the control both through the \emph{switching mode}, $\xi_t:=\sum_{j=1}^N\beta_j\ett_{[\tau_j,\tau_{j+1})}(t)$, but also through the time of the last switch from off to operating in each of the different production units. As opposed to the situation in the standard optimal switching problem, the switching mode and the operation mode may thus differ due to the lag.

We will consider the problem of finding a strategy $u$ that maximizes
\begin{align}
J(u):=\E\bigg[\Upsilon_{\alpha^u_T}(\theta_T^{u})+\int_0^T\psi_{\alpha^u_s}(s,\theta_s^{u})ds -\sum_{j}c_{\beta_{j-1},\beta_j}(\tau_j)\bigg].\label{ekv:objfun}
\end{align}
where the process $\theta^u$ is such that the $i^{\rm th}$ component gives the elapsed time in the present ``on''-cycle for Plant $i$. The process $\theta^u$ will allow us to model increased production costs during startup or lower production during ramp-up periods (see \eg \cite{ObactMMOR} for a situation where ramping is important). The results presented will be derived under the assumption that the $\psi_\veca$, $\Upsilon_\veca$ and $c_{\vecb,\vecb'}$ are adapted w.r.t.~a filtration generated by a Brownian motion. However, these results readily extend to more general (quasi-left continuous) filtrations, \eg a filtration generated by a Brownian motion and an independent Poisson random measure.

The remainder of the article is organized as follows. In the next section we state the problem, set the notation used throughout the article and detail the set of assumptions that are made. Then, in Section~\ref{sec:VERthm} a verification theorem is derived. This verification theorem is an extension of the original verification theorem for the multi-modes optimal switching problem developed in~\cite{BollanMSwitch1}. In Section~\ref{sec:exist} we show that there exists a family of processes that satisfies the requirements of the verification theorem, thus proving existence of an optimal control for the optimal switching problem with random lag. Then, in Section~\ref{sec:ValFun} we focus on the case when the underlying uncertainty in the processes $\psi_\veca$ and $\Upsilon_\veca$ can be modeled by an SDE and derive a dynamic programming relation for the corresponding value functions.

\section{Preliminaries}
We consider the finite horizon problem and thus assume that the terminal time $T$ is fixed with $T<\infty$. We will assume that turning off a unit gives immediate results on the operation mode and we have $\alpha_t\leq \xi_t$ for all $t\in[0,T]$. The state space for $(\alpha,\xi)$ is then $\mcAI:=\{(\veca,\vecb)\in\mcA\times\mcI:\veca\leq \vecb\}$. Furthermore, we define the following sets\footnote{Throughout, we assume that $\leq$, $\wedge$ and $^+$ are defined componentwise so that, for any two vectors $x,y\in\R^n$, $x\leq y$ implies that $x_i\leq y_i$, $[x\wedge y]_i=\min(x_i,y_i)$ and $[x^+]_i=\max(0,x_i)$, for $i=1,\ldots,n$.}:
\begin{itemize}
  \item For each $\vecb\in\mcI$, we let $\mcA_{\vecb}:=\{\veca\in\mcA:\veca\leq\vecb\}$ and for each $\veca\in\mcA$ we let $\mcI_{\veca}:=\{\vecb\in\mcI:\vecb\geq \veca\}$.
  \item For each $(\veca,\vecb)\in\mcAI$ we let $\mcA_{\veca,\vecb}:=\{\veca'\in\mcA_\vecb: |a_i'|\geq |a_i|\:\:{\rm and}\:\:a'_i = a_i\:\:{\rm when}\:\: a_i\in\{-1,-b_i\}\}$.
  \item For each $\vecb\in\mcI$ we let $\mcI^{-\vecb}:=\mcI\setminus\{\vecb\}$ and for each $\veca'\in\mcA_{\veca,\vecb}$ we let $\mcA^{-\veca'}_{\veca,\vecb}:=\mcA_{\veca,\vecb}\setminus \{\veca'\}$
  \item For each $\vecb\in\mcI$ we let $\mcA_{\vecb}^{\rm abs}:=\Pi_{i=1}^{n}\{\{-b_i\}\cup \{-1\}\}$.
  \item We introduce the set\footnote{For $\vecb=(b_1,\ldots,b_n)\in\mcI$ we let $[0,T]^{\vecb}:=[0,b_1T]\times\cdots\times[0,b_nT]$} $\bigSET:=[0,T] \times \cup_{(\veca,\vecb)\in\mcAI}[0,T]^{\vecb} \times[0,T]^{\veca^+}\times(\veca,\vecb)$ and let $\mcD_A:=[0,T] \times \cup_{(\veca,\vecb)\in\mcAI}[0,T]^{\vecb}\times(\veca,\vecb)$ and $\mcD_\lambda:= \cup_{\vecb\in\mcI}[0,T]^{\vecb}\times\{\vecb\}$. Furthermore, for each $(\veca,\vecb)\in\mcJ$ we let $\mcD_{(\veca,\vecb)}:=[0,T]\times [0,T]^\vecb \times[0,T]^{\veca^+}$.
\end{itemize}
Note here that $\mcA_{\veca,\vecb}$ is the set of all $\veca'\in\mcA_\vecb$ that the operation mode may transition to from $\veca$ when $\xi=\vecb$ and $\mcA_{\vecb}^{\rm abs}$ is the set of all states in $\mcA$ that are absorbing for $\alpha$ when $\xi=\vecb$.\\

We let $(\Omega,\mcG,\bbG,\Prob)$ be a probability space endowed with a $d$-dimensional Brownian motion $(B_t:0\leq t\leq T)$ whose augmented natural filtration is $\bbF:=(\mcF_t)_{0\leq t\leq T}$. For all $(t,\nu,\veca,\vecb)\in\mcD_A$ let $(A^{t,\nu,\veca,\vecb}_s:0\leq s\leq T)$ be a mixed Markov chain (sometimes also referred to as a stochastic hybrid system~\cite{LygerosSHS}) with state-space $\mcA_{\veca,\vecb}$. We assume that $A^{t,\nu,\veca,\vecb}_s=\veca$ for $s\in[0,t\vee \max_i \nu_i]$ and that on $(t\vee \max_i \nu_i,T]$, the transition rate from $\veca$ to $\veca'$, $\lambda^{\nu,\vecb}_{\veca,\veca'}(s)$, is $\bbF$-progressively measurable for all $\veca'\in\mcA_{\veca,\vecb}$. For $\veca'\notin\mcA_{\veca,\vecb}$ we let $\lambda^{\nu,\vecb}_{\veca,\veca'}(s)\equiv 0$. We assume that $\bbG:=(\mcG_t)_{0\leq t\leq T}$ is the augmented natural filtration generated by $B$ and the family $((A^{t,\nu,\veca,\vecb}_s)_{0\leq s\leq T}:(t,\nu,\veca,\vecb)\in\mcD_A)$, satisfying the usual conditions (more information about enlargement of filtrations can be found in \eg Chapter 6 of~\cite{Protter}).\\

Recall here the concept of left continuity in expectation: A process $(X_t:0\leq t\leq T)$ is left continuous in expectation (LCE) if for each stopping time $\gamma$ and each sequence of stopping times $\gamma_k\nearrow\gamma$ we have $\lim\limits_{k\to\infty}\E\left[X_{\gamma_k}\right] = \E\left[X_\gamma\right]$.\\

\noindent Throughout we will use the following notation:
\begin{itemize}
  \item $\mcP_{\bbF}$ (resp. $\mcP_{\bbG}$) is the $\sigma$-algebra of $\bbF$-progressively measurable ($\bbG$-progressively measurable) subsets of $[0,T]\times \Omega$.
  \item We let $\mcS^{2}$ be the set of all $\R$-valued, $\mcP_{\bbG}$-measurable, \cadlag processes $(Z_t: 0\leq t\leq T)$ such that $\E\left[\sup_{t\in[0,T]} |Z_t|^2\right]<\infty$. We let $\mcS_{e}^{2}$ be the subset of processes that are non-negative and LCE.
  \item We let $\mcS_\bbF^{2}$ be the subset of $\mcS^{2}$ of processes that are $\mcP_{\bbF}$-measurable and let $\mcS_{\bbF,c}^{2}$ be the subset of $\mcS_\bbF^{2}$ of all processes that are continuous.
  \item We let $\mcH^{2}$ denote the set of all $\R$-valued $\mcP_{\bbF}$-measurable processes $(Z_t: 0\leq t\leq T)$ such that $\E\left[\int_0^T |Z_t|^2 dt\right]<\infty$.
  \item We let $\mcH^{\infty}_\bbF$ denote the set of all $\R$-valued $\mcP_{\bbF}$-measurable processes $(Z_t: 0\leq t\leq T)$ such that $|Z_t|<\infty$, $d\Prob\otimes dt$-a.e.
  \item We let $\mcT$ ($\mcT^{\bbF}$) be the set of all $\bbG$-($\bbF$-)stopping times and for each $\gamma\in\mcT$ ($\mcT^{\bbF}$) we let $\mcT_\gamma$ ($\mcT^{\bbF}_\gamma$) be the subset of stopping times $\tau$ such that $\tau\geq \gamma$, $\Prob$-a.s.
  \item We let $\mcU$ be the set of all $u=(\tau_1,\ldots,\tau_N;\beta_1,\ldots,\beta_N)$, where $(\tau_j)_{j=1}^N$ is an increasing sequence of $\bbG$-stopping times and $\beta_j\in\mcI^{-\beta_{j-1}}$ is $\mcG_{\tau_j}$-measurable.
  \item We let $\mcU^f$ be the subset of controls $u\in\mcU$ for which $N$ is finite $\Prob$-a.s.~(\ie $\mcU^f:=\{u\in\mcU:\: \Prob\left[\{\omega\in\Omega : N(\omega)>k, \:\forall k>0\}\right]=0\}$) and $\mcU^k$ the subset of controls for which $N\leq k$, $\Prob$-a.s. For $\gamma\in\mcT$ let $\mcU_\gamma$ (resp.~$\mcU_\gamma^f$ and $\mcU_\gamma^k$) be the subset of $\mcU$ (resp.~$\mcU^f$ and $\mcU^k$) with $\tau_1\in\mcT_\gamma$.
  %\item For each $(\nu,\vecb)\in\cup_{\vecb\in\mcI}[0,T]^\vecb\times\{\vecb\}$ and each $t\in [0,T]$, we let $\Lambda^{t,\nu,\vecb}(s):=\int_t^s\lambda^{\nu,\vecb}(r)dr$ for $0\leq s \leq T$.
\end{itemize}
Our problem will be characterised by four objects:
\begin{itemize}
  \item A collection $(\Upsilon_{\veca}:\Omega\times [0,T]^{\veca^+}\to\R)_{\veca\in\mcA}$ of $\mcF_T\otimes \mcB([0,T]^{\veca^+})$-measurable maps.
  \item A collection $(\psi_{\veca}:\Omega\times [0,T]\times [0,T]^{\veca^+}\to\R)_{\veca\in\mcA}$ where $\psi_{\veca}$ is a $\mcP_{\bbF}\otimes \mcB([0,T]^{\veca^+})$-measurable map.
  \item A cost process $C^u_t:=\sum_{\tau_j\leq t}c_{\beta_{j-1},\beta_j}(\tau_j)$, where $(c_{\vecb,\vecb'}:\Omega\times [0,T]\to\R)_{(\vecb,\vecb')\in\mcI^2}$ is a collection of $\mcP_{\bbF}$-measurable processes.
  \item A family $(((\lambda^{\nu,\vecb}_{\veca,\veca'}(s):0\leq s\leq T)_{(\veca,\veca')\in\mcA_{\vecb}\times\mcA_{\vecb}})_{\nu\in[0,T]^\vecb})_{\vecb\in\mcI}$ of $\R$-valued, $\mcP_{\bbF}$-measurable transition intensities, \ie $\lambda^{\nu,\vecb}_{\veca,\veca'}\geq 0$ for $\veca'\neq\veca$ and $\lambda^{\nu,\vecb}_{\veca,\veca}=-\sum_{\veca'\in\mcA_{\veca,\vecb}^{-\veca}}\lambda^{\nu,\vecb}_{\veca,\veca'}$.
\end{itemize}

We make the following assumptions:%(where all statements hold $\Prob$-a.s.~unless stated otherwise)
\begin{ass}\label{ass:prelim}
\begin{enumerate}[(i)]
  \item\label{ass:onpsi} For each $\veca\in\mcA$, $\psi_\veca(\cdot,0)\in \mcS_{\bbF}^{2}$ and $\E[\Upsilon_\veca^2]<\infty$. Furthermore, we assume that there are constants $k_\psi>0$ and $k_\Upsilon>0$ such that, for each $(z,z')\in [0,T]^{\veca}\times [0,T]^{\veca}$,
  \begin{equation*}
    |\psi_\veca(t,z)-\psi_\veca(t,z')|\leq k_\psi|z-z'|,
  \end{equation*}
  for all $t\in [0,T]$ and
  \begin{equation*}
    |\Upsilon_{\veca}(z)-\Upsilon_{\veca}(z')|\leq k_\Upsilon |z-z'|,
  \end{equation*}
  $\Prob$-a.s. (where the exception set does not depend on the tuple $(t,z,z')$).
  \item\label{ass:onc} The switching costs $(c_{\vecb,\vecb'})_{\vecb,\vecb'\in \mcI}\in (\mcS_{\bbF,c}^2)^{m\times m}$ are such that, $\Prob$-a.s.,
  \begin{enumerate}[(a)]
    \item $\inf_{t\in [0,T]} c_{\vecb,\vecb'}(t)\geq 0$,\quad for all $(\vecb,\vecb')\in \mcI\times\mcI$
    \item $c_{\vecb_{1},\vecb_2}(t_1)+c_{\vecb_2,\vecb_3}(t_2)+\ldots+c_{\vecb_{k-1},\vecb_k}(t_{k-1})+c_{\vecb_{k},\vecb_1}(t_k)\geq\epsilon>0$,\quad for all $0\leq t_1\leq t_2\leq\cdots\leq t_k\leq T$ and $(\vecb_1,\ldots,\vecb_k)\in \mcI^k$.
  \end{enumerate}
  \item\label{ass:@end} For all $\veca\in \mcA$ and $z\in [0,T]^{\veca^+}$, $\Upsilon_{\veca}(z)>\max_{(\vecb,\vecb')\in \mcI_{\veca}\times\mcI} \{\Upsilon_{\veca\wedge\vecb'}(z\wedge T\vecb')-c_{\vecb,\vecb'}(T)\}$, $\Prob$-a.s.
  \item\label{ass:lambda} For each $(\nu,\vecb)\in\mcD_\lambda$ and all $(\veca,\veca')\in\mcA_{\vecb}$ the process $\lambda^{\nu,\vecb}_{\veca,\veca'}(\cdot)\in \mcH_\bbF^{\infty}$ and we assume that there is a constant $K_\lambda>0$, such that
      \begin{equation*}
        |\lambda^{\nu,\vecb}_{\veca,\veca'}(s)|\leq K_\lambda,
      \end{equation*}
      $d\Prob\otimes ds$-a.e. Furthermore, we assume that each element of $\lambda^{\nu,\vecb}$ is Lipschitz continuous in $\nu$:
      \begin{equation*}
        |\lambda^{\nu,\vecb}(s)-\lambda^{\nu',\vecb}(s)|\leq k_\lambda |\nu-\nu'|,
      \end{equation*}
      for all $s\in [0,T]$, $\Prob$-a.s. (where the exception set does not depend on the tuple $(s,\nu,\nu')$).
\end{enumerate}
\end{ass}

The above assumptions are mainly standard assumptions for optimal switching problems. Assumptions \ref{ass:onpsi} and \ref{ass:onc}.a together imply that the expected maximal reward is finite. Assumption \ref{ass:onc}.b implies that there is always a positive switching cost associated to making a loop of switches and \ref{ass:@end} implies that it is never optimal to switch at time $T$.\\

Each control $u=(\tau_1,\ldots,\tau_N;\beta_1,\ldots,\beta_N)$ defines the switching mode starting in $\vecb\in\mcI$ which is a process $(\xi^{\vecb}_t: 0\leq t\leq T)$ given by
\begin{equation*}
\xi^{\vecb}_t:=\vecb\ett_{[0,\tau_1)}(t)+\sum_{j=1}^N \beta_j\ett_{[\tau_{j},\tau_{j+1})}(t),
\end{equation*}
with $\tau_{N+1}=\infty$ (for notational simplicity we will write $\xi$ for $\xi^{0}$). The switching mode thus, in some sense, tells us the preferred operation state. For each initial mode $\vecb$ and each vector $\nu\in [0,T]^{\vecb}$, we let the control $u$ define the sequence\footnote{For two vectors $x,y\in\R^n$ we define the product $xy$ as $[xy]_i=x_iy_i$, where $[x]_i$ denotes the $i^{\rm th}$ component of the vector $x$.} $(\vartheta^{\nu,\vecb}_0,\ldots,\vartheta^{\nu,\vecb}_N)$ with $\vartheta^{\nu,\vecb}_0:=\nu$, $\vartheta^{\nu,\vecb}_1:=\nu\beta_1+\tau_1(\beta_1-\vecb)^+$ and then recursively $\vartheta^{\nu,\vecb}_{j}:=\vartheta^{\nu,\vecb}_{j-1}\beta_{j}+\tau_{j}(\beta_{j}-\beta_{j-1})^+$ for $j=2,\ldots,N$.\\

Given that the operating mode at time $t$ is $\veca$, the switching mode is $\vecb\in \mcI$ and given the vector of activation times $\nu$, such that $(t,\nu,\veca,\vecb)\in\mcD_A$, the family of mixed Markov-chains $(A^{\cdot,\cdot,\cdot,\cdot}_s:0\leq s\leq T)$ defines the sequence of operating modes, $\bar\alpha_0,\ldots,\bar\alpha_N$, at the intervention times as $\bar\alpha_0:=\veca$ and then recursively
\begin{align*}
\bar\alpha_j:=A^{\tau_{j-1},\vartheta^{\nu,\vecb}_{j-1},\bar\alpha_{j-1},\beta_{j-1}}_{\tau_j}\wedge \beta_j,
\end{align*}
for $j=1,\ldots,N$, with $\tau_0=t$ and $\beta_0=\vecb$. This leads us to define the operating mode $(\alpha^{t,\nu,\veca,\vecb,u}(s):0\leq s\leq T)$ as
\begin{align*}
\alpha^{t,\nu,\veca,\vecb,u}_s:=A^{t,\nu,\veca,\vecb}_s\ett_{[0,\tau_{1})}(s)+\sum_{j=1}^N A^{\tau_j,\vartheta^{\nu,\vecb}_{j},\bar\alpha_j,\beta_j}_s\ett_{[\tau_j,\tau_{j+1})}(s).
\end{align*}
For notational simplicity we use the same shorthand as above and write $\alpha^u$ for $\alpha^{0,0,0,0,u}$.\\

%\begin{example}
%Consider a situation with two units (\ie $\mcI=\{(0,0),(0,1),(1,0),(1,1)\}$): Plant 1 that has a constant rate 3 of starting and failure rate 0.1 once in operation and Plant 2 that has starting rate 2 times ``time since the unit was turned on'' and never fails. Then $\lambda^{\cdot,\cdot}_{\cdot,(\cdot,-1)}\equiv 0$ as Plant 2 never fails and $\lambda^{\nu,(b_1,1)}_{(a_1,0),(a_1,1)}(r)=2(r-\nu_2)^+$ for $-1\leq a_1\leq b_1$.
%\begin{table}[h!]\centering \setlength{\tabcolsep}{2pt}
%\begin{tabular}{| r | c | c | c | c | c | c |}\hline
%  $\veca\setminus \veca'$ & $(-1,0)$ & $(-1,1)$ & $(0,0)$ & $(0,1)$ & $(1,0)$ & $(1,1)$  \\
%  \hline
%   $(-1,0)$           &     $-2(r-\nu_2)^+$    & $2(r-\nu_2)^+$ &     0   &    0    &     0   &     0    \\ \hline
%   $(-1,1)$               &     0    &     0    &     0   &    0    &     0   &     0    \\ \hline
%   $(0,0)$                &     0    &     0    &     $-3-2(r-\nu_2)^+$   &   $2(r-\nu_2)^+$    &     3   &     0    \\ \hline
%   $(0,1)$                &     0    &     0    &     0   &    $-3$    &     0   &     3    \\ \hline
%   $(1,0)$                &    0.1   &     0    &     0   &    0    &    $-0.1-2(r-\nu_2)^+$   &     $2(r-\nu_2)^+$    \\ \hline
%   $(1,1)$                &     0    &    0.1    &     0   &    0    &     0   &     $-0.1$    \\ \hline
%\end{tabular}
%\caption{Transition rates $\lambda^{\nu,(1,1)}_{\veca,\veca'}(r)$.}\label{tab:exDATA1}
%\end{table}
%\end{example}

\begin{example}
Consider the case when transitions for different plants are independent and Plant $i$ has a failure rate $r^{\rm fail}_i\geq 0$ if operating and 0 otherwise and a startup rate $r^{\rm start}_i:[0,T]\to\R_ +$, where the input is the time that has elapsed since the unit was turned on. Then\footnote{Where $\vece_i:=(0,\ldots,0,1,0,\ldots,0)$ with a 1 in the $i^{\rm th}$ position.},
\begin{equation*}
\lambda^{\nu,\vecb}_{\veca,\veca'}(s)=\sum_{i=1}^n(\ett_{[\veca'=\veca+\vece_i]}r^{\rm start}_i(s-\nu_ i) + \ett_{[\veca'=\veca-2\vece_i]}r^{\rm fail}_i)
\end{equation*}
when $\veca'\in\mcA_{\veca,\vecb}^{-\veca}$ and $\lambda^{\nu,\vecb}_{\veca,\veca'}(s)=-\sum_{\veca'\in \mcA_{\veca,\vecb}^{-\veca}}\lambda^{\nu,\vecb}_{\veca,\veca'}(s)$.\qed
\end{example}

We let $((\theta^{t,\nu,z,\veca,\vecb}_s)_{0\leq s\leq T}:(t,\nu,z,\veca,\vecb)\in \bigSET)$ be given by
\begin{align*}
\theta^{t,\nu,z,\veca,\vecb}_s:=\Big(z+\int_t^{s\vee t} (A^{t,\nu,\veca,\vecb}_r)^+dr\Big)\wedge T\alpha^{t,\nu,\veca,\vecb}_s.
\end{align*}
To define the \emph{time in present on-mode} for a control $u=(\tau_1,\ldots,\tau_N;\beta_1,\ldots,\beta_N)\in\mcU_t$ starting in $z$ at time $t$ we let $\theta^0_s:=\theta^{t,\nu,z,\veca,\vecb}_s$ and then recursively define
\begin{equation*}
\theta^j_s:=\theta^{\tau_j,\vartheta^{\nu,\vecb}_{j},\theta^{j-1}_{\tau_j}\wedge T\beta_j,\bar\alpha_j,\beta_j}_s.
\end{equation*}
This allows us to define
\begin{align*}
\theta^{t,\nu,z,\veca,\vecb,u}_s:=\theta^0_s\ett_{[0,\tau_1)}(s)+\sum_{j=1}^N\theta^j_s\ett_{[\tau_j,\tau_{j+1})}(s)
\end{align*}
and again we let $\theta^u:=\theta^{0,0,0,0,0,u}$.\\

\begin{rem}
Note that, with $\vartheta_t:=\sum_{j=1}^N \vartheta_j\ett_{[\tau_{j},\tau_{j+1})}$ the set $\mcD$ is the state-space for $(s,\vartheta_s,\theta_s,\alpha_s,\xi_s)_{0\leq s\leq T}$, $\mcD_A$ is the state-space for $(s,\vartheta_s,\alpha_s,\xi_s)_{0\leq s\leq T}$ and $\mcD_\lambda$ is the state-space for $(\vartheta_s,\xi_s)_{0\leq s\leq T}$.
\end{rem}

%We end this section by showing a stability result:
%\begin{lem}
%\begin{equation}
%|\E\bigg[\int_0^TA^{t,\nu,\veca,\vecb}_r-A^{t',\nu',\veca,\vecb}_rdr\bigg]|\leq k_A(|t-t'|+|\nu-\nu'|)
%\end{equation}
%\end{lem}
%
%\noindent\emph{Proof.} We have
%\begin{align*}
%|\E\bigg[\int_0^TA^{t,\nu,\veca,\vecb}_r-A^{t',\nu',\veca,\vecb}_rdr\bigg]|&\leq n\sum_{\veca'\in\mcA_{\veca,\vecb}}|\E\bigg[\int_0^T\Big(\ett_{[A^{t,\nu,\veca,\vecb}_r=\veca']}-\ett_{[A^{t',\nu',\veca,\vecb}_r=\veca']}\Big)dr\bigg]|
%\\
%&=n\sum_{\veca'\in\mcA_{\veca,\vecb}}|\int_0^T\Prob\Big[A^{t,\nu,\veca,\vecb}_r=\veca'\Big]-\Prob\Big[A^{t',\nu',\veca,\vecb}_r=\veca'\Big]|
%\\
%&=n\sum_{\veca'\in\mcA_{\veca,\vecb}}|\int_0^T\E\Big[\big[e^{\Lambda^{t,\nu,\veca,\vecb}(r)}-e^{\Lambda^{t',\nu',\veca,\vecb}(r)}\big]_{\veca,\veca'}\Big]dr|
%\\
%&=n\sum_{\veca'\in\mcA_{\veca,\vecb}}|\int_0^T\E\Big[\big[e^{\Lambda^{t,\nu,\veca,\vecb}(r)-\Lambda^{t',\nu',\veca,\vecb}(r)}(e^{\Lambda^{t',\nu',\veca,\vecb}(r)}-e^{-\Lambda^{t,\nu,\veca,\vecb}(r)+2\Lambda^{t',\nu',\veca,\vecb}(r)})\big]_{\veca,\veca'}\Big]dr|
%\\
%&\leq C\sum_{\veca'\in\mcA_{\veca,\vecb}}\int_0^T\Big(\E\Big[\big[e^{2\Lambda^{t,\nu,\veca,\vecb}(r)-2\Lambda^{t',\nu',\veca,\vecb}(r)}\big]_{\veca,\veca'}\Big]\Big)^{1/2}dr
%\\
%&\leq C
%\end{align*}

\noindent We are now ready to state the optimal switching problem with random lag:\\

\noindent\textbf{Problem 1.} Find $u^*\in\mcU$, such that
\begin{equation}\label{ekv:OPTprob}
J(u^*)=\sup_{u\in\mcU} J(u).
\end{equation}
\qed

\bigskip

\begin{rem}\label{rem:pos}
Note that we have
\begin{align*}
J(u)=&\:\E\bigg[\Upsilon_{\alpha^u_T}(\theta_T^{u}) -\min_{\veca\in\mcA}\min_{z\in [0,T]^{\veca^+}}\Upsilon_{\veca}(z)+ \int_0^T(\psi_{\alpha^u_s}(s,\theta_s^{u}) -\min_{\veca\in\mcA}\min_{z\in [0,T]^{\veca^+}}\psi_{\veca}(s,z))ds
\\
&-\sum_{j}c_{\beta_{j-1},\beta_j}(\tau_j)\bigg]+\E\bigg[\min_{\veca\in\mcA}\min_{z\in [0,T]^{\veca^+}}\Upsilon_{\veca}(z) + \int_0^T\min_{\veca\in\mcA}\min_{z\in [0,T]^{\veca^+}}\psi_{\veca}(s,z)ds\bigg].
\end{align*}
Hence, we can without loss of generality assume that for each $\veca\in\mcA$, $\Upsilon_{\veca}$ and $\psi_{\veca}$ are both non-negative.
\end{rem}

The following proposition is a standard result for optimal switching problems and is due to the ``no-free-loop'' condition.
\begin{prop}\label{prop:finSTRAT} Suppose that there is a $u^*\in\mcU$ such that $J(u^*)\geq J(u)$ for all $u\in\mcU$. Then $u^*\in\mcU^f$.
\end{prop}

\noindent\emph{Proof.} Assume that $u\in \mcU\setminus \mcU^f$ and let $B:=\{\omega \in\Omega: N(\omega)>k, \:\forall k>0\}$, then $\Prob[B]>0$. Furthermore, if $B$ holds then the switching mode must make an infinite number of loops and we have
\begin{align*}
J(u)&\leq \E\bigg[ \int_0^T \max_{\veca\in\mcA} \sup_{z\in [0,T]^{\veca}} |\psi_{\veca}(s,z)|ds- \infty \ett_B \epsilon\bigg]
=-\infty,%\min_{0\leq t_1\leq t_2\leq\cdots\leq t_k}\min_{{\vecb_1,\ldots,\vecb_k}}(c_{\vecb_1,\vecb_2}(t)+\cdots+c_{\vecb_{k-1},\vecb_k}(t)+c_{\vecb_{k},\vecb_1}(t))
\end{align*}
by Assumption~\ref{ass:prelim} (\ref{ass:onpsi}) and (\ref{ass:onc}). Now, by the above non-negativity assumption on $\Upsilon$ and $\psi$ we have $J(u)\geq 0$ for $u=\emptyset$ and the assertion follows.\qed

\bigskip

We end this section with two useful lemmas:
\begin{lem}\label{lemma:AisLCE}
Let $(\gamma_m)_{m\geq 1}$ be a sequence of $\bbG$-stopping times such that $\gamma_m\nearrow\gamma$ $\Prob$-a.s.~for some $\gamma\in\mcT$. Then, for any $(t,\nu,\veca,\vecb)\in \mcD_A$
\begin{equation*}
\lim_{m\to\infty}\E\left[g(A^{t,\nu,\veca,\vecb}_{\gamma_m})\right]= \E\left[g(A^{t,\nu,\veca,\vecb}_{\gamma})\right]
\end{equation*}
for all $\mcF_\gamma$-measurable functions $g$ such that $\sum_{\veca\in\mcA}\E[|g(\veca)|^2]<\infty$.
\end{lem}

\noindent\emph{Proof.} We have\footnote{Throughout, $C\in(0,\infty)$ is a constant that may change value from line to line.}
\begin{align*}
\E\left[|g(A^{t,\nu,\veca,\vecb}_{\gamma_m})-g(A^{t,\nu,\veca,\vecb}_{\gamma})|\right]&\leq C\E\left[\sum_{\veca'\in\mcA}|g(\veca')|^2\right]^{1/2}\Prob(A^{t,\nu,\veca,\vecb}_{\gamma_m}\neq A^{t,\nu,\veca,\vecb}_{\gamma})^{1/2}
\\
&\leq C\E\left[\sum_{\veca'\in\mcA}|g(\veca')|^2\right]^{1/2}\E\left[\int_{\gamma_m}^{\gamma}\max_{\veca'\in\mcA}|\lambda_{\veca',\veca'}^{\nu,\vecb}(s)| ds\right]^{1/2}
\\
&\leq C\E\left[{\gamma}-{\gamma_m}\right]^{1/2}
\end{align*}
and the last part goes to zero as $m\to\infty$.\qed\\

\begin{lem}\label{lemma:condINTisL1cont}
For any $(t,\nu,z,\veca,\vecb)\in \bigSET$ and any $s\in[t,T]$ we have
\begin{equation*}
\lim_{(\nu',z')\to(\nu,z)}\E\left[\int_t^s (\psi_{A^{t,\nu',\veca,\vecb}_{r}}(r,\theta^{t,\nu',z',\veca,\vecb}_r) - \psi_{A^{t,\nu,\veca,\vecb}_{r}}(r,\theta^{t,\nu,z,\veca,\vecb}_r))dr\big|\mcG_t\right]=0,
\end{equation*}
$\Prob$-a.s.
\end{lem}

\noindent\emph{Proof.} First note that by definition we have $\theta^{t,\nu,z',\veca,\vecb}_r-\theta^{t,\nu,z,\veca,\vecb}_r= z'-z$. Assumption~\ref{ass:prelim}.\ref{ass:onpsi} now gives
\begin{align*}
|\psi_{A^{t,\nu',\veca,\vecb}_{r}}(r,\theta^{t,\nu',z',\veca,\vecb}_r) - \psi_{A^{t,\nu,\veca,\vecb}_{r}}(r,\theta^{t,\nu,z,\veca,\vecb}_r)| &\leq |\psi_{A^{t,\nu',\veca,\vecb}_{r}}(r,\theta^{t,\nu',z,\veca,\vecb}_r) - \psi_{A^{t,\nu,\veca,\vecb}_{r}}(r,\theta^{t,\nu,z,\veca,\vecb}_r)|
\\
&\quad  + k_\psi |z'-z|,
\end{align*}
$\Prob$-a.s.~for all $r\in[0,T]$. If $\veca\in\mcAabs_\vecb$ then $A^{t,\nu,\veca,\vecb}_r=A^{t,\nu',\veca,\vecb}_r =\veca$, $\Prob$-a.s.~for all $r\in[0,T]$ and the result follows. Assume instead that $\veca\notin\mcAabs_\vecb$ and let $\eta$ and $\eta'$ be the first transition times of $A^{t,\nu,\veca,\vecb}$ and $A^{t,\nu',\veca,\vecb}$, respectively. For all $(t,\nu,z,\veca,\vecb)\in \bigSET$ and all $r\in[t,s]$ we let
\begin{equation*}
\Gamma^{t,\nu,z,\veca,\vecb}_r:=\E\bigg[\int_r^s\psi_{A^{t,\nu,\veca,\vecb}_{r}}(r,\theta^{t,\nu,z,\veca,\vecb}_r)dr\big|\mcG_r\bigg].
\end{equation*}
Then, with $\theta_r^{t,z,\veca}:=z+(r-t)^+\veca^+$, we have
\begin{align*}
&\Gamma^{t,\nu',z,\veca,\vecb}_t-\Gamma^{t,\nu,z,\veca,\vecb}_t =\E\left[\int_{\eta'}^s \psi_{A^{t,\nu',\veca,\vecb}_{r}}(r,\theta^{t,\nu',z,\veca,\vecb}_r)dr- \int_{\eta}^s \psi_{A^{t,\nu,\veca,\vecb}_{r}}(r,\theta^{t,\nu,z,\veca,\vecb}_r)dr\big|\mcG_t\right]
\\
&=\sum_{\veca'\in\mcA^{-\veca}_{\veca,\vecb}}\E\left[\int_t^s (\lambda_{\veca,\veca'}^{\nu',\vecb}(r)e^{\int_t^r\lambda_{\veca,\veca}^{\nu',\vecb}(v)dv} \Gamma^{r,\nu',\theta_r^{t,z,\veca}\wedge T\veca',\veca',\vecb}_r - \lambda_{\veca,\veca'}^{\nu,\vecb}(r) e^{\int_t^r\lambda_{\veca,\veca}^{\nu,\vecb}(v)dv} \Gamma^{r,\nu,\theta_r^{t,z,\veca}\wedge T\veca',\veca',\vecb}_r)dr\big|\mcG_t\right]
\end{align*}
where we have used the fact that $\Prob(\eta\geq r)=e^{\int_t^r\lambda_{\veca,\veca}^{\nu',\vecb}(v)dv}$ and that the transition-rate of $A^{t,\nu,\veca,\vecb}_{\cdot}$ from $\veca$ to $\veca'$ is $\lambda_{\veca,\veca'}^{\nu,\vecb}(\cdot)$. Using the identity $ab-a'b'=1/2((a-a')(b+b')+(a+a')(b-b'))$ we find that
\begin{align*}
&|\Gamma^{t,\nu',z,\veca,\vecb}_t-\Gamma^{t,\nu,z,\veca,\vecb}_t|
\\
&\leq \sum_{\veca'\in\mcA^{-\veca}_{\veca,\vecb}}\E\bigg[\int_t^s \Big\{ |\lambda_{\veca,\veca'}^{\nu',\vecb}(r)e^{\int_t^r\lambda_{\veca,\veca}^{\nu',\vecb}(v)dv}-\lambda_{\veca,\veca'}^{\nu,\vecb}(r) e^{\int_t^r\lambda_{\veca,\veca}^{\nu,\vecb}(v)dv}| |\Gamma^{r,\nu',\theta_r^{t,z,\veca}\wedge T\veca',\veca',\vecb}_r + \Gamma^{r,\nu,\theta_r^{t,z,\veca}\wedge T\veca',\veca',\vecb}_r|
\\
&\quad +|\lambda_{\veca,\veca'}^{\nu',\vecb}(r)e^{\int_t^r\lambda_{\veca,\veca}^{\nu',\vecb}(v)dv}+\lambda_{\veca,\veca'}^{\nu,\vecb}(r) e^{\int_t^r\lambda_{\veca,\veca}^{\nu,\vecb}(v)dv}| |\Gamma^{r,\nu',\theta_r^{t,z,\veca}\wedge T\veca',\veca',\vecb}_r - \Gamma^{r,\nu,\theta_r^{t,z,\veca}\wedge T\veca',\veca',\vecb}_r|\Big\}dr\big|\mcG_t\bigg]
\\
&\leq C\E\bigg[|\nu'-\nu|\int_0^T \max_{\veca'\in \mcA}\max_{\zeta\in[0,T]^{(\veca')^+}}|\psi_{\veca'}(r,\zeta)|dr + \sum_{\veca'\in\mcA^{-\veca}_{\veca,\vecb}}\int_t^s|\Gamma^{r,\nu',\theta_r^{t,z,\veca}\wedge T\veca',\veca',\vecb}_r - \Gamma^{r,\nu,\theta_r^{t,z,\veca}\wedge T\veca',\veca',\vecb}_r|dr \big|\mcG_t\bigg].
\end{align*}
Now, as $\Gamma^{r,\nu',z,\veca,\vecb}_r-\Gamma^{r,\nu,z,\veca,\vecb}_r=0$, $\Prob$-a.s., whenever $\veca\in\mcAabs_{\vecb}$ can use an induction argument to deduce that
\begin{align*}
&|\Gamma^{t,\nu',z,\veca,\vecb}_t-\Gamma^{t,\nu,z,\veca,\vecb}_t|
\leq C|\nu'-\nu|\E\bigg[\int_0^T \max_{\veca'\in \mcA}\max_{\zeta\in[0,T]^{(\veca')^+}}|\psi_{\veca'}(r,\zeta)|dr\big|\mcG_t\bigg]
\end{align*}
and the assertion follows as the last term is $\Prob$-a.s.~bounded by Assumption~\ref{ass:prelim}.\ref{ass:onpsi} and Doob's maximal inequality.

%Using Doob's martingale inequality in combination with an induction argument and the fact that $\E\Big[\sup_{r\in [0,T]}|\Gamma^{r,\nu',z,\veca,\vecb}_r-\Gamma^{r,\nu,z,\veca,\vecb}_r|^2\Big]=0$ for $\veca\in\mcAabs_{\vecb}$ we get
%\begin{align*}
%&\E\Big[\sup_{r\in [0,T]}|\Gamma^{r,\nu',z',\veca,\vecb}_r-\Gamma^{r,\nu,z,\veca,\vecb}_r|^2\Big]\leq C(|\nu'-\nu|+|z'-z|),
%\end{align*}
%for all $(\veca,\vecb)\in\mcAI$. As $\Gamma^{r,\nu',\zeta',\veca,\vecb}_r$ is $L^2$-bounded by Assumption~\ref{ass:prelim}.\ref{ass:onpsi}, the assertion now follows from the Vitali convergence theorem.\qed\\

%%%%%%%%%%%%%%%%%%%%%%%%%%%%%%%%%%%%%%%%%%%%%%%%%%%%%%%%%%%%%%%%%%%%%%%%%%%%%%%%%%%%%%%%%%%%%%%%%%%%%%%%%%%%%%%%%%%%%%%%%%%%%%%%%%%%%%%%%%%%%%%%%%

\subsection{The Snell envelope}
In this section we gather the main results concerning the Snell envelope that will be useful later on. When presenting the theory we introduce an auxiliary probability space $(\tilde\Prob,\tilde\Omega,\tilde\mcF,(\tilde\mcF_t)_{0\leq t\leq T})$ that we assume satisfies the usual conditions. For any $\tilde\bbF$-stopping time $\eta$, we let $\tilde\mcT_\eta$ be the set of $\tilde\bbF$-stopping times $\tau\geq \eta$, with $\tilde\bbF:=(\tilde\mcF_t)_{0\leq t\leq T}$ and recall that a progressively measurable process $U_t$ is of class [D] if the set of random variables $\{X_\tau:\tau\in\tilde\mcT_0\}$ is uniformly integrable.

\begin{thm}[The Snell envelope]\label{thm:Snell}
Let $U=(U_t)_{0\leq t\leq T}$ be an $\tilde \bbF$-adapted, $\R$-valued, \cadlag process of class [D]. Then there exists a unique (up to indistinguishability), $\R$-valued \cadlag process $Z=(Z_t)_{0\leq t\leq T}$ called the Snell envelope, such that $Z$ is the smallest supermartingale that dominates $U$. Furthermore, the following holds:
\begin{enumerate}[(i)]
  \item\label{Snell:sup} For any stopping time $\gamma$,
    \begin{equation}\label{ekv:SnellZ}
      Z_{\gamma}=\esssup_{\tau\in \tilde\mcT_{\gamma}}\E\left[U_\tau\big|\tilde\mcF_\gamma\right].
    \end{equation}
  \item\label{Snell:DoobMeyer} The Doob-Meyer decomposition of the supermartingale $Z$ implies the existence of a triple $(M,K^c,K^d)$ where $(M_t:0\leq t\leq T)$ is a uniformly integrable right-continuous martingale, $(K^c_t:0\leq t\leq T)$ is a non-decreasing, predictable, continuous process with $K^c_0=0$ and $K^d_t$ is non-decreasing purely discontinuous predictable with $K^d_0=0$, such that
      \begin{equation}\label{ekv:DoobMeyerDec}
        Z_t=M_t-K^c_t-K^d_t.
      \end{equation}
      Furthermore, $\{\Delta_t K^d>0\}\subset \{\Delta_t U<0\}\cup\{Z_{t^-}=U_{t^-}\}$ for all $t\in[0,T]$.
  \item\label{Snell:att} Let $\eta\in\tilde\mcT_0$ and assume that $U$ is non-negative, $L^2$-bounded, right continuous and that $\lim_{m\to\infty}\E[U_{\gamma_m}]\leq\E[U_{\gamma}]$ whenever $\gamma_m\nearrow\gamma\in\tilde\mcT_\eta$. Then, the stopping time $\tau^*_{\eta}$ defined by $\tau^*_{\eta}:=\inf\{s\geq\eta:Z_s=U_s\}\wedge T$ is optimal after $\eta$, \ie
    \begin{equation*}
      Z_{\eta}=\E\left[U_{\tau^*_\eta}\big|\tilde\mcF_\eta\right].
    \end{equation*}
    Furthermore, in this setting the Snell envelope, $Z$, is left continuous in expectation, \ie $K_d\equiv 0$, and $Z$ is a martingale on $[\eta,\tau^*_\eta]$.
  \item\label{Snell:lim} Let $U^k$ be a sequence of \cadlag processes converging pointwisely to the \cadlag process $U$ and let $Z^k$ be the Snell envelope of $U^k$. Then the sequence $Z^k$ converges pointwisely to a process $Z$ and $Z$ is the Snell envelope of $U$.
  \item\label{Snell:DynP} We have the following dynamic programming relation: For any $\gamma\in\tilde\mcT_0$ and any $\eta\in\tilde\mcT_{\gamma}$ we have
  \begin{align*}
    Z_\gamma=\esssup_{\tau\in\tilde\mcT_{\gamma}}\E\left[\ett_{[\eta\leq\tau]}Z_\eta+\ett_{[\tau>\eta]}U_\tau \big|\tilde\mcF_\gamma\right].
  \end{align*}
\end{enumerate}
\end{thm}
In the above theorem (\ref{Snell:sup})-(\ref{Snell:att}) are standard. Proofs can be found in \cite{ElKarouiLN} (see \cite{Latifa} for an English version), Appendix D in~\cite{KarShreve2}, \cite{HamRefBSDE} and in the appendix of~\cite{CvitKar}. Statement (\ref{Snell:lim}) was proved in \cite{BollanMSwitch1}. The last statement follows by noting that
\begin{align*}
\esssup_{\tau\in\tilde\mcT_{\gamma}}\E\left[\ett_{[\eta\leq\tau]}Z_\eta+\ett_{[\tau>\eta]}U_\tau \big|\tilde\mcF_\gamma\right]&=\esssup_{\tau\in\tilde\mcT_{\gamma}}\E\left[\ett_{[\eta\leq\tau]}\esssup_{\tau'\in \tilde\mcT_{\eta}}\E\left[U_{\tau'}\big|\tilde\mcF_\eta\right]+\ett_{[\tau>\eta]}U_\tau \big|\tilde\mcF_\gamma\right]
\\
&=\esssup_{\tau\in\tilde\mcT_{\gamma}}\E\left[\ett_{[\eta\leq\tau]}\E\left[U_\tau\big|\tilde\mcF_\eta\right]+\ett_{[\tau>\eta]}U_\tau \big|\tilde\mcF_\gamma\right]
\\
&=\esssup_{\tau\in\tilde\mcT_{\gamma}}\E\left[\ett_{[\eta\leq\tau]}U_\tau+\ett_{[\tau>\eta]}U_\tau \big|\tilde\mcF_\gamma\right]=Z_\gamma.
\end{align*}

The Snell envelope will be the main tool in showing that Problem 1 has a unique solution.

%%%%%%%%%%%%%%%%%%%%%%%%%%%%%%%%%%%%%%%%%%%%%%%%%%%%%%%%%%%%%%%%%%%%%%%%%%%%%%%%%%%%%%%%%%%%%%%%%%%%%%%%%%%%%%%%%%%%%%%%%%%%%%%%%%%%%%%%%%%%%%%%%%
%%%%%%%%%%%%%%%%%%%%%%%%%%%%%%%%%%%%%%%%%%%%%%%%%%%%%%%%%%%%%%%%%%%%%%%%%%%%%%%%%%%%%%%%%%%%%%%%%%%%%%%%%%%%%%%%%%%%%%%%%%%%%%%%%%%%%%%%%%%%%%%%%%
%%%%%%%%%%%%%%%%%%%%%%%%%%%%%%%%%%%%%%%%%%%%%%%%%%%%%%%%%%%%%%%%%%%%%%%%%%%%%%%%%%%%%%%%%%%%%%%%%%%%%%%%%%%%%%%%%%%%%%%%%%%%%%%%%%%%%%%%%%%%%%%%%%

\section{A verification theorem\label{sec:VERthm}}
The method for solving Problem 1 will be based on deriving an optimal control under the assumption that a specific family of processes exists, and then (in  the next section) showing that the family indeed does exist. We will refer to any such family of processes as a verification family.
\begin{defn}\label{def:vFAM}
We define a \emph{verification family} to be a family of \cadlag supermartingales $((Y^{t,\nu,z,\veca,\vecb}_s)_{0\leq s\leq T}: (t,\nu,z,\veca,\vecb)\in \bigSET)$ such that:
\begin{enumerate}[a)]
  \item\label{vFAM:inS2e} For every $(t,\nu,z,\veca,\vecb)\in \bigSET$ we have $Y^{t,\nu,z,\veca,\vecb}\in\mcS_e^2$ and $(Y^{s,\nu,z,\veca,\vecb}_s:0\leq s\leq T)\in\mcS_e^2$.
  \item\label{vFAM:bnd} The family is bounded in the sense that $\E[\sup\limits_{(t,\nu,z,\veca,\vecb)\in\bigSET}\sup\limits_{s\in[0,T]}|Y^{t,\nu,z,\veca,\vecb}_s|^2]<\infty$.
  \item\label{vFAM:cont} The family is continuous in $(\nu,z)$ in the sense that
  \begin{equation*}
    \lim_{(p,q)\to(0,0)}\E\Big[\sup_{(t,\nu,z)\in \mcD_{(\veca,\vecb)}}|Y^{t,(\nu+p)^+\wedge \vecb T,(z+q)^+\wedge T\veca^+,\veca,\vecb}_{t}-Y^{t,\nu,z,\veca,\vecb}_t|^2\Big]= 0,
  \end{equation*}
  for every $(\veca,\vecb)\in\mcJ$.
  \item\label{vFAM:recur} The family satisfies the recursion
  \begin{align}\nonumber
  Y^{t,\nu,z,\veca,\vecb}_s=\esssup_{\tau \in \mcT_{s}} \E\bigg[&\int_s^{\tau\wedge T}\psi_{A^{t,\nu,\veca,\vecb}_r}\left(r,\theta_r^{t,\nu,z,\veca,\vecb}\right)dr+\ett_{[\tau \geq T]}\Upsilon_{A^{t,\nu,\veca,\vecb}_T}\left(\theta_T^{t,\nu,z,\veca,\vecb}\right)
  \\
  &+\ett_{[\tau < T]}\max_{\beta\in\mcI}\left\{-c_{\vecb,\beta}(\tau)+Y^{\tau,\beta\nu + \tau(\beta-\vecb)^+,\theta_\tau^{t,\nu,z,\veca,\vecb}\wedge T\beta,A^{t,\nu,\veca,\vecb}_\tau\wedge\beta,\beta}_\tau\right\}\Big| \mcG_s\bigg].\label{ekv:Ydef}
  \end{align}
\end{enumerate}
\end{defn}

The purpose of the present section is to reduce the solution of Problem 1 to showing existence of a verification family. This is done in the verification theorem below. First we give a lemma that will be used in the proof of the verification theorem:
\begin{lem}\label{lemma:contEXT}
Let $((Y^{t,\nu,z,\veca,\vecb}_s)_{0\leq s\leq T}: (t,\nu,z,\veca,\vecb)\in \bigSET)$ be a verification family, then
\begin{align*}
\lim_{(p,q)\to(0,0)}\E\Big[\sup_{(t,\nu,z)\in \mcD_{(\veca,\vecb)}}\sup_{s\in[0,T]}|\E[Y^{t,(\nu+p)^+\wedge \vecb T,(z+q)^+\wedge T\veca^+,\veca,\vecb}_{s}-Y^{t,\nu,z,\veca,\vecb}_s|\mcF_t]|^2\Big]= 0.
\end{align*}
\end{lem}

\noindent\emph{Proof.} For $s\in[0,t]$ the result follows immediately from property \emph{\ref{vFAM:cont})} and Doob's maximal inequality. We thus let $s\in(t,T]$ and note that if $\veca\in\mcAabs$ then
\begin{align*}
\E[Y^{t,(\nu+p)^+\wedge \vecb T,(z+q)^+\wedge T\veca^+,\veca,\vecb}_{s}-Y^{t,\nu,z,\veca,\vecb}_s|\mcF_t]=\E[Y^{s,(\nu+p)^+\wedge \vecb T,\theta_s^{t,(z+q)^+\wedge T\veca^+,\veca},\veca,\vecb}_{s}-Y^{s,\nu,\theta_s^{t, z,\veca},\veca,\vecb}_s|\mcF_t]
\end{align*}
and by Doob's maximal inequality we get, since $\theta_s^{t,(z+q)^+\wedge T\veca^+,\veca}-\theta_s^{t, z,\veca}=(z+q)^+\wedge T\veca^+-z$, that
\begin{align*}
&\E\Big[\sup_{(t,\nu,z)\in \mcD_{(\veca,\vecb)}}\sup_{s\in[0,T]}|\E[Y^{t,(\nu+p)^+\wedge \vecb T,(z+q)^+\wedge T\veca^+,\veca,\vecb}_{s}-Y^{t,\nu,z,\veca,\vecb}_s|\mcF_t]|^2\Big]
\\
&\leq C\E\Big[\sup_{(t,\nu,z)\in \mcD_{(\veca,\vecb)}}\sup_{s\in[0,T]}|Y^{s,(\nu+p)^+\wedge \vecb T,\theta_s^{t,(z+q)^+\wedge T\veca^+,\veca},\veca,\vecb}_{s}-Y^{s,\nu,\theta_s^{t, z,\veca},\veca,\vecb}_s|^2\Big]
\\
&\leq C\E\Big[\sup_{(s,\nu,z)\in \mcD_{(\veca,\vecb)}}\sup_{s\in[0,T]}|Y^{s,(\nu+p)^+\wedge \vecb T,(z+q)^+\wedge T\veca^+,\veca,\vecb}_{s}-Y^{s,\nu,z,\veca,\vecb}_s|^2\Big],
\end{align*}
which tends to 0 as $(p,q)\to(0,0)$ by property \emph{\ref{vFAM:cont})}. For $\veca\in\mcA_{\vecb}\setminus \mcAabs_\vecb$ we find, by arguing as in the proof of Lemma~\ref{lemma:condINTisL1cont}, that for $(t,\nu,z)$ and $(t,\nu',z')$ in $\mcD_{(\veca,\vecb)}$,
\begin{align*}
\E\Big[Y^{t,\nu',z',\veca,\vecb}_s-Y^{t,\nu,z,\veca,\vecb}_s|\mcF_t\Big]&=\E\Big[\ett_{[\eta'\leq s]}Y^{\eta',\nu',\theta_{\eta'}^{t,z',\veca}\wedge T A^{t,\nu',\veca,\vecb}_{\eta'},A^{t,\nu',\veca,\vecb}_{\eta'},\vecb}_{s}-\ett_{[\eta\leq s]}Y^{\eta,\nu,\theta_\eta^{t,z,\veca}\wedge TA^{t,\nu,\veca,\vecb}_{\eta},A^{t,\nu,\veca,\vecb}_\eta,\vecb}_s
\\
&\quad+\ett_{[\eta' > s]}Y^{s,\nu',\theta_s^{t,z',\veca},\veca,\vecb}_s-\ett_{[\eta > s]}Y^{s,\nu,\theta_s^{t,z,\veca},\veca,\vecb}_s\big| \mcF_t\Big]
\\
%%%%%%%%%%%%%%%%%%%%%%%%%%%%%%%%%%%%%%%%%%%%%%%%%%%%%%%%%%%%%%%%%%%%%%%%%%%%%%%%%%%%%%%%%%%%%%%%%%%%%%%%%%%%%%%%%%%%%%%%%%%%%%%%%%%%%%%%%%%%%%%%%%
&=\E\bigg[\sum_{\veca'\in\mcA^{-\veca}_{\veca,\vecb}}\Big(\int_{t}^{s}\lambda^{\nu',\vecb}_{\veca,\veca'}(r)e^{\int_{t}^r \lambda^{\nu',\vecb}_{\veca,\veca}(v)dv} Y^{r,\nu',\theta_{r}^{t,z',\veca}\wedge T\veca',\veca',\vecb}_sdr
\\
&\quad-\int_{t}^{s}\lambda^{\nu,\vecb}_{\veca,\veca'}(r)e^{\int_{t}^r\lambda^{\nu,\vecb}_{\veca,\veca}(v)dv} Y^{r,\nu,\theta_{r}^{t,z,\veca}\wedge T\veca',\veca',\vecb}_s dr\Big)
\\
&\quad + e^{\int_{t}^s\lambda^{\nu',\vecb}_{\veca,\veca}(v)dv} Y^{s,\nu',\theta_{s}^{t,z',\veca},\veca,\vecb}_s - e^{\int_{t}^s\lambda^{\nu,\vecb}_{\veca,\veca}(v)dv} Y^{s,\nu,\theta_{s}^{t,z,\veca},\veca,\vecb}_s\Big| \mcF_{t}\bigg],
\end{align*}
where $\eta$ and $\eta'$ are the first transition times of $A^{t,\nu,\veca,\vecb}$ and $A^{t,\nu',\veca,\vecb}$, respectively. Using the relation $ab-a'b'=1/2((a-a')(b+b')+(a+a')(b-b'))$ we get that
\begin{align*}
\E\Big[Y^{t,\nu',z',\veca,\vecb}_s-Y^{t,\nu,z,\veca,\vecb}_s|\mcF_t\Big]
%%%%%%%%%%%%%%%%%%%%%%%%%%%%%%%%%%%%%%%%%%%%%%%%%%%%%%%%%%%%%%%%%%%%%%%%%%%%%%%%%%%%%%%%%%%%%%%%%%%%%%%%%%%%%%%%%%%%%%%%%%%%%%%%%%%%%%%%%%%%%%%%%%
&\leq\E\bigg[\sum_{\veca'\in\mcA^{-\veca}_{\veca,\vecb}}\Big(\int_{t}^{s}\big(|\lambda^{\nu',\vecb}_{\veca,\veca'}(r)e^{\int_{t}^r \lambda^{\nu',\vecb}_{\veca,\veca}(v)dv}- \lambda^{\nu,\vecb}_{\veca,\veca'}(r)e^{\int_{t}^r\lambda^{\nu,\vecb}_{\veca,\veca}(v)dv}|
\\
&\quad \cdot |Y^{r,\nu',\theta_{r}^{t,z',\veca}\wedge T\veca',\veca',\vecb}_s+Y^{r,\nu,\theta_{r}^{t,z,\veca}\wedge T\veca',\veca',\vecb}_s|
\\
&\quad+2K_\lambda |Y^{r,\nu',\theta_{r}^{t,z',\veca}\wedge T\veca',\veca',\vecb}_s-Y^{r,\nu,\theta_{r}^{t,z,\veca}\wedge T\veca',\veca',\vecb}_s|\big)dr\Big)
\\
&\quad + |e^{\int_{t}^s\lambda^{\nu',\vecb}_{\veca,\veca}(v)dv}-e^{\int_{t}^s\lambda^{\nu,\vecb}_{\veca,\veca}(v)dv}| |Y^{s,\nu',\theta_{s}^{t,z',\veca},\veca,\vecb}_s + Y^{s,\nu,\theta_{r}^{t,z,\veca},\veca,\vecb}_s|
\\
&\quad + 2 |Y^{s,\nu',\theta_{s}^{t,z',\veca},\veca,\vecb}_s - Y^{s,\nu,\theta_{r}^{t,z,\veca},\veca,\vecb}_s|\Big| \mcF_{t}\bigg]
\\
%%%%%%%%%%%%%%%%%%%%%%%%%%%%%%%%%%%%%%%%%%%%%%%%%%%%%%%%%%%%%%%%%%%%%%%%%%%%%%%%%%%%%%%%%%%%%%%%%%%%%%%%%%%%%%%%%%%%%%%%%%%%%%%%%%%%%%%%%%%%%%%%%%
&\leq C\E\bigg[|\nu'-\nu|\sup_{(t,\nu,z,\veca,\vecb)\in\bigSET} \sup_{s\in[0,T]}|Y^{t,\nu,z,\veca,\vecb}_s|
\\
&\quad+ \sum_{\veca'\in\mcA^{-\veca}_{\veca,\vecb}}\Big\{\sup_{r\in[0,T]}(|Y^{r,\nu',\theta_{r}^{t,z',\veca}\wedge T\veca',\veca',\vecb}_s - Y^{r,\nu,\theta_{r}^{t,z,\veca}\wedge T\veca',\veca',\vecb}_s|)\Big\}
\\
&\quad + |Y^{s,\nu',\theta_{s}^{t,z',\veca},\veca,\vecb}_s-Y^{s,\nu,\theta_{r}^{t,z,\veca},\veca,\vecb}_s| \Big| \mcF_{t}\bigg]
%%%%%%%%%%%%%%%%%%%%%%%%%%%%%%%%%%%%%%%%%%%%%%%%%%%%%%%%%%%%%%%%%%%%%%%%%%%%%%%%%%%%%%%%%%%%%%%%%%%%%%%%%%%%%%%%%%%%%%%%%%%%%%%%%%%%%%%%%%%%%%%%%%
\end{align*}
By again applying Doob's maximal inequality and an induction argument in the $\veca$-component, combined with the square integrability property in \emph{\ref{vFAM:bnd})}, the assertion follows.\qed\\

We have the following verification theorem:
\begin{thm}\label{thm:vfc}
Assume that there exists a verification family $((Y^{t,\nu,z,\veca,\vecb}_s)_{0\leq s\leq T}: (t,\nu,z,\veca,\vecb)\in \bigSET)$. Then $((Y^{t,\nu,z,\veca,\vecb}_s)_{0\leq s\leq T}: (t,\nu,z,\veca,\vecb)\in \bigSET)$ is unique (\ie there is at most one verification family, up to indistinguishability) and:
\begin{enumerate}[(i)]
  \item Satisfies $Y_0^{0,0,0,0,0}=\sup_{u\in \mcU} J(u)$.
  \item Defines the optimal strategy, $u^*=(\tau_1^*,\ldots,\tau_{N^*}^*;\beta_1^*,\ldots,\beta_{N^*}^*)$, for Problem 1, where $(\tau_j^*)_{1\leq j\leq {N^*}}$ is a sequence of $\bbG$-stopping times given by
  \begin{align}\nonumber
  \tau^*_j:=\inf \Big\{&s \geq \tau^*_{j-1}:\:Y_s^{\tau^*_{j-1},\vartheta^*_{j-1},z^*_{j-1},\veca^*_{j-1},\beta^*_{j-1}}=\max_{\beta\in \mcI}\Big\{-c_{\beta^*_{j-1},\beta}(s)
  \\
  &+Y^{s,\beta\vartheta^*_{j-1}+s(\beta-\beta^*_{j-1})^+,\theta_s^{\tau_{j-1}^*,\vartheta^*_{j-1},z_{j-1}^*,\veca^*_{j-1},\beta_{j-1}^*}\wedge T\beta,A^{\tau^*_{j-1},\vartheta^*_{j-1},\veca^*_{j-1},\beta^*_{j-1}}_{s}\wedge\beta,\beta}_s\Big\}\Big\},\label{ekv:taujDEF}
  \end{align}
  $(\beta_j^*)_{1\leq j\leq {N^*}}$ is defined as a measurable selection of
  \begin{equation*}
  \beta^*_j\in\mathop{\arg\max}_{\beta\in \mcI}\Big\{-c_{\beta^*_{j-1},\beta}(\tau^*_j) + Y^{\tau^*_j,\beta\vartheta^*_{j-1}+\tau^*_j(\beta-\beta^*_{j-1})^+, \theta_{\tau^*_j}^{\tau_{j-1}^*,\vartheta^*_{j-1},z_{j-1}^*,\veca^*_{j-1},\beta_{j-1}^*}\wedge T\beta,A^{\tau^*_{j-1},\vartheta^*_{j-1},\veca^*_{j-1},\beta^*_{j-1}}_{\tau^*_j}\wedge\beta,\beta}_{\tau^*_j}\Big\},
  \end{equation*}
  where $\vartheta^*_{j}=\beta^*_j\vartheta^*_{j-1}+\tau^*_j(\beta^*_j-\beta^*_{j-1})^+$, $z^*_{j}:=\theta_{\tau_{j}^*}^{\tau^*_{j-1},\vartheta^*_{j-1},z_{j-1}^*,\veca^*_{j-1},\beta^*_{j-1}}\wedge T\beta^*_{j}$ and $\veca^*_{j}:=A^{\tau^*_{j-1},\vartheta^*_{j-1},\veca^*_{j-1},\beta_{j-1}^*}_{\tau^*_{j}}\wedge\beta_j^*$, with $(\tau^*_0,\vartheta^*_0,z^*_0,\veca^*_0,\beta^*_{0}):=(0,0,0,0,0)$ and $N^*:=\max\{j:\tau_j^*<T\}$.
%  \begin{equation}\label{ekv:tau_ijDEF}
%  \tau^i_j:=\min_{k\geq 1} \{\tau_k > \tau^i_{j-1}:\: (\vecb^*_{k})_i\neq (\vecb^*_{k-1})_i\}.
%  \end{equation}
\end{enumerate}
\end{thm}

\noindent\emph{Proof.} Note that the proof amounts to showing that for all $(t,\nu,z,\veca,\vecb)\in \bigSET$, we have
\begin{align*}
Y^{t,\nu,z,\veca,\vecb}_s=\esssup_{u\in \mcU_s^{f}} \E\bigg[&\int_s^{T}\psi_{\alpha_r^{t,\nu,\veca,\vecb,u}}(r,\theta^{t,\nu,z,\veca,\vecb,u}_r)dr
+\Upsilon_{\alpha_T^{t,\nu,\veca,\vecb,u}}(\theta^{t,\nu,z,\veca,\vecb,u}_T)
-\sum_{j=1}^N c_{\beta_{j-1},\beta_{j}}(\tau_j)\Big|\mcG_s\bigg],
\end{align*}
for all $s\in [t,T]$, where $\beta_0=\vecb$. Then uniqueness is immediate, $(i)$ follows from Proposition \ref{prop:finSTRAT} and $(ii)$ follows from repeated use of Theorem~\ref{thm:Snell}.\ref{Snell:att}.\\

\noindent {\bf Step 1} We start by showing that for all $(t,\nu,z,\veca,\vecb)\in\bigSET$ the recursion \eqref{ekv:Ydef} can be written in terms of stopping times. From \eqref{ekv:Ydef} we have that, for each $(t,\nu,z,\veca,\vecb)\in \bigSET$,
\begin{align*}
Z^{t,\nu,z,\veca,\vecb}:=\Big(Y_s^{t,\nu,z,\veca,\vecb}+\int_0^s\psi_{A^{t,\nu,\veca,\vecb}_r}(r,\theta_r^{t,\nu,z,\veca,\vecb})dr:\: 0\leq s\leq T\Big)
\end{align*}
is the smallest supermartingale that dominates the process
\begin{align*}
&\bigg(\int_0^{s}\psi_{A^{t,\nu,\veca,\vecb}_r}(r,\theta_r^{t,\nu,z,\veca,\vecb})dr+\ett_{[s = T]}\Upsilon_{A^{t,\nu,\veca,\vecb}_T}(\theta_T^{t,\nu,z,\veca,\vecb})
\\
&+\ett_{[s < T]}\max_{\beta\in\mcI}\left\{-c_{\vecb,\beta}(s)+Y^{s,\beta\nu+s(\beta-\vecb)^+,\theta_s^{t,\nu,z,\veca,\vecb}\wedge T\beta,A^{t,\nu,\veca,\vecb}_s\wedge\beta,\beta}_s\right\}:0\leq s\leq T\bigg).
\end{align*}
We will show that under assumptions \emph{\ref{vFAM:inS2e})}-\emph{\ref{vFAM:cont})} the dominated process satisfies the assumptions of Theorem~\ref{thm:Snell}.\ref{Snell:att} after which the assertion follows. Fix $\vecb'\in\mcI^{-\vecb}$ and note that for all $s\leq s'\leq T$ we have the following trivial relation
\begin{align*}
&Y^{s,\vecb'\nu+s(\vecb'-\vecb)^+,\theta_s^{t,\nu,z,\veca,\vecb}\wedge T\vecb',A^{t,\nu,\veca,\vecb}_{s}\wedge\vecb',\vecb'}_s - Y^{s',\vecb'\nu+s'(\vecb'-\vecb)^+,\theta_{s'}^{t,\nu,z,\veca,\vecb}\wedge T\vecb',A^{t,\nu,\veca,\vecb}_{s'}\wedge\vecb',\vecb'}_{s'}
\\
&= (Y^{s,\vecb'\nu+s(\vecb'-\vecb)^+,\theta_s^{t,\nu,z,\veca,\vecb}\wedge T\vecb',A^{t,\nu,\veca,\vecb}_{s}\wedge\vecb',\vecb'}_s - Y^{s',\vecb'\nu+s'(\vecb'-\vecb)^+,\theta_{s'}^{t,\nu,z,\veca,\vecb}\wedge T\vecb',A^{t,\nu,\veca,\vecb}_{s}\wedge\vecb',\vecb'}_{s'})
\\
&\quad + (Y^{s',\vecb'\nu+s'(\vecb'-\vecb)^+,\theta_{s'}^{t,\nu,z,\veca,\vecb}\wedge T\vecb',A^{t,\nu,\veca,\vecb}_{s}\wedge\vecb',\vecb'}_{s'} - Y^{s',\vecb'\nu+s'(\vecb'-\vecb)^+,\theta_{s'}^{t,\nu,z,\veca,\vecb}\wedge T\vecb',A^{t,\nu,\veca,\vecb}_{s'}\wedge\vecb',\vecb'}_{s'}).
\end{align*}
If $\gamma_m$ is a sequence of stopping times such that $\gamma_m\nearrow \gamma\in\mcT$, $\Prob$-a.s., we thus have
\begin{align*}
\lim_{m\to\infty}\E[Y^{{\gamma_m},\vecb'\nu+(\vecb'-\vecb)^+{\gamma_m},\theta_{{\gamma_m}}^{t,\nu,z,\veca,\vecb}\wedge T\vecb',A^{t,\nu,\veca,\vecb}_{{\gamma_m}}\wedge\vecb',\vecb'}_{\gamma_m}-Y^{{\gamma},\vecb'\nu+(\vecb'-\vecb)^+ \gamma,\theta_{\gamma}^{t,\nu,z,\veca,\vecb}\wedge T\vecb',A^{t,\nu,\veca,\vecb}_{\gamma}\wedge\vecb',\vecb'}_{\gamma}]
\\
=\lim_{m\to\infty}\E[Y^{{\gamma},\vecb'\nu+(\vecb'-\vecb)^+{\gamma},\theta_{{\gamma}}^{t,\nu,z,\veca,\vecb}\wedge T\vecb',A^{t,\nu,\veca,\vecb}_{{\gamma_m}}\wedge\vecb',\vecb'}_{\gamma}-Y^{{\gamma},\vecb'\nu+(\vecb'-\vecb)^ +\gamma,\theta_{\gamma}^{t,\nu,z,\veca,\vecb}\wedge T\vecb',A^{t,\nu,\veca,\vecb}_{\gamma}\wedge\vecb',\vecb'}_{\gamma}]=0,
\end{align*}
where the first equality follows by \emph{\ref{vFAM:inS2e})} and \emph{\ref{vFAM:cont})} and the second equality follows from the $L^2$-boundedness assumed in \emph{\ref{vFAM:bnd})} in combination with Lemma~\ref{lemma:AisLCE}.\\

The dominated process is thus $L^2$-bounded by Assumption~\ref{ass:prelim} and \emph{\ref{vFAM:bnd})}, positive by Remark~\ref{rem:pos} and LCE on $[0,T)$. At time $T$ it may have a jump but the jump has to be positive by Assumption~\ref{ass:prelim}.\ref{ass:@end}. Theorem~\ref{thm:Snell}.\ref{Snell:att} now implies that, for each $\gamma\in\mcT$, there is a stopping time, $\tau_{\gamma}\in \mcT_\gamma$, such that:
\begin{align*}
Y^{t,\nu,z,\veca,\vecb}_\gamma=\E\bigg[&\int_\gamma^{\tau_\gamma\wedge T}\psi_{A^{t,\nu,\veca,\vecb}_r}\left(r,\theta_r^{t,\nu,z,\veca,\vecb}\right)dr +\ett_{[\tau_\gamma \geq T]}\Upsilon_{A^{t,\nu,\veca,\vecb}_T}\left(\theta_T^{t,\nu,z,\veca,\vecb}\right)
\\
&+\ett_{[\tau_\gamma < T]}\max_{\beta\in\mcI}\left\{-c_{\vecb,\beta}(\tau_\gamma)+Y^{\tau_\gamma,\beta\nu+\tau_\gamma(\beta-\vecb)^+, \theta_{\tau_\gamma}^{t,\nu,z,\veca,\vecb}\wedge T\beta,A^{t,\nu,\veca,\vecb}_{\tau_\gamma}\wedge\beta,\beta}_{\tau_\gamma}\right\}\Big| \mcG_\gamma\bigg].
\end{align*}

\noindent {\bf Step 2} We now show that if $u^*\in\mcU^f$ then $Y^{0,0,0,0,0}_0=J(u^*)$. First, define
\begin{align*}
Z_s:=Y_s^{0,0,0,0,0}+\int_0^s\psi_0(r,0)dr.
\end{align*}
Then by Theorem~\ref{thm:Snell}, $Z_s$ is the smallest supermartingale that dominates
\begin{align*}
\Big(&\int_0^{s}\psi_0\left(r,0\right)dr+\ett_{[s=T]}\Upsilon_0\left(0\right)
+ \ett_{[s < T]}\max_{\beta\in \mcI^{-0}}\left\{-c_{0,\beta}+Y^{s,s\beta,0,0,\beta}_s\right\}:0\leq s\leq T\Big)
\end{align*}
and by step 1
\begin{align*}
Y_0^{0,0,0,0,0}&=\esssup_{\tau \in \mcT} \E\bigg[\int_0^{\tau\wedge T}\psi_0\left(r,0\right)dr+\ett_{[\tau\geq T]}\Upsilon_0\left(0\right)\ett_{[\tau < T]}\max_{\beta\in \mcI^{-0}}\left\{-c_{0,\beta}(\tau)+Y^{\tau,\tau\beta,0,0,\beta}_\tau\right\}\bigg]
\\
&=\E\bigg[\int_0^{\tau^*_1\wedge T}\psi_0\left(r,0\right)dr+\ett_{[\tau^*_1\geq T]}\Upsilon_0\left(0\right)+ \ett_{[\tau^*_1< T]}\max_{\beta\in \mcI^{-0}}\left\{-c_{0,\beta}(\tau^*_1)+Y^{\tau^*_1,\tau^*_1\beta,0,0,\beta}_{\tau^*_1}\right\}\bigg]
\\
&=\E\bigg[\int_0^{\tau^*_1\wedge T}\psi_0\left(r,0\right)dr+\ett_{[\tau^*_1\geq T]}\Upsilon_0\left(0\right)+ \ett_{[\tau^*_1 < T]}\left\{-c_{0,\beta^*_1}(\tau^*_1)+Y^{\tau^*_1,\nu_1^*,z_1^*,\veca^*_1,\beta^*_1}_{\tau^*_1}\right\}\bigg].
\end{align*}
Now suppose that, for some $j'>0$ we have, for all $j\leq j'$,
\begin{align*}
Y&_s^{\tau^*_{j-1},\vartheta^*_{j-1},z^*_{j-1},\veca^*_{j-1}\beta^*_{j-1}}=\E\bigg[\int_s^{\tau^*_j\wedge T}\psi_{A^{\tau^*_{j-1},\vartheta^*_{j-1},\veca^*_{j-1}\beta^*_{j-1}}_r}\left(r,\theta_r^{\tau^*_{j-1}, \vartheta^*_{j-1},z^*_{j-1},\veca^*_{j-1}\beta^*_{j-1}}\right)dr
\\
&+\ett_{[\tau^*_j \geq T]}\Upsilon_{A^{\tau^*_{j-1},\vartheta^*_{j-1},\veca^*_{j-1}\beta^*_{j-1}}_T} \left(\theta_T^{\tau^*_{j-1},\vartheta^*_{j-1},z^*_{j-1},\veca^*_{j-1}\beta^*_{j-1}}\right)
+ \ett_{[\tau^*_j < T]} \left\{-c_{\beta^*_{j-1},\beta^*_{j}}(\tau^*_j)+Y^{\tau^*_j,\vartheta^*_{j},z^*_{j},\veca^*_j,\beta^*_{j}}_{\tau^*_j}\right\}\Big| \mcG_s\bigg],
\end{align*}
$\Prob$-a.s., for each $\tau^*_{j-1}\leq s\leq T$.\\

For all $M\geq 1$, let $(G_k^M)_{1\leq k \leq M}$ be an partition of $[0,T]^n$ into sets of the type $G_k^M=[z_{M,k,1},z_{M,k,1}')\times\ldots\times[z_{M,k,n},z_{M,k,n}')$ where $z_{M,k,i}< z_{M,k,i}'$ and $\max_{k,i}|z_{M,k,i}-z_{M,k,i}'|\to 0$ as $M\to\infty$ and let $(\kappa_{k}^M)_{1\leq k\leq M}$ be the sequence of points in $\R^n$ given by $[\kappa_k^M]_i:=z_{M,k,i}$. For $M\geq 1$ and $t\geq \tau^*_{j'}$, we define the process
\begin{equation*}
\hat Y^{t,M}_{s}:=\sum_{(\veca,\vecb)\in \mcJ}\ett_{[A_{t}^{\tau^*_{j'},\vartheta^*_{j'},\veca^*_{j'},\beta^*_{j'}}=\veca]}\ett_{[\beta^*_{j'}=\vecb]}\sum_{k=1}^M\ett_{[\vartheta^*_{j'} \in G_k^M]}\sum_{l=1}^M\ett_{[\theta_t^{\tau^*_{j'},\vartheta^*_{j'},z^*_{j'},\veca^*_{j'},\beta^*_{j'}}\in G_l^M]} Y_s^{t,\kappa_{k}^M,\kappa_{l}^M,\veca,\vecb},
\end{equation*}
for all $s\in[t,T]$. Now, $\ett_{[A_{t}^{\tau^*_{j'},\vartheta^*_{j'},\veca^*_{j'},\beta^*_{j'}}=\veca]}\ett_{[\beta^*_{j'}=\vecb]}\ett_{[\vartheta^*_{j'}\in G_k^M]}\ett_{[z^*_{j'}\in G_l^M]}\ett_{[\theta_t^{\tau^*_{j'},\vartheta^*_{j'},z^*_{j'},\veca^*_{j'},\beta^*_{j'}}\in G_l^M]}$\\
$\cdot\Big(Y_s^{t,\kappa_{k}^M,\kappa_{l}^M,\veca,\vecb} + \int_{t}^s\psi_{A^{t,\kappa_{k}^M,\veca,\vecb}_r}(r,\theta_r^{t,\kappa_{k}^M,\kappa_{l}^M,\veca,\vecb})dr\Big)$
is the product of a $\mcG_{t}$--measurable positive r.v.~and a supermartingale, thus, it is a supermartingale for $s\geq t$. Hence, as
\begin{align*}
\Big(&\hat Y^{t,M}_s+\sum_{(\veca,\vecb)\in \mcJ}\ett_{[A_{t}^{\tau^*_{j'},\vartheta^*_{j'},\veca^*_{j'},\beta^*_{j'}}=\veca]}\ett_{[\beta^*_{j'}=\vecb]}\sum_{k=1}^M\ett_{[\vartheta^*_{j'}\in G_k^M]}\sum_{l=1}^M\ett_{[\theta_t^{\tau^*_{j'},\vartheta^*_{j'},z^*_{j'},\veca^*_{j'},\beta^*_{j'}}\in G_l^M]}
\\
&\cdot\int_{t}^s\psi_{A^{t,\kappa_{k}^M,\veca,\vecb}_r}(r,\theta_r^{t,\kappa_{k}^M,\kappa_{l}^M,\veca,\vecb})dr: t\leq s\leq T\Big)
\end{align*}
is the sum of a finite number of supermartingales it is also a supermartingale.

%%%%%%%%%%%%%%%%%%%%%%%%%%%%%%%%%% Ny del %%%%%%%%%%%%%%%%%%%%%%%%%%%%%%%%%%%
By Lemma~\ref{lemma:contEXT} we have that
\begin{align*}
\E\Big[\sup_{t\in[\tau^*_{j'},T]}\sup_{s\in[t,T]}| \E\big[Y_s^{\tau^*_{j'},\vartheta^*_{j'},z^*_{j'},\veca^*_{j'},\beta^*_{j'}} - Y^{t,M}_s|\mcG_t\big]\big|^2\Big]\to 0,
\end{align*}
as $M\to\infty$. This implies that there is a subsequence $(M_\iota)_{\iota\geq 1}$ such that
\begin{align*}
\sup_{t\in[\tau^*_{j'},T]}\sup_{s\in[t,T]}\E\big[Y_s^{\tau^*_{j'},\vartheta^*_{j'},z^*_{j'},\veca^*_{j'},\beta^*_{j'}} - Y^{t,M_\iota}_s\big|\mcG_t\big]\to 0,
\end{align*}
$\Prob$-a.s.~as $\iota\to\infty$. In particular, we note that $(Y^{t,M_\iota}_t:\tau^*_{j'}\leq t\leq T)$ is a sequence of \cadlag processes that converges $\Prob$-a.s.~to $(Y_t^{\tau^*_{j'},\vartheta^*_{j'},z^*_{j'},\veca^*_{j'},\beta^*_{j'}}:\tau^*_{j'}\leq t\leq T)$ uniformly in $t$ and we conclude that $(Y_t^{\tau^*_{j'},\vartheta^*_{j'},z^*_{j'},\veca^*_{j'},\beta^*_{j'}}:\tau^*_{j'}\leq t\leq T)$ is a \cadlag process.
%%%%%%%%%%%%%%%%%%%%%%%%%%%%%%%%%%%%%%%%%%%%%%%%%%%%%%%%%%%%%%%%%%%%%%%%%%%%%
Furthermore, by Lemma \ref{lemma:condINTisL1cont} and dominated convergence we get
\begin{align*}
&\E\Big[Y_s^{\tau^*_{j'},\vartheta^*_{j'},z^*_{j'},\veca^*_{j'},\beta^*_{j'}} + \int_{t}^s\psi_{A^{\tau^*_{j'},\vartheta^*_{j'},\veca^*_{j'},\beta^*_{j'}}_r} (r,\theta_r^{\tau^*_{j'},\vartheta^*_{j'},z^*_{j'},\veca^*_{j'},\beta^*_{j'}})dr\big|\mcG_t\Big]
\\
&=\mathop{\lim}_{\iota\to\infty}\E\Big[\hat Y^{t,{M_\iota}}_s+\sum_{(\veca,\vecb)\in \mcJ}\ett_{[A_{t}^{\tau^*_{j'},\vartheta^*_{j'},\veca^*_{j'},\beta^*_{j'}}=\veca]}\ett_{[\beta^*_{j'}=\vecb]}\sum_{k=1}^{M_\iota} \ett_{[\vartheta^*_{j'}\in G_k^{M_\iota}]}\sum_{l=1}^{M_\iota}\ett_{[\theta_t^{\tau^*_{j'},\vartheta^*_{j'},z^*_{j'},\veca^*_{j'},\beta^*_{j'}}\in G_l^{M_\iota}]}
\\
&\quad\cdot\int_{t}^s\psi_{A^{t,\kappa_{k}^{M_\iota},\veca,\vecb}_r}(r,\theta_r^{t,\kappa_{k}^{M_\iota},\kappa_{l}^{M_\iota},\veca,\vecb}) dr\big|\mcG_t\Big],
\end{align*}
$\Prob$-a.s., for all $s\in[t,T]$. This implies that for all $\tau^*_{j'}\leq t\leq s $ we have
\begin{align*}
Y_t^{\tau^*_{j'},\vartheta^*_{j'},z^*_{j'},\veca^*_{j'},\beta^*_{j'}}&=\mathop{\lim}_{\iota\to\infty}\hat Y^{t,{M_\iota}}_t
\\
&\geq \mathop{\lim}_{k\to\infty}\E\Big[\hat Y^{t,{{M_\iota}}}_s+\sum_{(\veca,\vecb)\in \mcJ}\ett_{[A_{t}^{\tau^*_{j'},\vartheta^*_{j'},\veca^*_{j'},\beta^*_{j'}}=\veca]}\ett_{[\beta^*_{j'}=\vecb]}\sum_{k=1}^{{M_\iota}} \ett_{[\vartheta^*_{j'} \in G_k^{M_\iota}]}\sum_{l=1}^{M_\iota} \ett_{[\theta_t^{\tau^*_{j'},\vartheta^*_{j'},z^*_{j'},\veca^*_{j'},\beta^*_{j'}}\in G_l^{M_\iota}]}
\\
&\quad\cdot\int_{t}^s\psi_{A^{t,\kappa_{k}^{M_\iota},\veca,\vecb}_r}(r,\theta_r^{t,\kappa_{k}^{M_\iota},\kappa_{l}^{M_\iota},\veca,\vecb})dr \big|\mcG_t\Big]
\\
&=E\Big[Y_s^{\tau^*_{j'},\vartheta^*_{j'},z^*_{j'},\veca^*_{j'},\beta^*_{j'}} + \int_{t}^s\psi_{A^{\tau^*_{j'},\vartheta^*_{j'},\veca^*_{j'},\beta^*_{j'}}_r}(r,\theta_r^{\tau^*_{j'},\vartheta^*_{j'},z^*_{j'},\veca^*_{j'}, \beta^*_{j'}})dr\big|\mcG_t\Big]
\end{align*}
$\Prob$-a.s.~where we have used the supermartingale property to reach the inequality. Hence, $\Big(Y_s^{\tau^*_{j'},\vartheta^*_{j'},z^*_{j'},\veca^*_{j'},\beta^*_{j'}} + \int_{\tau^*_{j'}}^s\psi_{A^{\tau^*_{j'},\veca^*_{j'},\veca^*_{j'},\beta^*_{j'}}_r}(r,\theta_r^{\tau^*_{j'},\vartheta^*_{j'},z^*_{j'}, \veca^*_{j'}, \beta^*_{j'}})dr:\: \tau^*_{j'}\leq s\leq T\Big)$ is a \cadlag supermartingale that dominates
\begin{align}\nonumber
&\Big(\int_{\tau^*_{j'}}^s\psi_{A^{\tau^*_{j'},\vartheta^*_{j'},\veca^*_{j'},\beta^*_{j'}}_r}(r,\theta_r^{\tau^*_{j'},\vartheta^*_{j'},z^*_{j'}, \veca^*_{j'},\beta^*_{j'}})dr+\ett_{[s = T]}\Upsilon_{A^{\tau^*_{j'},\vartheta^*_{j'},\veca^*_{j'},\beta^*_{j'}}_T}(\theta_T^{\tau^*_{j'},\vartheta^*_{j'},z^*_{j'},\veca^*_{j'}, \beta^*_{j'}})
\\
&+\ett_{[s < T]}\max_{\beta\in\mcI}\Big\{-c_{\beta^*_{j'},\beta}(s)+Y^{s,\vartheta^*_{j'}\beta+s(\beta-\beta^*_{j'})^+,\theta_s^{\tau^*_{j'}, \vartheta^*_{j'},z^*_{j'},\veca^*_{j'},\beta^*_{j'}}\wedge T\beta,A^{\tau^*_{j'},\vartheta^*_{j'},\veca^*_{j'},\beta^*_{j'}}_s\wedge \beta,\beta}_s\Big\}:\: \tau^*_{j'}\leq s\leq T\Big).\label{ekv:dominerad}
\end{align}
It remains to show that it is the smallest supermartingale with this property. Let $(Z_s:\:0\leq s\leq T)$ be a supermartingale that dominates \eqref{ekv:dominerad} for all $s\in[\tau^*_{j'},T]$. Then for each $(t,\nu,z,\veca,\vecb)\in \bigSET$ and $s\geq t$, we have
\begin{align*}
&\ett_{[\tau^*_{j'}=t]}\ett_{[\vartheta^*_{j'}=\nu]}\ett_{[z^*_{j'}=z]}\ett_{[\veca^*_{j'}=\veca]}\ett_{[\beta^*_{j'}=\vecb]}Z_s
\\
&\geq \ett_{[\tau^*_{j'}=t]}\ett_{[\vartheta^*_{j'}=\nu]}\ett_{[z^*_{j'}=z]}\ett_{[\veca^*_{j'}=\veca]}\ett_{[\beta^*_{j'}=\vecb]} \Big(\int_{t}^s\psi_{A^{t,\nu,\veca,\vecb}_r}(r,\theta_r^{t,\nu,z,\veca,\vecb})dr+\ett_{[s = T]}\Upsilon_{A^{t,\nu,\veca,\vecb}_T}(\theta_T^{t,\nu,z,\veca,\vecb})
\\
&\quad+\ett_{[s < T]}\max_{\beta\in\mcI}\Big\{-c_{\beta^*_{j'},\beta}(s)+Y^{s,\nu\beta+s(\beta-\vecb)^+,\theta_s^{t,\nu,z,\veca,\vecb}\wedge T\beta,A^{t,\nu,\veca,\vecb}_s\wedge\beta,\beta}_s\Big\}\Big),
\end{align*}
which by \eqref{ekv:Ydef} gives that
\begin{align*}
&\ett_{[\tau^*_{j'}=t]}\ett_{[\vartheta^*_{j'}=\nu]}\ett_{[z^*_{j'}=z]}\ett_{[\veca^*_{j'}=\veca]}\ett_{[\beta^*_{j'}=\vecb]}Z_s
\\
&\geq \ett_{[\tau^*_{j'}=t]}\ett_{[\vartheta^*_{j'}=\nu]}\ett_{[z^*_{j'}=z]}\ett_{[\veca^*_{j'}=\veca]}\ett_{[\beta^*_{j'}=\vecb]} \Big(Y^{t,\nu,z,\veca,\vecb}_s+\int_{t}^s\psi_{A^{t,\nu,\veca,\vecb}_r}(r,\theta_r^{t,\nu,z,\veca,\vecb})dr\Big\}\Big).
\end{align*}
Since this holds for all $(t,\nu,z,\veca,\vecb)\in \bigSET$ we get
\begin{equation*}
Z_s\geq Y^{\tau^*_{j'},\vartheta^*_{j'},z^*_{j'},\veca^*_{j'},\beta^*_{j'}}_s+\int_{\tau^*_{j'}}^s \psi_{A^{\tau^*_{j'},\vartheta^*_{j'},\veca^*_{j'},\beta^*_{j'}}_r}(r,\theta_r^{\tau^*_{j'},\vartheta^*_{j'},z^*_{j'},\veca^*_{j'},\beta^*_{j'}})dr
\end{equation*}
for all $s\geq \tau^*_{j'}$. Hence, $\Big(Y_s^{\tau^*_{j'},\vartheta^*_{j'},z^*_{j'},\veca^*_{j'},\beta^*_{j'}} + \int_{\tau^*_{j'}}^s\psi_{A^{\tau^*_{j'},\vartheta^*_{j'},\veca^*_{j'},\beta^*_{j'}}_r} (r,\theta_r^{\tau^*_{j'},\vartheta^*_{j'},z^*_{j'},\veca^*_{j'},\beta^*_{j'}})dr:\: \tau^*_{j'}\leq s\leq T\Big)$ is the Snell envelope of \eqref{ekv:dominerad} and by Theorem~\ref{thm:Snell}.\ref{Snell:att} and step 1
\begin{align*}
Y_s^{\tau^*_{j'},\vartheta^*_{j'},z^*_{j'},\veca^*_{j'},\beta^*_{j'}}=\E\bigg[&\int_s^{\tau^*_{j'+1}\wedge T}\psi_{A^{\tau^*_{j'},\vartheta^*_{j'},\veca^*_{j'},\beta^*_{j'}}_r} (r,\theta_r^{\tau^*_{j'},\vartheta^*_{j'},z^*_{j'},\veca^*_{j'},\beta^*_{j'}})dr+\ett_{[\tau^*_{j'+1}\geq T]} \Upsilon_{A^{\tau^*_{j'},\vartheta^*_{j'},\veca^*_{j'},\beta^*_{j'}}_T}(\theta_T^{\tau^*_{j'},\vartheta^*_{j'},z^*_{j'},\veca^*_{j'},\beta^*_{j'}})
\\
&+ \ett_{[\tau^*_{j'+1} < T]} \left\{-c_{\beta^*_{j'},\beta^*_{j'+1}}(\tau^*_{j'+1}) + Y^{\tau^*_{j'+1},\vartheta^*_{j'+1},z^*_{j'+1},\veca^*_{j'+1},\beta^*_{j'+1}}_{\tau^*_{j'+1}}\right\}\Big| \mcG_s\bigg].
\end{align*}
$\Prob$-a.s. By induction we get that for each $K\geq 0$
\begin{align*}
Y^{0,0,0,0,0}_0=\E\bigg[&\int_0^{\tau^*_K \wedge T}\sum_{j=0}^{K\wedge N^*} \ett_{[\tau^*_{j}\leq r < \tau^*_{j+1}]} \psi_{A^{\tau^*_{j},\vartheta^*_{j},\veca^*_{j},\beta^*_{j}}_r}(r,\theta_r^{\tau^*_{j},\vartheta^*_{j},z^*_{j},\veca^*_{j},\beta^*_{j}})dr
\\
&+\sum_{j=0}^{K\wedge N^*} \ett_{[\tau^*_{j}< T\leq \tau^*_{j+1}]}\Upsilon_{A^{\tau^*_{j},\vartheta^*_{j},\veca^*_{j},\beta^*_{j}}_T} (\theta_T^{\tau^*_{j},\vartheta^*_{j},z^*_{j},\veca^*_{j},\beta^*_{j}})
\\
&-\sum_{j=1}^{K\wedge N^*}c_{\beta^*_{j-1},\beta^*_{j}}(\tau^*_j)+\ett_{[\tau^*_{K+1} < T]}\{-c_{\beta^*_{K},\beta^*_{K+1}}(\tau^*_{K+1})+Y^{\tau^*_{K+1},\vartheta^*_{K+1},z^*_{K+1},\veca^*_{K+1},\beta^*_{K+1}}_{\tau^*_{K+1}}\}\bigg],
\end{align*}
where $\tau^*_{N^*+1}=\tau^*_{N^*+2}=\cdots=\infty$. Letting $N\to\infty$ while assuming that $u^*\in\mcU^f$ we find that $Y^{0,0,0,0,0}_0=J(u^*)$.\\

\noindent {\bf Step 3} It remains to show that the strategy $u^*$ is optimal. To do this we pick any other strategy $\hat u:=(\hat\tau_1,\ldots,\hat\tau_{\hat N};\hat\beta_1,\ldots,\hat\beta_{\hat N})\in\mcU^f$ and let the triple $(\hat\vartheta_j,\hat z_j,\hat\veca_j)_{1\leq j\leq \hat N}$ be defined by the recursions $\hat\vartheta_j:=\hat\vartheta_{j-1}\hat\beta_j+(\hat\beta_j-\hat\beta_{j-1})^+\hat\tau_j$, $\hat z_{j}:=\theta_{\hat \tau_{j}}^{\hat \tau_{j-1},\hat\vartheta_{j-1},\hat z_{j-1},\hat\veca_{j-1},\hat \beta_{j-1}}\wedge T\hat \beta_{j}$ and $\hat\veca_{j}:=A^{\hat \tau_{j-1},\hat\vartheta_{j-1},\hat\veca_{j-1},\hat \beta_{j-1}}_{\hat \tau_{j}}\wedge\hat \beta_j$, with $(\hat\tau_0,\hat\vartheta_0,\hat z_0,\hat\veca_0,\hat\beta_0):=(0,0,0,0,0)$. By the definition of $Y^{0,0,0,0,0}_0$ in~\eqref{ekv:Ydef} we have
\begin{align*}
Y^{0,0,0,0,0}_0 &\geq \E\bigg[\int_0^{\hat \tau_1}\psi_0\left(r,0\right)dr
+\ett_{[\hat \tau_1\geq T]}\Upsilon_0\left(0\right) + \ett_{[\hat \tau_1 < T]}\max_{\beta\in \mcI^{-0}}\left\{-c_{0,\beta}(\hat\tau_1)+Y^{\hat \tau_1,\hat\tau_1\beta,0,0,\beta}_{\hat \tau_1}\right\}\bigg]
\\
&\geq\E\bigg[\int_0^{\hat \tau_1}\psi_0\left(r,0\right)dr+\ett_{[\hat \tau_1\geq T]}\Upsilon_0\left(0\right) + \ett_{[\hat \tau_1 < T]}\left\{-c_{0,\hat\beta_1}(\hat \tau_1)+Y^{\hat \tau_1,\hat\vartheta_1,\hat z_1,\hat\veca_1,\hat \beta_1}_{\hat \tau_1}\right\}\bigg]
\end{align*}
but in the same way
\begin{align*}
Y^{\hat \tau_1,\hat\vartheta_1,\hat z_1,\hat\veca_1,\hat \beta_1}_{\hat \tau_1}\geq\E\bigg[&\int_{\hat\tau_1}^{\hat\tau_{2}}\psi_{A^{\hat \tau_1,\hat\vartheta_1,\hat z_1,\hat \beta_1}_r}(r,\theta_r^{\hat \tau_1,\hat\vartheta_1,\hat z_1,\hat\veca_1,\hat \beta_1})dr+\ett_{[\hat\tau_{2}\geq T]}\Upsilon_{A^{\hat \tau_1,\hat\vartheta_1,\hat z_1,\hat\veca_1,\hat \beta_1}_T}(\theta_T^{\hat \tau_1,\hat\vartheta_1,\hat z_1,\hat\veca_1,\hat \beta_1})
\\
&+ \ett_{[\hat\tau_{2} < T]}\left\{-c_{\hat\beta_{1},\hat\beta_{2}}(\hat\tau_{2})+Y^{\hat \tau_2,\hat\vartheta_2,\hat z_2,\hat\veca_2,\hat \beta_2}_{\hat\tau_2}\right\}\Big| \mcG_{\hat\tau_1}\bigg],
\end{align*}
$\Prob$--a.s. By repeating this argument and using the dominated convergence theorem we find that $J(u^*)\geq J(\hat u)$ which proves that $u^*$ is in fact optimal and thus belongs to $\mcU^f$.\qed

\bigskip

\begin{rem}
Note that the above proof can be trivially extended to arbitrary initial conditions $(z^*_0,\veca^*_0,\beta^*_0)\in\cup_{(\veca,\vecb)\in\mcAI}[0,T]^{\veca^+}\times(\veca,\vecb)$.
\end{rem}

%\bigskip
%
%\begin{rem}\label{rem:vfgenF}
%Note that we may replace \emph{\ref{vFAM:cont})} with the weaker assumption:
%\begin{enumerate}[\emph{c')}]
%  \item \emph{For each $(\veca,\vecb)\in\mcJ$ and each $(\nu,z)\in [0,T]^\vecb\times [0,T]^{\veca^+}$ the process $(Y_t^{t,\nu,z,\veca,\vecb}:0\leq t\leq T)$ is right continuous and LCE and for each $(t,s)\in [0,T]\times [0,T]$ the family $(\E[Y_s^{t,\nu,z,\veca,\vecb}|\mcG_t]:(\nu,z)\in [0,T]^\vecb\times [0,T]^{\veca^+})$ is $\Prob$-a.s.~continuous in $(\nu,z)$.}
%\end{enumerate}
%Furthermore, we can allow the filtration $\bbF$ the be more general, \eg generated by a Brownian motion and an independent Poisson random measure.
%\end{rem}

%%%%%%%%%%%%%%%%%%%%%%%%%%%%%%%%%%%%%%%%%%%%%%%%%%%%%%%%%%%%%%%%%%%%%%%%%%%%%%%%%%%%%%%%%%%%%%%%%%%%%%%%%%%%%%%%%%%%%%%%%%%%%%%%%%%%%%%%%%%%%%%%%%
%%%%%%%%%%%%%%%%%%%%%%%%%%%%%%%%%%%%%%%%%%%%%%%%%%%%%%%%%%%%%%%%%%%%%%%%%%%%%%%%%%%%%%%%%%%%%%%%%%%%%%%%%%%%%%%%%%%%%%%%%%%%%%%%%%%%%%%%%%%%%%%%%%
%%%%%%%%%%%%%%%%%%%%%%%%%%%%%%%%%%%%%%%%%%%%%%%%%%%%%%%%%%%%%%%%%%%%%%%%%%%%%%%%%%%%%%%%%%%%%%%%%%%%%%%%%%%%%%%%%%%%%%%%%%%%%%%%%%%%%%%%%%%%%%%%%%

\section{Existence\label{sec:exist}}
Theorem~\ref{thm:vfc} presumes existence of the family $((Y^{t,\nu,z,\veca,\vecb}_s)_{0\leq s\leq T}: (t,\nu,z,\veca,\vecb)\in \bigSET)$. To obtain a satisfactory solution to Problem 1, we thus need to establish that there exists a family of processes satisfying properties \emph{a)-d)} in the definition of a verification family. We will follow the standard existence proof which goes by applying a Picard iteration (see \cite{CarmLud,BollanMSwitch1,HamZhang}). We thus define a sequence $((Y^{t,\nu,z,\veca,\vecb,k}_s)_{0\leq s\leq T}: (t,\nu,z,\veca,\vecb)\in \bigSET)_{k\geq 0}$ of families of processes as
\begin{align}
Y^{t,\nu,z,\veca,\vecb,0}_s:=\E\bigg[\int_{s}^{T}\psi_{A^{t,\nu,\veca,\vecb}_r}(r,\theta_r^{t,\nu,z,\veca,\vecb})dr + \Upsilon_{A^{t,\nu,\veca,\vecb}_T}(\theta_T^{t,\nu,z,\veca,\vecb})\Big| \mcG_s\bigg]\label{ekv:Y0def}
\end{align}
and
\begin{align}\nonumber
Y^{t,\nu,z,\veca,\vecb,k}_s:=\esssup_{\tau \in \mcT_{s}} \E\bigg[&\int_s^{\tau\wedge T}\psi_{A^{t,\nu,\veca,\vecb}_r}\left(r,\theta_r^{t,\nu,z,\veca,\vecb}\right)dr+\ett_{[\tau \geq T]}\Upsilon_{A^{t,\nu,\veca,\vecb}_T}\left(\theta_T^{t,\nu,z,\veca,\vecb}\right)
\\
&+\ett_{[\tau < T]}\max_{\beta\in\mcI}\left\{-c_{\vecb,\beta}(\tau)+Y^{\tau,\beta\nu+\tau(\beta-\vecb)^+,\theta_\tau^{t,\nu,z,\veca,\vecb}\wedge T\beta,A^{t,\nu,\veca,\vecb}_\tau\wedge\beta,\beta,k-1}_\tau\right\}\Big| \mcG_s\bigg]\label{ekv:Ykdef}
\end{align}
for $k\geq 1$.

In this section we will show that the limiting family, $((\tilde Y^{t,\nu,z,\veca,\vecb}_s)_{0\leq s\leq T}: (t,\nu,z,\veca,\vecb)\in \bigSET)$, obtained when letting $k\to\infty$ is a verification family, thus proving existence of an optimal control for Problem 1. This will be done over a number of steps where we start by showing that for each $k$ the family defined by the above recursions satisfy properties \emph{a)-c)}. We then show that property \emph{d)} follows from Theorem~\ref{thm:Snell}.\ref{Snell:lim}. However, we start by showing that the above defined family is uniformly $L^2$-bounded. We let $\bar\psi:=\max_{\veca'\in\mcA}\max_{z\in[0,T]^{\veca'}}|\psi_{\veca'}(\cdot,z)|$, $\bar\Upsilon:=\max_{\veca'\in\mcA}\max_{z\in[0,T]^{\veca'}}|\Upsilon_{\veca'}(z)|$ and define
\begin{equation*}
\bar Y_t:=E\bigg[\int_0^{T}\bar\psi(r)dr+\bar\Upsilon\Big|\mcF_t\bigg].
\end{equation*}
We have the following:

\begin{prop}\label{prop:YkL2bnd}
For each $k\geq 0$, the family of processes $\big(Y^{t,\nu,z,\veca,\vecb,k}: (t,\nu,z,\veca,\vecb)\in \bigSET\big)$ is $L^2$-bounded in the sense that there is a constant $K_Y>0$ such that
\begin{equation}\label{ekv:YkL2bnd}
\E\Big[\sup_{(t,\nu,z,\veca,\vecb)\in \bigSET}\sup_{s\in[0,T]}|Y^{t,\nu,z,\veca,\vecb,k}_s|^2\Big]\leq K_Y<\infty,
\end{equation}
for all $k\geq 0$.
\end{prop}
\noindent \emph{Proof.} Let $\bar Y^k=\sup_{(t,\nu,z,\veca,\vecb)\in \bigSET}|Y^{t,\nu,z,\veca,\vecb,k}|$. Since $Y^{t,\nu,z,\veca,\vecb,k}\geq 0$, applying an induction argument gives
\begin{align*}
\bar Y^k_s&\leq \E\bigg[\int_s^{T}\max_{\veca\in\mcA}\max_{z\in [0,T]^{\veca^+}}|\psi_{\veca}(r,z)|dr
+\max_{\veca\in\mcA}\max_{z\in [0,T]^{\veca^+}}|\Upsilon_{\veca}(z)|\Big|\mcF_s\bigg]\leq \bar Y_s.
\end{align*}
Hence, by Doob's maximal inequality we get
\begin{align*}
\E\Big[\sup_{s\in [0,T]}|\bar Y^k_s|^2\Big]\leq C\E\Big[|\bar Y_0|^2 \Big]\leq C\E\bigg[&\int_0^{T}\bar\psi(r)^2dr
+\bar\Upsilon^2 \bigg].
\end{align*}
where the right hand side is bounded by Assumption~\ref{ass:prelim}.\qed\\

\bigskip

It should be noted that the above bound is uniform in $k$ which implies that the limit family (if it exists) satisfies the same inequality. In particular, we conclude that property \emph{\ref{vFAM:bnd})} holds for all $k$. Properties \emph{\ref{vFAM:inS2e})} and \emph{\ref{vFAM:cont})} will be shown by induction and we make the following induction hypothesis:\\

\noindent{\bf H.k.} \emph{The family $((Y^{t,\nu,z,\veca,\vecb,k}_s)_{0\leq s\leq T}: (t,\nu,z,\veca,\vecb)\in \bigSET)$ is such that:
\begin{enumerate}[i)]
  \item\label{hyp:inS2e} For every $(t,\nu,z,\veca,\vecb)\in \bigSET$ we have $Y^{t,\nu,z,\veca,\vecb,k}\in\mcS_e^2$ and $(Y^{s,\nu,z,\veca,\vecb,k}_s:0\leq s\leq T)\in\mcS_e^2$.
  \item\label{hyp:cont} For every $(\veca,\vecb)\in\mcJ$, we have to following continuity property
  \begin{equation*}
    \lim_{(p,q)\to(0,0)}\E\Big[\sup_{(t,\nu,z)\in \mcD_{(\veca,\vecb)}}|\E[Y^{t,(\nu+p)^+\wedge \vecb T,(z+q)^+\wedge T\veca^+,\veca,\vecb,k}_{t}-Y^{t,\nu,z,\veca,\vecb,k}_t|\mcF_t]|^2\Big]= 0.
  \end{equation*}
\end{enumerate}}

\bigskip

We note that, under the induction hypotheses {\bf H.0}-{\bf H.k}, arguing as in the proof of Theorem~\ref{thm:vfc} we have
\begin{align}
Y^{t,\nu,z,\veca,\vecb,k+1}_s=\esssup_{u\in \mcU_s^{k+1}} \E\bigg[&\int_s^{T}\psi_{\alpha_r^{t,\nu,\veca,\vecb,u}}(r,\theta^{t,\nu,z,\veca,\vecb,u}_r)dr
+\Upsilon_{\alpha_T^{t,\nu,\veca,\vecb,u}}(\theta^{t,\nu,z,\veca,\vecb,u}_T)
-\sum_{j=1}^N c_{\beta_{j-1},\beta_{j}}(\tau_j)\Big|\mcG_s\bigg],\label{ekv:Ykalt}
\end{align}
where the supremum is attained by a control $u^{*,k+1}\in \mcU_s^{k+1}$. Furthermore, the characterisation of $Y^{t,\nu,z,\veca,\vecb,k}_t$ in terms of the recursion \eqref{ekv:Y0def}-\eqref{ekv:Ykdef} can be further simplified by noting that the possible transitions of $A^{t,\nu,\veca,\vecb}_\cdot$ forms paths in a directed acyclic graph where the leafs are the members of the set $\mcA^{\rm abs}_{\vecb}$. Letting $\eta\in\mcT_t$ denote the first transition time of $A^{t,\nu,\veca,\vecb}_\cdot$, Theorem~\ref{thm:Snell}.\ref{Snell:DynP} allows us to write (recall that $\theta_s^{t,z,\veca}=z+\veca^+(s-t)^+$)
\begin{align}\nonumber
Y^{t,\nu,z,\veca,\vecb,k}_t=\esssup_{\tau \in \mcT_{t}} \E\bigg[&\int_t^{\eta\wedge\tau\wedge T}\psi_{\veca}(r,\theta_r^{t,z,\veca})dr+\ett_{[\eta> T]}\ett_{[\tau \geq T]}\Upsilon_{\veca}(\theta_T^{t,z,\veca})
\\ \nonumber
&+\ett_{[\eta \leq\tau\wedge T]}Y^{\eta,\nu,\theta_\eta^{t,z,\veca}\wedge TA^{t,\nu,\veca,\vecb}_\eta,A^{t,\nu,\veca,\vecb}_\eta,\vecb,k}_\eta
\\
&+\ett_{[\tau < T\wedge \eta]}\max_{\beta\in\mcI}\left\{-c_{\vecb,\beta}(\tau)+Y^{\tau,\beta\nu+\tau(\beta-\vecb)^+,\theta_\tau^{t,z,\veca}\wedge T\beta,\veca\wedge\beta,\beta,k-1}_\tau\right\}\Big| \mcF_t\bigg]\label{ekv:Ykdef3}
\end{align}
and
\begin{align}
Y^{t,\nu,z,\veca,\vecb,0}_t=\E\bigg[&\int_t^{\eta \wedge T}\psi_{\veca}(r,\theta_r^{t,z,\veca})dr+\ett_{[\eta> T]}\Upsilon_{\veca}(\theta_T^{t,z,\veca})+\ett_{[\eta\leq T]}Y^{\eta,\nu,\theta_\eta^{t,z,\veca}\wedge TA^{t,\nu,\veca,\vecb}_\eta,A^{t,\nu,\veca,\vecb}_\eta,\vecb,0}_\eta\Big| \mcF_t\bigg].\label{ekv:Y0def3}
\end{align}
Actually, since $A^{t,\nu,\veca,\vecb}$ is a pure jump Markov process, we have that $\sigma(A^{t,\nu,\veca,\vecb}_s:0\leq s\leq \cdot)$ is generated by sets of the type $\{s<\eta\}$ on $[0,\eta)$. Hence, for each $\tau\in\mcT_t$ there is a $\tilde\tau\in\mcT^{\bbF}_{t}$ such that $\tilde\tau\wedge\eta=\tau\wedge\eta$, $\Prob$-a.s.~and we only need to take the essential supremum over $\mcT^{\bbF}_t$ (the set of $\bbF$-stopping times $\tau\geq t$) in \eqref{ekv:Ykdef3}.

Furthermore, we note that when $\veca\in\mcA^{\rm abs}_\vecb$, then by the definition of $A^{t,\nu,\veca,\vecb}$ we have $\eta=\infty$ and thus
\begin{align}\nonumber
Y^{t,\nu,z,\veca,\vecb,k}_t=\esssup_{\tau \in \mcT^{\bbF}_{t}} \E\bigg[&\int_t^{\tau\wedge T}\psi_{\veca}(r,\theta_r^{t,z,\veca})dr+\ett_{[\tau \geq T]}\Upsilon_{\veca}(\theta_T^{t,z,\veca})
\\
&+\ett_{[\tau < T]}\max_{\beta\in\mcI}\left\{-c_{\vecb,\beta}(\tau)+Y^{\tau,\beta\nu+\tau(\beta-\vecb)^+,\theta_\tau^{t,z,\veca}\wedge T\beta,\veca\wedge\beta,\beta,k-1}_\tau\right\}\Big| \mcF_t\bigg]\label{ekv:Ykdef4}
\end{align}
and
\begin{align}
Y^{t,\nu,z,\veca,\vecb,0}_t=\E\bigg[&\int_t^{T}\psi_{\veca}(r,\theta_r^{t,z,\veca})dr+\Upsilon_{\veca}(\theta_T^{t,z,\veca})\Big| \mcF_t\bigg].\label{ekv:Y0def5}
\end{align}

%%%%%%%%%%%%%%%%%%%%%%%%%%%%%%%%%%%%%%%%%%%%%%%%%%%%%%%%%%%%%%%%%%%%%%%%%%%%%%%%%%%%%%%%%%%%%%%%%%%%%%%%%%%%%%%%%%%%%%%%%%%%%%%%%%%%%%%%%%%%%%%%%%

We start by showing that the induction hypothesis holds for $k=0$:

\begin{prop}\label{prop:H0}
The induction hypothesis {\bf H.0} holds.
\end{prop}

\noindent \emph{Proof.} We first show continuity in $z$. This follows by the Lipschitz assumptions on $\psi$ and $\Upsilon$ together with the definition of $\theta^{t,\nu,z,\veca,\vecb}$. Indeed we have
\begin{equation*}
|\theta^{t,\nu,z',\veca,\vecb}_s-\theta^{t,\nu,z,\veca,\vecb}_s|\leq |z'-z|,
\end{equation*}
for all $s\in[0,T]$. It thus follows immediately from \eqref{ekv:Y0def} that
\begin{align*}
\sup_{(t,\nu,z)\in\mcD_{(\veca,\vecb)}}|Y^{t,\nu,(z+q)^+\wedge T,\veca,\vecb,0}_{t}-Y^{t,\nu,z,\veca,\vecb,0}_{t}|\leq C|q|,
\end{align*}
$\Prob$-a.s.~and in particular we have
\begin{align*}
\E\bigg[\sup_{(t,\nu,z)\in\mcD_{(\veca,\vecb)}}|Y^{t,\nu,(z+q)^+\wedge T,\veca,\vecb,0}_{t}-Y^{t,\nu,z,\veca,\vecb,0}_{t}|^2\bigg]\leq C|q|.
\end{align*}

\bigskip

%%%%%%%%%%%%%%%%%%%%%%%%%%%%%%%%%%%%%%%%%%%%%%%%%%%%%%%%%%%%%%%%%%%%%%%%%%%%%%%%%%%%%%%%%%%%%%%%%%%%%%%%%%%%%%%%%%%%%%%%%%%%%%%%%%%%%%%%%%%%%%%%%%

Next we show continuity in $\nu$. We have, with $\eta$ the first transition time of $A^{t,\nu,\veca,\vecb}$ and $\eta'$ the first transition time of $A^{t,\nu',\veca,\vecb}$,
\begin{align*}
Y^{t,\nu,z,\veca,\vecb,0}_{t}-Y^{t,\nu',z,\veca,\vecb,0}_{t}&=\E\bigg[\int_{t}^{\eta \wedge T}\psi_{\veca}(s,\theta_s^{t,z,\veca})ds-\int_{t}^{\eta' \wedge T}\psi_{\veca}(s,\theta_s^{t,z,\veca})ds+\ett_{[\eta> T]}\Upsilon_{\veca}(\theta_T^{t,z,\veca})
\\
&-\ett_{[\eta'> T]}\Upsilon_{\veca}(\theta_T^{t,z,\veca})+\ett_{[\eta\leq T]}Y^{\eta,\nu,\theta_{\eta}^{t,z,\veca}\wedge TA^{t,\nu,\veca,\vecb}_\eta,A^{t,\nu,\veca,\vecb}_\eta,\vecb,0}_\eta
\\
&-\ett_{[\eta'\leq T]}Y^{\eta',\nu',\theta_{\eta'}^{t,z,\veca}\wedge TA^{t,\nu',\veca,\vecb}_{\eta'},A^{t,\nu',\veca,\vecb}_{\eta'},\vecb,0}_{\eta'}\Big| \mcF_{t}\bigg]
\\%%%%%%%%%%%%%%%%%%%%%%%%%%%%%%%%%%%%%%%%%%%%%%%%%
&=\E\bigg[\int_{t}^{T}\Big(e^{\int_{t}^s\lambda^{\nu,\vecb}_{\veca,\veca}(r)dr} -e^{\int_{t}^s\lambda^{\nu',\vecb}_{\veca,\veca}(r)dr}\Big)\psi_{\veca}(s,\theta_s^{t,z,\veca})ds
\\
&+(e^{\int_{t}^T\lambda^{\nu,\vecb}_{\veca,\veca}(r)dr} - e^{\int_{t}^T\lambda^{\nu',\vecb}_{\veca,\veca}(r)dr}) \Upsilon_{\veca}(\theta_T^{t,z,\veca})
\\
&+\sum_{\veca'\in\mcA^{-\veca}_{\veca,\vecb}}\int_{t}^{T}\Big(\lambda^{\nu,\vecb}_{\veca,\veca'}(s) e^{\int_{t}^s\lambda^{\nu,\vecb}_{\veca,\veca}(r)dr}Y^{s,\nu,\theta_{s}^{t,z,\veca}\wedge T\veca',\veca',\vecb,0}_s
\\
&-\lambda^{\nu',\vecb}_{\veca,\veca'}(s)e^{\int_{t}^s\lambda^{\nu',\vecb}_{\veca,\veca}(r)dr}Y^{s,\nu',\theta_{s}^{t,z,\veca}\wedge T\veca',\veca',\vecb,0}_s\Big) ds \Big| \mcF_{t}\bigg]
\end{align*}
Using the identity $ab-a'b'=1/2((a-a')(b+b')+(a+a')(b-b'))$ gives
\begin{align*}
Y^{t,\nu,z,\veca,\vecb,0}_{t}-Y^{t,\nu',z,\veca,\vecb,0}_{t}&\leq \E\bigg[\int_{t}^{T}|e^{\int_{t}^s\lambda^{\nu,\vecb}_{\veca,\veca}(r)dr}-e^{\int_{t}^s\lambda^{\nu',\vecb}_{\veca,\veca}(r)dr}|\cdot |\psi_{\veca}(s,\theta_s^{t,z,\veca})|ds
\\
&+|e^{\int_{t}^T\lambda^{\nu,\vecb}_{\veca,\veca}(r)dr}-e^{\int_{t}^T\lambda^{\nu',\vecb}_{\veca,\veca}(r)dr}|\cdot |\Upsilon_{\veca}(\theta_T^{t,z,\veca})|
\\
&+\sum_{\veca'\in\mcA^{-\veca}_{\veca,\vecb}}\int_{t}^{T}\Big(|\lambda^{\nu,\vecb}_{\veca,\veca'}(s)e^{\int_{t}^s\lambda^{\nu,\vecb}_{\veca,\veca}(r)dr}-\lambda^{\nu',\vecb}_{\veca,\veca'}(s)e^{\int_{t}^s\lambda^{\nu',\vecb}_{\veca,\veca}(r)dr}|
\\
&\cdot |Y^{s,\nu,\theta_{s}^{t,z,\veca}\wedge T\veca',\veca',\vecb,0}_s+Y^{s,\nu',\theta_{s}^{t,z,\veca}\wedge T\veca',\veca',\vecb,0}_s|
\\
&+2K_{\lambda}|Y^{s,\nu,\theta_{s}^{t,z,\veca}\wedge T\veca',\veca',\vecb,0}_s-Y^{s,\nu',\theta_{s}^{t,z,\veca}\wedge T\veca',\veca',\vecb,0}_s|\Big)ds\Big| \mcF_{t}\bigg].
\end{align*}
Noting that the same equality holds for $Y^{t,\nu',z,\veca,\vecb,0}_{t}-Y^{t,\nu,z,\veca,\vecb,0}_{t}$ an induction argument gives that
\begin{align*}
|Y^{t,\nu,z,\veca,\vecb,0}_{t}-Y^{t,\nu',z,\veca,\vecb,0}_{t}|&\leq C|\nu-\nu'|\bar Y_t,
\end{align*}
$\Prob$-a.s. Applying Doob's maximal inequality, we find that, for $p\in \R^{n}$
\begin{align*}
\E\bigg[\sup_{(t,\nu,z)\in\mcD_{(\veca,\vecb)}}|Y^{t,(\nu+p)^+\wedge T\vecb,z,\veca,\vecb,0}_{t}-Y^{t,\nu,z,\veca,\vecb,0}_{t}|^2\bigg] \leq C|p|.
\end{align*}
Put together, this implies that \textbf{H.0}.\ref{hyp:cont} holds.\\

%%%%%%%%%%%%%%%%%%%%%%%%%%%%%%%%%%%%%%%%%%%%%%%%%%%%%%%%%%%%%%%%%%%%%%%%%%%%%%%%%%%%%%%%%%%%%%%%%%%%%%%%%%%%%%%%%%%%%%%%%%%%%%%%%%%%%%%%%%%%%%%%%%

To show that \textbf{H.0}.\ref{hyp:inS2e} holds we note that
\begin{align*}
Y^{t,\nu,z,\veca,\vecb,0}_s&=\E\bigg[\int_0^{T}\psi_{\alpha_r^{t,\nu,\veca,\vecb}}(r,\theta_r^{t,z,\veca,\vecb})dr +\Upsilon_{\alpha_T^{t,\nu,\veca,\vecb}}(\theta_T^{t,z,\veca,\vecb})\Big| \mcG_s\bigg]-\int_0^s\psi_{\alpha_r^{t,\nu,\veca,\vecb}}(r,\theta_r^{t,z,\veca,\vecb})dr
\end{align*}
is the sum of a continuous process and a martingale in a quasi-left continuous filtration, hence, it has a version which is \cadlag and LCE and, thus, belongs to $\mcS^2_e$ by Proposition~\ref{prop:YkL2bnd}.

Let $(t,\nu,z)\in \mcD_{(\veca,\vecb)}$ and assume that $t' \in [t,T]$. To evaluate $|Y^{t,\nu,z,\veca,\vecb,0}_{t}-Y^{t',\nu,z,\veca,\vecb,0}_{t'}|$ we note that during $[t,t']$ we may have a transition in $A^{t,\nu,\veca,\vecb}$ (the probability of which is bounded by $1-e^{-K_\lambda(t'-t)}$) and thus have
\begin{align*}
|Y^{t,\nu,z,\veca,\vecb,0}_t-&Y^{t',\nu,z,\veca,\vecb,0}_{t'}|\leq \E\Big[\int_t^{t'}\bar\psi(s)ds\big|\mcF_t\Big]+2(1-e^{-K_\lambda(t'-t)})\bar Y_t
\\
& \quad+ |\E\Big[Y^{t',\nu,z,\veca,\vecb,0}_{t'}\big|\mcF_t\Big] -Y^{t',\nu,z,\veca,\vecb,0}_{t'}|+C|t'-t|.
\end{align*}
Note that the first two terms on the right hand side in the above equation satisfies:
\begin{align*}
\E\Big[\int_{t'}^{t}\bar\psi(s)ds\big|\mcF_t\Big]+2(1-e^{-K_\lambda |t'-t|})\bar Y_t \to 0,
\end{align*}
$\Prob$-a.s.~as $|t'-t|\to 0$ by $\Prob$-a.s.~boundedness of $\bar\psi$ and $\bar Y$.\\

Concerning the third term let us first consider the case when $\veca\in\mcAabs_\vecb$. Define $\mcK_{t,t'}(X):=|\E\Big[X\big|\mcF_t\Big]-\E\Big[X\big|\mcF_{t'}\Big]|$, then $\mcK_{t,t'}$ is a subadditive operator and
\begin{align*}
&|\E\Big[Y^{t',\nu,z,\veca,\vecb,0}_{t'}\big|\mcF_t\Big] -Y^{t',\nu,z,\veca,\vecb,0}_{t'}|\leq \sup_{(r,p,q)\in\mcD_{(\veca,\vecb)}} \mcK_{t,t'}(Y^{r,p,q,\veca,\vecb,0}_{r})
\\
&\leq \int_0^T\sup_{q\in[0,T]^{\veca^+}}\mcK_{t,t'}(\psi_{\veca}(s,q))ds + \sup_{q\in[0,T]^{\veca^+}}\mcK_{t,t'}(\Upsilon_{\veca}(q))
\\
&\leq \sum_{l=1}^{M}\bigg(\int_0^T\mcK_{t,t'}(\psi_{\veca}(s,z^{M}_l))ds + \mcK_{t,t'}(\Upsilon_{\veca}(z^{M}_l))\bigg)+C\delta(M),
\end{align*}
where $\delta(M)$ is the diameter of a partition $(G^M_l)_{l=1}^M$ of $[0,T]^{\veca^+}$ and $z^M_l\in G^M_l$. The last of the above inequalities follows by noting that for all $z,z'\in [0,T]^{\veca^+}$
\begin{align*}
\mcK_{t,t'}(\psi_{\veca}(s,z))&=\mcK_{t,t'}(\psi_{\veca}(s,z')+(\psi_{\veca}(s,z)-\psi_{\veca}(s,z')))
\\
&\leq \mcK_{t,t'}(\psi_{\veca}(s,z'))+\mcK_{t,t'}(\psi_{\veca}(s,z)-\psi_{\veca}(s,z')))
\\
&\leq \mcK_{t,t'}(\psi_{\veca}(s,z'))+2k_\psi|z'-z|.
\end{align*}
Now, by the martingale representation theorem there is, for each $(s,l)\in [0,T]\times \{1,\ldots,M\}$, a process $Z^{\psi,\veca,s,l,M}\in\mcH^2_{\bbF}$ such that $\E\Big[\psi_{\veca}(s,z^M_l)\big|\mcF_t\Big]=\E\Big[\psi_{\veca}(s,z^M_l)\Big]+\int_0^t Z^{\psi,\veca,s,l,M}_rdW_r$ and a process $Z^{\Upsilon,\veca,l,M}\in\mcH^2_{\bbF}$ such that $\E\Big[\Upsilon_{\veca}(z^M_l)\big|\mcF_t\Big]=\E\Big[\Upsilon_{\veca}(z^M_l)\Big]+\int_0^t Z^{\Upsilon,\veca,l,M}_rdW_r$. We thus have
\begin{align*}
&|\E\Big[Y^{t',\nu,z,\veca,\vecb,0}_{t'}\big|\mcF_t\Big] -Y^{t',\nu,z,\veca,\vecb,0}_{t'}|\leq
\sum_{l=1}^M\bigg(\int_0^T|\int_{t}^{t'}Z^{\psi,\veca,s,l,M}_rdW_r|ds + |\int_{t}^{t'}Z^{\Upsilon,\veca,l,M}_rdW_r|\bigg)+C\delta(M).
\end{align*}
By dominated convergence we thus get that
\begin{align*}
\lim_{h\to 0}\sup_{(t,\nu,z)\in \mcD_{(\veca,\vecb)}} |Y^{(t+h)^+\wedge T,\nu,z,\veca,\vecb,0}_{(t+h)^+\wedge T}-Y^{t,\nu,z,\veca,\vecb,0}_t|\leq C\delta(M),
\end{align*}
$\Prob$-a.s., where $\delta(M)$ can be made arbitrarily small. Now assume that $\veca\notin\mcAabs_\vecb$ and let $(G_l^M)_{l=1}^M$ be a partition of $\mcD_{(\veca,\vecb)}$ with $\max_l {\rm diam}(G_{l}^M)=\delta(M)\to 0$ as $M\to\infty$ and let $(t_l^M,\nu^M_l,z^M_l)\in G_l^M$ then
\begin{align*}
&|\E\Big[Y^{t',\nu,z,\veca,\vecb,0}_{t'}\big|\mcF_t\Big] -Y^{t',\nu,z,\veca,\vecb,0}_{t'}|
\\
%%%%%%%%%%%%%%%%%%%%%%%%%%%%%%%%%%%%%%%%%%%%%%%%%%%%%%%%%%%%%%%%%%%%%%%%%%%%%%%%%%%%%%%%%%%%%%%%%%%%%%%%%%%%%%%%%%%%%%%%%%%%%%%%%%%%%%%%%%%%%%%%%%
&\leq \int_0^T\sup_{(r,p,q)\in \mcD_{(\veca,\vecb)}}\mcK_{t,t'}(e^{\int_r^s\lambda^{p,\vecb}_{\veca,\veca}(v)dv}\psi_{\veca}(s,q))ds + \sup_{(r,p,q)\in \mcD_{(\veca,\vecb)}}\mcK_{t,t'}(e^{\int_r^T\lambda^{p,\vecb}_{\veca,\veca}(v)dv}\Upsilon_{\veca}(q))
\\
&\quad + \sum_{\veca'\in\mcA^{-\veca}_{\veca,\vecb}}\int_0^T\sup_{(r,p,q)\in \mcD_{(\veca,\vecb)}}\mcK_{t,t'} (\lambda^{p,\vecb}_{\veca,\veca'}(s)e^{\int_r^s\lambda_{\veca,\veca}^{p,\vecb}(v)dv} Y^{s,p,q,\veca',\vecb,0}_{s})ds
\\
%%%%%%%%%%%%%%%%%%%%%%%%%%%%%%%%%%%%%%%%%%%%%%%%%%%%%%%%%%%%%%%%%%%%%%%%%%%%%%%%%%%%%%%%%%%%%%%%%%%%%%%%%%%%%%%%%%%%%%%%%%%%%%%%%%%%%%%%%%%%%%%%%%
&\leq \sum_{l=1}^M\bigg(\int_0^T\mcK_{t,t'}(e^{\int_{t^M_l}^s\lambda^{\nu_l^M,\vecb}_{\veca,\veca}(v)dv}\psi_{\veca}(s,z^M_l))ds + \mcK_{t,t'}(e^{\int_{t^M_l}^T\lambda^{\nu_l^M,\vecb}_{\veca,\veca}(v)dv}\Upsilon_{\veca}(z^M_l))
\\
&\quad + \sum_{\veca'\in\mcA^{-\veca}_{\veca,\vecb}} \int_0^T\mcK_{t,t'} (\lambda^{\nu_l^M,\vecb}_{\veca,\veca'}(s)e^{\int_{t^M_l}^{s}\lambda^{\nu_l^M,\vecb}_{\veca,\veca}(v)dv} Y^{s,\nu_l^M,z^M_l,\veca',\vecb,0}_{s})ds\bigg)+C(1+\bar Y_t+\bar Y_{t'})\delta(M).
%\\
%&\quad + C \sum_{\veca'\in\mcA^{-\veca}_{\veca,\vecb}} \sup_{r\in[0,T]}\E\big[\sup_{(s,p,q)\in \mcD_{(\veca,\vecb)}}\sup_{|(p',q')-(p,q)|\leq \delta(M)}|Y^{s,p,q,\veca',\vecb,0}_{s} - Y^{s,p',q',\veca',\vecb,0}_{s}|\big|\mcF_r\big]
\end{align*}
where, to arrive at the last inequality, we have used that for $(r,\nu,z),(r',\nu',z')\in\mcD_{(\veca,\vecb)}$, with $r\leq r'$,
\begin{align*}
&\mcK_{t,t'}(e^{\int_{r}^s\lambda^{\nu,\vecb}_{\veca,\veca}(v)dv}\psi_{\veca}(s,z)) \leq \mcK_{t,t'}(e^{\int_{r}^s\lambda^{\nu,\vecb}_{\veca,\veca}(v)dv}\psi_{\veca}(s,z')) + 2k_\psi|z'-z|
\\
%%%%%%%%%%%%%%%%%%%%%%%%%%%%%%%%%%%%%%%%%%%%%%%%%%%%%%%%%%%%%%%%%%%%%%%%%%%%%%%%%%%%%%%%%%%%%%%%%%%%%%%%%%%%%%%%%%%%%%%%%%%%%%%%%%%%%%%%%%%%%%%%%%
&\leq\mcK_{t,t'}((e^{\int_{r}^{s}\lambda^{\nu,\vecb}_{\veca,\veca}(v)dv}-e^{\int_{r'}^s\lambda^{\nu',\vecb}_{\veca,\veca}(v)dv}) \psi_{\veca}(s,z') + e^{\int_{r'}^s\lambda^{\nu',\vecb}_{\veca,\veca}(v)dv}\psi_{\veca}(s,z')) + 2k_\psi|z'-z|
\\
%&\leq \mcK_{t,t'}(e^{\int_{r'}^s\lambda^{\nu,\vecb}_{\veca,\veca}(v)dv}\psi_{\veca}(s,z')) + K_\lambda |r'-r| (\E\big[|\psi_{\veca}(s,z')|\big|\mcF_t\big] +\E\big[|\psi_{\veca}(s,z')|\big|\mcF_{t'}\big])  + 2k_\psi|z'-z|
%\\
%&=\mcK_{t,t'}((e^{\int_{r'}^s\lambda^{\nu,\vecb}_{\veca,\veca}(v)dv}-e^{\int_{r'}^s\lambda^{\nu',\vecb}_{\veca,\veca}(v)dv})\psi_{\veca}(s,z') + e^{\int_{r'}^s\lambda^{\nu',\vecb}_{\veca,\veca}(v)dv}\psi_{\veca}(s,z'))
%\\
%&\quad + K_\lambda |r'-r| (\E\big[|\psi_{\veca}(s,z')|\big|\mcF_t\big] +\E\big[|\psi_{\veca}(s,z')|\big|\mcF_{t'}\big])  + 2k_\psi|z'-z|
%\\
%%%%%%%%%%%%%%%%%%%%%%%%%%%%%%%%%%%%%%%%%%%%%%%%%%%%%%%%%%%%%%%%%%%%%%%%%%%%%%%%%%%%%%%%%%%%%%%%%%%%%%%%%%%%%%%%%%%%%%%%%%%%%%%%%%%%%%%%%%%%%%%%%%
&\leq \mcK_{t,t'}( e^{\int_{r'}^s\lambda^{\nu',\vecb}_{\veca,\veca}(v)dv}\psi_{\veca}(s,z')) + (K_\lambda |r'-r|+Tk_\lambda |\nu'-\nu|) (\E\big[|\psi_{\veca}(s,z')|\big|\mcF_t\big] +\E\big[|\psi_{\veca}(s,z')|\big|\mcF_{t'}\big])
\\
&\quad + 2k_\psi|z'-z|
\end{align*}
and
\begin{align*}
&\mcK_{t,t'}(\lambda^{\nu,\vecb}_{\veca,\veca'}(s)e^{\int_{r}^s\lambda^{\nu,\vecb}_{\veca,\veca}(v)dv}Y^{s,\nu,z,\veca',\vecb,k}_{s})
\\
%%%%%%%%%%%%%%%%%%%%%%%%%%%%%%%%%%%%%%%%%%%%%%%%%%%%%%%%%%%%%%%%%%%%%%%%%%%%%%%%%%%%%%%%%%%%%%%%%%%%%%%%%%%%%%%%%%%%%%%%%%%%%%%%%%%%%%%%%%%%%%%%%%
&=\mcK_{t,t'}((\lambda^{\nu,\vecb}_{\veca,\veca'}(s)e^{\int_{r}^s\lambda^{\nu,\vecb}_{\veca,\veca}(v)dv} - \lambda^{\nu',\vecb}_{\veca,\veca'}(s)e^{\int_{r'}^s\lambda^{\nu',\vecb}_{\veca,\veca}(v)dv})Y^{s,\nu,z,\veca',\vecb,k}_{s}
\\
&\quad + \lambda^{\nu',\vecb}_{\veca,\veca'}(s)e^{\int_{r'}^s\lambda^{\nu',\vecb}_{\veca,\veca}(v)dv}(Y^{s,\nu,z,\veca',\vecb,k}_{s} -Y^{s,\nu',z',\veca',\vecb,k}_{s}) + \lambda^{\nu',\vecb}_{\veca,\veca'}(s)e^{\int_{r'}^s\lambda^{\nu',\vecb}_{\veca,\veca}(v)dv}Y^{s,\nu',z',\veca',\vecb,k}_{s})
\\
%%%%%%%%%%%%%%%%%%%%%%%%%%%%%%%%%%%%%%%%%%%%%%%%%%%%%%%%%%%%%%%%%%%%%%%%%%%%%%%%%%%%%%%%%%%%%%%%%%%%%%%%%%%%%%%%%%%%%%%%%%%%%%%%%%%%%%%%%%%%%%%%%%
&\leq \mcK_{t,t'}(\lambda^{\nu',\vecb}_{\veca,\veca'}(s)e^{\int_{r'}^s\lambda^{\nu',\vecb}_{\veca,\veca}(v)dv}Y^{s,\nu',z',\veca',\vecb,k}_{s}) + (K_\lambda |r'-r|+Tk_\lambda |\nu'-\nu|) (\bar Y^k_t + \bar Y^k_{t'})
\\
&\quad + K_\lambda (\E\big[|Y^{s,\nu,z,\veca',\vecb,k}_{s} - Y^{s,\nu',z',\veca',\vecb,k}_{s}|\big|\mcF_t\big] +\E\big[|Y^{s,\nu,z,\veca',\vecb,k}_{s} -Y^{s,\nu',z',\veca',\vecb,k}_{s}|\big|\mcF_{t'}\big]).
\end{align*}
while noting that, by the above results, $|Y^{s,\nu,z,\veca',\vecb,0}_{s} -Y^{s,\nu',z',\veca',\vecb,0}_{s}|\leq C(|z'-z|+|\nu'-\nu|\bar Y_s)$. We conclude that
\begin{align*}
\lim_{h\to 0}\sup_{(t,\nu,z)\in \mcD_{(\veca,\vecb)}} |Y^{(t+h)^+\wedge T,\nu,z,\veca,\vecb,0}_{(t+h)^+\wedge T}-Y^{t,\nu,z,\veca,\vecb,0}_t|\leq C\delta(M),
\end{align*}
$\Prob$-a.s.~for all $(\veca,\vecb)\in\mcJ$. This implies that $(Y^{t,\nu,z,\veca,\vecb}_t:0\leq t\leq T)$ is in fact continuous and Proposition~\ref{prop:YkL2bnd} guarantees that it belongs to $\mcS^2_e$.\qed\\

%%%%%%%%%%%%%%%%%%%%%%%%%%%%%%%%%%%%%%%%%%%%%%%%%%%%%%%%%%%%%%%%%%%%%%%%%%%%%%%%%%%%%%%%%%%%%%%%%%%%%%%%%%%%%%%%%%%%%%%%%%%%%%%%%%%%%%%%%%%%%%%%%%
%%%%%%%%%%%%%%%%%%%%%%%%%%%%%%%%%%%%%%%%%%%%%%%%%%%%%%%%%%%%%%%%%%%%%%%%%%%%%%%%%%%%%%%%%%%%%%%%%%%%%%%%%%%%%%%%%%%%%%%%%%%%%%%%%%%%%%%%%%%%%%%%%%
%%%%%%%%%%%%%%%%%%%%%%%%%%%%%%%%%%%%%%%%%%%%%%%%%%%%%%%%%%%%%%%%%%%%%%%%%%%%%%%%%%%%%%%%%%%%%%%%%%%%%%%%%%%%%%%%%%%%%%%%%%%%%%%%%%%%%%%%%%%%%%%%%%

\begin{prop}\label{prop:Hk}
The induction hypothesis {\bf H.k} holds for all $k\geq 0$.
\end{prop}

\noindent\emph{Proof.} Since we have already shown that the induction hypothesis holds for $k=0$ we will assume that \textbf{H.0}-\textbf{H.k} hold for some $k\geq 0$. Now, as noted above, this implies the existence of a control $u^\diamond:=(\tau^\diamond_1,\ldots,\tau^\diamond_{N^\diamond};\beta_1^\diamond,\ldots,\beta_{N^\diamond}^\diamond)$, with ${N^\diamond}\leq k+1$, $\Prob$-a.s.~such that
\begin{align*}
Y^{t,\nu,z,\veca,\vecb,k+1}_t= \E\bigg[&\int_s^{T}\psi_{\alpha_r^{t,\nu,\veca,\vecb,u^\diamond}}(r,\theta^{t,\nu,z,\veca,\vecb,u^\diamond}_r)dr
+\Upsilon_{\alpha_T^{t,\nu,\veca,\vecb,u^\diamond}}(\theta^{t,\nu,z,\veca,\vecb,u^\diamond}_T)
-\sum_{j=1}^{N^\diamond} c_{\beta^\diamond_{j-1},\beta^\diamond_{j}}(\tau^\diamond_j)\Big|\mcF_t\bigg].
\end{align*}

Fix $(\veca,\vecb)\in\mcJ$. We will prove the Lemma in 3 different steps:\\

\textbf{Step 1} We first show continuity in $z$. Whenever $(t,\nu',z')\in\mcD_{\veca,\vecb}$ we have
\begin{align*}
Y^{t,\nu',z',\veca,\vecb,k+1}_t \geq \E\bigg[&\int_s^{T}\psi_{\alpha_r^{t,\nu',\veca,\vecb,u^\diamond}}(r,\theta^{t,\nu',z',\veca,\vecb,u^\diamond}_r)dr
+\Upsilon_{\alpha_T^{t,\nu',\veca,\vecb,u^\diamond}}(\theta^{t,\nu',z',\veca,\vecb,u^\diamond}_T)
-\sum_{j=1}^{N^\diamond} c_{\beta^\diamond_{j-1},\beta^\diamond_{j}}(\tau^\diamond_j)\Big|\mcF_t\bigg].
\end{align*}

Now, by definition we have $|\theta^{t,\nu,z',\veca,\vecb}_s-\theta^{t,\nu,z,\veca,\vecb}_s|\leq |z'-z|$, for all $s\in[0,T]$ which gives that \begin{align*}
\theta^{t,\nu,z',\veca,\vecb,(\tau_1^\diamond,\beta_1^\diamond)}_{\tau_1^\diamond} -\theta^{t,\nu,z,\veca,\vecb,(\tau_1^\diamond,\beta_1^\diamond)}_{\tau_1^\diamond}&=(z'+\int_t^{\tau_1^\diamond}A^{t,\nu,\veca,\vecb}_rdr)\wedge T(\beta_1^\diamond\wedge A^{t,\nu,\veca,\vecb}_{\tau_1^\diamond})
\\
&\quad-(z+\int_t^{\tau_1^\diamond}A^{t,\nu,\veca,\vecb}_rdr)\wedge T(\beta_1^\diamond \wedge A^{t,\nu,\veca,\vecb}_{\tau_1^\diamond}),
\end{align*}
and in particular $|\theta^{t,\nu,z',\veca,\vecb,(\tau_1^\diamond,\beta_1^\diamond)}_{s} -\theta^{t,\nu,z,\veca,\vecb,(\tau_1^\diamond,\beta_1^\diamond)}_{s}|\leq |z'-z|$, for all $s\in[0,T]$.
Repeating this argument and noting that $k$ is finite we conclude that
\begin{equation*}
|\theta^{t,\nu,z',\veca,\vecb,u^\diamond}_s-\theta^{t,\nu,z,\veca,\vecb,u^\diamond}_s|\leq |z'-z|,
\end{equation*}
for all $s\in[0,T]$. By the above we find that, $\Prob$-a.s.~for all $t\in[0,T]$,
\begin{align*}
Y^{t,\nu,z,\veca,\vecb,k+1}_t-Y^{t,\nu,(z+q)^+\wedge T,\veca,\vecb,k+1}_{t}\leq C|q|.
\end{align*}
Furthermore, since the reversed argument applies to $Y^{t,\nu,(z+q)^+\wedge T,\veca,\vecb,k+1}_{t}-Y^{t,\nu,z,\veca,\vecb,k+1}_t$ we find that
\begin{align*}
\sup_{(t,\nu,z)\in\mcD_{(\veca,\vecb)}}|Y^{t,\nu,(z+q)^+\wedge T,\veca,\vecb,k+1}_{t}-Y^{t,\nu,z,\veca,\vecb,k+1}_t|\leq C|q|,
\end{align*}
$\Prob$-a.s.~and
\begin{align*}
\E\bigg[\sup_{(t,\nu,z)\in\mcD_{(\veca,\vecb)}}|Y^{t,\nu,(z+q)^+\wedge T,\veca,\vecb,k}_{t}-Y^{t,\nu,z,\veca,\vecb,k}_{t}|^2\bigg]\leq C|q|,
\end{align*}
where we note that the constant $C$ does not depend on $k$.

\bigskip

%%%%%%%%%%%%%%%%%%%%%%%%%%%%%%%%%%%%%%%%%%%%%%%%%%%%%%%%%%%%%%%%%%%%%%%%%%%%%%%%%%%%%%%%%%%%%%%%%%%%%%%%%%%%%%%%%%%%%%%%%%%%%%%%%%%%%%%%%%%%%%%%%%

\textbf{Step 2} Next we show continuity in $\nu$. We have
\begin{align*}
Y^{t,\nu,z,\veca,\vecb,k+1}_{t}-Y^{t,\nu',z,\veca,\vecb,k+1}_{t}&\leq \E\bigg[\int_{t}^{\tau^\diamond_1 \wedge T}(\psi_{A^{t,\nu,\veca,\vecb}_s}(s,\theta_s^{t,\nu,z,\veca,\vecb})-\psi_{A^{t,\nu',\veca,\vecb}_s}(s,\theta_s^{t,\nu',z,\veca,\vecb}))ds
\\
&\quad+\ett_{[\tau^\diamond_1\geq T]}(\Upsilon_{A^{t,\nu,\veca,\vecb}_T}(\theta_T^{t,\nu,z,\veca,\vecb})-\Upsilon_{A^{t,\nu',\veca,\vecb}_T}(\theta_T^{t,\nu',z,\veca,\vecb}))
\\
&\quad+\ett_{[\tau^\diamond_1 < T]}(Y^{\tau^\diamond_1,\vecb\nu+\tau^\diamond_1(\beta^\diamond_1-\vecb)^+,\theta_{\tau^\diamond_1}^{t,\nu,z,\veca,\vecb}\wedge T\beta^\diamond_1,A^{t,\nu,\veca,\vecb}_{\tau^\diamond_1}\wedge\beta^\diamond_1,\beta^\diamond_1,k}_{\tau^\diamond_1}
\\
&\quad-Y^{\tau^\diamond_1,\vecb\nu'+\tau^\diamond_1(\beta^\diamond_1-\vecb)^+,\theta_{\tau^\diamond_1}^{t,\nu',z,\veca,\vecb}\wedge T\beta^\diamond_1,A^{t,\nu',\veca,\vecb}_{\tau^\diamond_1}\wedge\beta^\diamond_1,\beta^\diamond_1,k}_{\tau^\diamond_1})\Big| \mcF_t\bigg]
%\\
%&\leq \sup_{\tau\in\mcT_t}\E\bigg[\int_{t}^{\tau \wedge T}(\psi_{A^{t,\nu,\veca,\vecb}_s}(s,\theta_s^{t,\nu,z,\veca,\vecb})-\psi_{A^{t,\nu',\veca,\vecb}_s}(s,\theta_s^{t,\nu',z,\veca,\vecb}))ds
%\\
%&\quad+\ett_{[\tau\geq T]}(\Upsilon_{A^{t,\nu,\veca,\vecb}_T}(\theta_T^{t,\nu,z,\veca,\vecb})-\Upsilon_{A^{t,\nu',\veca,\vecb}_T}(\theta_T^{t,\nu',z,\veca,\vecb}))
%\\
%&\quad+\ett_{[\tau < T]}\max_{\beta\in\mcI^{-\vecb}}\{Y^{\tau,\vecb\nu+\tau(\beta-\vecb)^+,\theta_{\tau}^{t,\nu,z,\veca,\vecb}\wedge T\beta,A^{t,\nu,\veca,\vecb}_{\tau}\wedge\beta,\beta,k}_{\tau}
%\\
%&\quad-Y^{\tau,\vecb\nu'+\tau(\beta-\vecb)^+,\theta_{\tau}^{t,\nu',z,\veca,\vecb}\wedge T\beta,A^{t,\nu',\veca,\vecb}_{\tau}\wedge\beta,\beta,k}_{\tau}\}\Big| \mcF_t\bigg]
\end{align*}
Repeating the argument in the proof of Proposition~\ref{prop:H0} yields that $Y^{t,\nu,z,\veca,\vecb,k+1}_{t}-Y^{t,\nu',z,\veca,\vecb,k+1}_{t}\leq \Xi_t^{\nu,\nu',z,\veca,\vecb,k+1}$, where
\begin{align*}
\Xi^{\nu,\nu',z,\veca,\vecb,k+1}_t&:= \E\bigg[\int_{t}^{\tau^\diamond_1\wedge T}|e^{\int_{t}^s\lambda^{\nu,\vecb}_{\veca,\veca}(r)dr}-e^{\int_{t}^s\lambda^{\nu',\vecb}_{\veca,\veca}(r)dr}|\cdot |\psi_{\veca}(s,\theta_s^{t,z,\veca})|ds
\\
&+\ett_{[\tau^\diamond_1\geq T]}|e^{\int_{t}^T\lambda^{\nu,\vecb}_{\veca,\veca}(r)dr}-e^{\int_{t}^T\lambda^{\nu',\vecb}_{\veca,\veca}(r)dr}|\cdot |\Upsilon_{\veca}(\theta_T^{t,z,\veca})|
\\
&+\sum_{\veca'\in\mcA^{-\veca}_{\veca,\vecb}}\int_{t}^{\tau^\diamond_1\wedge T}\Big(|\lambda^{\nu,\vecb}_{\veca,\veca'}(s)e^{\int_{t}^s\lambda^{\nu,\vecb}_{\veca,\veca}(r)dr} -\lambda^{\nu',\vecb}_{\veca,\veca'}(s)e^{\int_{t}^s\lambda^{\nu',\vecb}_{\veca,\veca}(r)dr}|\bar Y_s
\\
&+K_{\lambda}\Xi^{\nu,\nu',\theta_s^{t,z,\veca}\wedge T\veca',\veca',\vecb,k+1}_s\Big)ds
\\
&+\ett_{[\tau^\diamond_1< T]}|e^{\int_{t}^{\tau^\diamond_1}\lambda^{\nu,\vecb}_{\veca,\veca}(r)dr} -e^{\int_{t}^{\tau^\diamond_1}\lambda^{\nu',\vecb}_{\veca,\veca}(r)dr}|\cdot \Xi^{\nu\beta^\diamond_1+{\tau^\diamond_1}(\beta^\diamond_1-\vecb)^+,\nu'\beta^\diamond_1+{\tau^\diamond_1}(\beta^\diamond_1-\vecb)^+ ,\theta_{\tau^\diamond_1}^{t,z,\veca}\wedge T\beta^\diamond_1 ,\veca\wedge\beta^\diamond_1,\beta^\diamond_1,k}_{\tau^\diamond_1}\Big| \mcF_{t}\bigg].
\end{align*}
Relying on an induction argument and the fact that the control $u^\diamond\in\mcU^f$ for all $k$, we conclude by symmetry that
\begin{align*}
|Y^{t,\nu,z,\veca,\vecb,k+1}_{t}-Y^{t,\nu',z,\veca,\vecb,k+1}_{t}|&\leq C|\nu-\nu'|\bar Y_t,
\end{align*}
for all $t\in[0,T]$, $\Prob$-a.s. Applying Doob's maximal inequality, we find that, for $p\in \R^{n}$
\begin{align*}
\E\bigg[\sup_{(t,\nu,z)\in\mcD_{(\veca,\vecb)}}|Y^{t,(\nu+p)^+\wedge T\vecb,z,\veca,\vecb,k+1}_{t}-Y^{t,\nu,z,\veca,\vecb,k+1}_{t}|^2\bigg] \leq C|p|,
\end{align*}
where again the constant $C$ does not depend on $k$. Put together, this implies that \textbf{H.k+1}.\ref{hyp:cont} holds.

%%%%%%%%%%%%%%%%%%%%%%%%%%%%%%%%%%%%%%%%%%%%%%%%%%%%%%%%%%%%%%%%%%%%%%%%%%%%%%%%%%%%%%%%%%%%%%%%%%%%%%%%%%%%%%%%%%%%%%%%%%%%%%%%%%%%%%%%%%%%%%%%%%

\bigskip

\textbf{Step 3} To show that \textbf{H.k+1}.\ref{hyp:inS2e} holds we note that
\begin{align*}
(Y^{t,\nu,z,\veca,\vecb,k+1}_s+\int_0^s\psi_{A_r^{t,\nu\veca,\vecb}}(r,\theta_s^{t,\nu,z,\veca,\vecb})dr:0\leq s\leq T)
\end{align*}
is the Snell envelope of a positive process in $\mcS^2_e$ and thus itself belongs to $\mcS^2_e$. Subtracting the continuous process $\int_0^\cdot \psi_{A_r^{t,\nu\veca,\vecb}}(r,\theta_s^{t,\nu,z,\veca,\vecb})dr$ we conclude that $Y^{t,\nu,z,\veca,\vecb,k+1}\in\mcS^2_e$.

Let $(t,\nu,z)\in \mcD_{(\veca,\vecb)}$ and assume that $t' \in [t,T]$. To evaluate $|Y^{t,\nu,z,\veca,\vecb,k+1}_{t}-Y^{t',\nu,z,\veca,\vecb,k+1}_{t'}|$ we note that during $[t,t']$ we may either have a transition in $A^{t,\nu,\veca,\vecb}$ (the probability of which is bounded by $1-e^{-K_\lambda(t'-t)}$), it may be optimal to switch to another mode or neither of the above. We conclude that
\begin{align}\nonumber
&|Y^{t,\nu,z,\veca,\vecb,k+1}_t-Y^{t',\nu,z,\veca,\vecb,k+1}_{t'}|\leq \E\Big[\int_t^{t'}\bar\psi(s)ds\big|\mcF_t\Big]+2(1-e^{-K_\lambda(t'-t)})\bar Y_t
\\
&+\sum_{\vecb'\in\mcI^{-\vecb}}\E\bigg[\sup_{s\in[t,t']}(|Y^{s,\vecb'\nu+s(\vecb'-\vecb)^+,z,\veca\wedge\vecb',\vecb',k}_{s} - Y^{t',\vecb'\nu+t'(\vecb'-\vecb)^+,z,\veca\wedge\vecb',\vecb',k}_{t'}|+|c_{\vecb,\vecb'}(s)-c_{\vecb,\vecb'}(t')|)\Big|\mcF_t\bigg]\nonumber
\\
&+ |\E\bigg[Y^{t',\nu,z,\veca,\vecb,k}_{t'}\Big|\mcF_t\bigg] -Y^{t',\nu,z,\veca,\vecb,k}_{t'}|+C|t'-t|,\label{ekv:Yttkbnd}
\end{align}
where the first three terms tend to 0 for all $|t'- t|\to 0$, $\Prob$-a.s.~(the first term by Assumption~\ref{ass:prelim}.\ref{ass:onpsi} and the third term by the induction hypothesis and continuity of the switching costs). Moving to the fourth term, we let $\tau^{r,\nu,z,\veca,\vecb,k+1}_s\in\mcT^{\bbF}_s$ be the first intervention time in an optimal control for $Y^{r,\nu,z,\veca,\vecb,k+1}_s$. By sublinearity of the operator $\mcK_{t,t'}$ introduced in the proof of Proposition~\ref{prop:H0} we have, for $(\zeta,\nu,z)\in\mcD_{(\veca,\vecb)}$ with $\zeta\leq t'$,
\begin{align*}
&|\E\Big[Y^{\zeta,\nu,z,\veca,\vecb,k+1}_{t'}\big|\mcF_t\Big] -Y^{\zeta,\nu,z,\veca,\vecb,k+1}_{t'}|
\\
&\leq \int_0^T\sup_{(r,p,q)\in \mcD_{(\veca,\vecb)}}\mcK_{t,t'}(\ett_{[s<\tau^{\zeta,\nu,z,\veca,\vecb,k+1}_{t'}]}e^{\int_r^s\lambda^{p,\vecb}_{\veca,\veca}(v)dv} \psi_{\veca}(s,q))ds
\\
&\quad+ \sup_{(r,p,q)\in \mcD_{(\veca,\vecb)}}\mcK_{t,t'}(\ett_{[T\leq\tau^{\zeta,\nu,z,\veca,\vecb,k+1}_{t'}]} e^{\int_r^T\lambda^{p,\vecb}_{\veca,\veca}(v)dv}\Upsilon_{\veca}(q))
\\
&\quad + \sum_{\veca'\in\mcA^{-\veca}_{\veca,\vecb}}\int_0^T\sup_{(r,p,q)\in \mcD_{(\veca,\vecb)}}\mcK_{t,t'} (\ett_{[s<\tau^{t',\nu,z,\veca,\vecb,k+1}_{t'}]}\lambda^{p,\vecb}_{\veca,\veca'}(s)e^{\int_r^s\lambda_{\veca,\veca}^{p,\vecb}(v)dv} Y^{s,p,q,\veca',\vecb,k+1}_{s})ds
\\
&\quad + \sum_{\vecb'\in\mcI^{-\vecb}} \sup_{(r,p,q)\in \mcD_{(\veca,\vecb)}}\mcK_{t,t'} (e^{\int_r^{\tau^{t',\nu,z,\veca,\vecb,k+1}_{t'}}\lambda_{\veca,\veca}^{p,\vecb}(v)dv} Y^{{\tau^{r,k+1}},p\vecb'+\tau(\vecb'-\vecb)^+,q,\veca',\vecb,k}_{{\tau^{r,k+1}}})+C|t'-\zeta|(\bar Y_t+\bar Y_{t'}),
\end{align*}
$\Prob$-a.s., where the last term represents the possibility of $A^{\zeta,\nu,z,\veca,\vecb}$ having a transition during $[\zeta,t')$.
Furthermore, arguing as in Step 1 and Step 2 we find that
\begin{align*}
|Y^{t',\nu,z,\veca,\vecb,k+1}_{t'}-Y^{\zeta,\nu',z',\veca,\vecb,k+1}_{t'}|\leq C(|\zeta-t'|(1+\bar Y_{t'})+|\nu'-\nu|+|z'-z|),
\end{align*}
$\Prob$-a.s. which implies that
\begin{align*}
\mcK_{t,t'}(Y^{t',\nu,z,\veca,\vecb,k+1}_{t'})\leq\mcK_{t,t'}(Y^{\zeta,\nu',z',\veca,\vecb,k+1}_{t'}) + C(|\zeta-t'|(1+\bar Y_{t}+\bar Y_{t'})+|\nu'-\nu|+|z'-z|),
\end{align*}
$\Prob$-a.s.~(where the exception set does not depend on $(t,t',\zeta,\nu,\nu',z,z')$). For each $M\geq 1$ we again partition the set $\mcD_{(\veca,\vecb)}$ into a partition with diameter $\delta(M)$ (now into a rectangular partition with a constant step-size $\Delta t$ in the time variable $t$) and note that when $t'\in[t_l^M,t_{l}^M+\Delta t]$ then $|\tau^{t_l^M,\nu,z,\veca,\vecb,k+1}_{t'}-\tau^{t_l^M,\nu,z,\veca,\vecb,k+1}_{t_l^M}|\wedge |\tau^{t_l^M,\nu,z,\veca,\vecb,k+1}_{t'}-\tau^{t_l^M,\nu,z,\veca,\vecb,k+1}_{t_{l}^M+\Delta t}|\leq \Delta t$, $\Prob$-a.s. With $\tau^{M,k+1}_{l,1}:=\tau^{t_l^M,\nu_l^M,z_l^M,\veca,\vecb,k+1}_{t_l^M}$ and $\tau^{M,k+1}_{l,2}:=\tau^{t_l^M,\nu_l^M,z_l^M,\veca,\vecb,k+1}_{t_{l}^M+\Delta t}$ we have
\begin{align*}
&|\E\Big[Y^{t',\nu,z,\veca,\vecb,k+1}_{t'}\big|\mcF_t\Big] -Y^{t',\nu,z,\veca,\vecb,k+1}_{t'}|
\\
%%%%%%%%%%%%%%%%%%%%%%%%%%%%%%%%%%%%%%%%%%%%%%%%%%%%%%%%%%%%%%%%%%%%%%%%%%%%%%%%%%%%%%%%%%%%%%%%%%%%%%%%%%%%%%%%%%%%%%%%%%%%%%%%%%%%%%%%%%%%%%%%%%
&\leq \sum_{i=1}^2\sum_{l=1}^M\bigg(\int_0^T\mcK_{t,t'}(\ett_{[s<\tau^{M,k+1}_{l,i}]}e^{\int_{t^M_l}^s\lambda^{\nu_l^M,\vecb}_{\veca,\veca}(v)dv} \psi_{\veca}(s,z^M_l))ds+ \mcK_{t,t'}(\ett_{[T<\tau^{M,k+1}_{l,i}]}e^{\int_{t^M_l}^T\lambda^{\nu_l^M,\vecb}_{\veca,\veca}(v)dv}\Upsilon_{\veca}(z^M_l))
\\
&\quad + \sum_{\veca'\in\mcA^{-\veca}_{\veca,\vecb}} \int_0^T\mcK_{t,t'} (\ett_{[s<\tau^{M,k+1}_{l,i}]}\lambda^{\nu_l^M,\vecb}_{\veca,\veca'}(s)e^{\int_{t^M_l}^{s}\lambda^{\nu_l^M,\vecb}_{\veca,\veca}(v)dv} Y^{s,\nu_l^M,z^M_l,\veca',\vecb,k+1}_{s})ds
\\
&\quad + \sum_{\vecb'\in\mcI^{-\vecb}}\mcK_{t,t'} (e^{\int_{t^M_l}^{\tau^{M,k+1}_{l,i}}\lambda^{\nu_l^M,\vecb}_{\veca,\veca}(v)dv} Y^{\tau^{M,k+1}_{l,i},\nu_l^M\vecb'+\tau^{M,k+1}_{l,i}(\vecb'-\vecb)^+,z^M_l,\veca,\vecb,k}_{\tau^{M,k+1}_{l,i}})\bigg)+C(1+\bar Y_t+\bar Y_{t'})\delta(M)
\\
&\quad + C \sum_{\veca'\in\mcA^{-\veca}_{\veca,\vecb}} \sup_{r\in[0,T]}\E\big[\max_{l\in\{1,\ldots,M\}}\sup_{(s,p,q)\in G_l^M}|Y^{s,p,q,\veca',\vecb,k+1}_{s} - Y^{s,\nu^M_l,z^M_l,\veca',\vecb,k+1}_{s}|\big|\mcF_r\big]
\\
&\quad + C \sum_{\vecb'\in\mcI^{-\vecb}} \sup_{r\in[0,T]}\E\big[\max_{l\in\{1,\ldots,M\}}\sup_{(s,p,q)\in G_l^M}|Y^{s,p,q,\veca'\wedge\vecb',\vecb',k}_{s} - Y^{t^M_l,\nu^M_l,z^M_l,\veca'\wedge\vecb',\vecb',k}_{t_l}|\big|\mcF_r\big].
\end{align*}
where we have used that
\begin{align*}
&\mcK_{t,t'}(e^{\int_{r}^s\lambda^{\nu,\vecb}_{\veca,\veca}(v)dv}Y^{s,\nu,z,\veca',\vecb,k}_{s})
\\
%%%%%%%%%%%%%%%%%%%%%%%%%%%%%%%%%%%%%%%%%%%%%%%%%%%%%%%%%%%%%%%%%%%%%%%%%%%%%%%%%%%%%%%%%%%%%%%%%%%%%%%%%%%%%%%%%%%%%%%%%%%%%%%%%%%%%%%%%%%%%%%%%%
&\leq \mcK_{t,t'}(e^{\int_{r'}^{s'}\lambda^{\nu',\vecb}_{\veca,\veca}(v)dv}Y^{s',\nu',z',\veca',\vecb,k}_{s'}) + (K_\lambda (|r'-r|+|s'-s|) + Tk_\lambda |\nu'-\nu|) (\bar Y_t + \bar Y_{t'})
\\
&\quad + \E\big[|Y^{s,\nu,z,\veca',\vecb,k}_{s} - Y^{s',\nu',z',\veca',\vecb,k}_{s'}|\big|\mcF_t\big] +\E\big[|Y^{s,\nu,z,\veca',\vecb,k}_{s} -Y^{s',\nu',z',\veca',\vecb,k}_{s}|\big|\mcF_{t'}\big].
\end{align*}

In the above inequality for $|\E\Big[Y^{t',\nu,z,\veca,\vecb,k+1}_{t'}\big|\mcF_t\Big] -Y^{t',\nu,z,\veca,\vecb,k+1}_{t'}|$ each of the terms of type $K_{t,t'}$ can be represented by a integral of type $\int_{t}^{t'}Z_rdW_r$. Now since $M$ was arbitrary and $\delta(M)\to 0$ as $M\to\infty$ we conclude from \eqref{ekv:Yttkbnd} that $(Y^{t,\nu,z,\veca,\vecb}_t:0\leq t\leq T)$ has a version such that
\begin{align*}
\sup_{(t,\nu,z)\in \mcD_{(\veca,\vecb)}} |Y^{(t+h)^+\wedge T,\nu,z,\veca,\vecb,k+1}_{(t+h)^+\wedge T}-Y^{t,\nu,z,\veca,\vecb,k+1}_t|\to 0,
\end{align*}
$\Prob$-a.s.~as $h\to 0$ for all $(\veca,\vecb)\in\mcJ$.

By an induction argument we conclude that \textbf{H.k} holds for all $k\geq 0$.\qed\\

\begin{rem}
For the more general case when $\bbF$ is generated by a Brownian motion and an independent Poisson random measure the induction hypothesis \textbf{H.k}.\ref{hyp:inS2e} has to be weakened to $\big(Y^{t,\nu,z,\veca,\vecb,k}_t: 0\leq t\leq T\big)$ being right continuous and LCE in $t$. Continuity in $(\nu,z)$ follows immediately from the above proof. The main difference is that the martingale representation theorem cannot be applied to show continuity in $t$. However, the terms $\mcK_{t,t+h}(\cdots)$ goes to zero $\Prob$-a.s.~as $h\searrow 0$ by right continuity of the filtration and right continuity of $\big(Y^{t,\nu,z,\veca,\vecb,k}_t: 0\leq t\leq T\big)$ follows. Furthermore, left continuity in expectation is immediate from~\eqref{ekv:Yttkbnd}.
\end{rem}

We are now ready to show that the limit family, $\lim_{k\to\infty}((Y^{t,\nu,z,\veca,\vecb,k}_s)_{0\leq s\leq T}: ({t,\nu,z,\veca,\vecb,k})\in \bigSET)$, exists and satisfies the properties of a verification family. We start with existence:
\bigskip

\begin{prop}\label{prop:Yklim}
For each $(t,\nu,z,\veca,\vecb)\in \bigSET$, the limit $\tilde Y^{t,\nu,z,\veca,\vecb}:=\lim_{k\to\infty}Y^{t,\nu,z,\veca,\vecb,k}$, exists as an increasing pointwise limit, $\Prob$-a.s. Furthermore, the process $(\tilde Y^{t,\nu,z,\veca,\vecb}_t:0\leq t\leq T)$ is continuous.
\end{prop}

\noindent\emph{Proof.} Since $\mcU^k_t\subset \mcU^{k+1}_t$ we have by \eqref{ekv:Ykalt} that, $\Prob$-a.s.,
\begin{align}
Y^{t,\nu,z,\veca,\vecb,k}_s \leq Y^{t,\nu,z,\veca,\vecb,k+1}_s\leq \E\bigg[\int_0^T \bar\psi(r)dr + \bar\Upsilon \Big|\mcG_s\bigg],\label{ekv:Ykbnd}
\end{align}
where the right hand side is bounded in $L^2$. Hence, the sequence $((Y^{t,\nu,z,\veca,\vecb,k}_s)_{0\leq s\leq T}: (t,\nu,z,\veca,\vecb)\in \bigSET)$ converges $\Prob$--a.s.~for all $s\in [0,T]$.\\

Concerning the second claim, note that by Proposition~\ref{prop:YkL2bnd} there is for each $\delta>0$ and $p\in(1,2)$ a constant $K>0$, such that the set \begin{align*}
B:=\Big\{\omega\in\Omega:\sup_{t\in[0,T]}\max\Big(\bar Y_t(\omega),\E\Big[\int_0^T\bar\psi^p(r)dr+\bar\Upsilon^p \big|\mcF_t\Big](\omega)\Big)\leq K\Big\}
\end{align*}
has probability $\Prob(B)\geq 1-\delta$. By the ``no-free-loop'' condition (Assumption~\ref{ass:prelim}.(\ref{ass:onc})) and the finiteness of $\mcI$ we get that for any control $(\tau_1,\ldots,\tau_N;\beta_1,\ldots,\beta_N)$,
\begin{align*}
\sum_{j=1}^{N}c_{\beta_j,\beta_{j-1}}(\tau_j)\geq \epsilon (N-m)/m,
\end{align*}
$\Prob$-a.s. Hence, there is a $\Prob$-null set $\mcN\subset\Omega$ such that for all $\omega\in B\setminus\mcN$,
\begin{align*}
0\leq Y^{t,\nu,z,\veca,\vecb,k}_t(\omega)&\leq \bar Y_t(\omega)-\epsilon(\E[N^k|\mcF_t](\omega)/m-1)
\\
&\leq K+\epsilon-\epsilon/m\E[ N^k |\mcF_t](\omega),
\end{align*}
where $(\tau^k_1,\ldots,\tau^k_{N^k};\beta^k_1,\ldots,\beta^k_{N^k})$ is a control corresponding to $Y^{t,\nu,z,\veca,\vecb,k}_t$. This implies that for $k'>0$ we have,
\begin{align*}
\Prob[N^k>k' |\mcF_t](\omega)\leq (Km/\epsilon+m)/k'.
\end{align*}
Now, for $0\leq k'\leq k$ we let $u^k_{k'}$ be the truncated control $u^k_{k'}:=(\tau^k_1,\ldots,\tau^k_{N^k\wedge k'};\beta^k_1,\ldots,\beta^k_{N^k\wedge k'})$. Then, clearly $u^k_{k'}\in\mcU^{k'}$ and we have
\begin{align*}
&\E\bigg[\int_t^{T}\psi_{\alpha_r^{t,\nu,\veca,\vecb,u^k_{k'}}}(r,\theta^{t,\nu,z,\veca,\vecb,u^k_{k'}}_r)dr
+\Upsilon_{\alpha_T^{t,\nu,\veca,\vecb,u^k_{k'}}}(\theta^{t,\nu,z,\veca,\vecb,u^k_{k'}}_T)
-\sum_{j=1}^{N^k\wedge k'} c_{\beta^k_{j-1},\beta^k_{j}}(\tau^k_j)\Big|\mcF_t\bigg]
\\
&\leq Y^{t,\nu,z,\veca,\vecb,k'}_t\leq Y^{t,\nu,z,\veca,\vecb,k}_t.
\end{align*}
Furthermore,
\begin{align*}
&Y^{t,\nu,z,\veca,\vecb,k}_t(\omega)-\E\bigg[\int_t^{T}\psi_{\alpha_r^{t,\nu,\veca,\vecb,u^k_{k'}}}(r,\theta^{t,\nu,z,\veca,\vecb,u^k_{k'}}_r)dr
+\Upsilon_{\alpha_T^{t,\nu,\veca,\vecb,u^k_{k'}}}(\theta^{t,\nu,z,\veca,\vecb,u^k_{k'}}_T)
-\sum_{j=1}^{N^k\wedge k'} c_{\beta^k_{j-1},\beta^k_{j}}(\tau^k_j)\Big|\mcF_t\bigg](\omega)
\\
&\leq \E\Big[\ett_{[N^k>k']}\big(\int_0^T\bar\psi(r)dr+\bar\Upsilon\big)\big|\mcF_t\Big](\omega)
\\
&\leq 2\E\big[\ett_{[N^k>k']}|\mcF_t\big]^{1/q}(\omega)\E\Big[\int_0^T\bar\psi^p(r)dr+\bar\Upsilon^p\big|\mcF_t\Big]^{1/p}(\omega)
\\
&\leq 2K^{1/p}((Km/\epsilon+m)/k')^{1/q},
\end{align*}
where $\frac{1}{p}+\frac{1}{q}=1$. There is thus a constant $C>0$ such that
\begin{align*}
Y^{t,\nu,z,\veca,\vecb,k}_t(\omega)-Y^{t,\nu,z,\veca,\vecb,k'}_t(\omega)\leq C(k')^{-1/q},
\end{align*}
for all $t\in[0,T]$ and all $0\leq k'\leq k$. We conclude that for all $\omega\in B\setminus\mcN$, the sequence $(Y^{t,\nu,z,\veca,\vecb,k}_t(\omega):0\leq t\leq T)_{k\geq 0}$ is a sequence of continuous functions that converges uniformly in $t$ which implies that the limit is continuous. Since $\delta>0$ was arbitrary we conclude that $\Prob$-almost all trajectories of $(\tilde Y^{t,\nu,z,\veca,\vecb,k}_t:0\leq t\leq T)$ are continuous.\qed

\bigskip

\begin{thm}\label{thm:limVT}
The limit family $((\tilde Y^{t,\nu,z,\veca,\vecb}_s)_{0\leq s\leq T}: (t,\nu,z,\veca,\vecb)\in \bigSET)$ is a verification family.
\end{thm}

\noindent \emph{Proof.} As noted above, property \emph{\ref{vFAM:bnd})} of a verification family follows immediately by Proposition~\ref{prop:YkL2bnd}. We now show that the limit satisfies the additional properties of the verification theorem as well, starting with the recursion.\\

\noindent\emph{\ref{vFAM:recur}) Limit satisfies \eqref{ekv:Ydef}.} From Proposition~\ref{prop:Yklim} and the proof of Proposition~\ref{prop:Hk} it follows that $Y^{t,\nu,z,\veca,\vecb}_t$ is $\Prob$-a.s.~jointly continuous in $(t,\nu,z)$ and in particular we note that $(Y^{s,\beta\nu+s(\beta-\vecb)^ + ,\theta_s^{t,z,\veca,\vecb}\wedge T\beta,A^{t,\nu,\veca,\vecb}_s\wedge\beta,\beta,k}_s:0\leq s\leq T)$
%\begin{align*}
%(Y^{s,\beta\nu+s(\beta-\vecb)^ + ,\theta_s^{t,z,\veca,\vecb}\wedge T\beta,A^{t,\nu,\veca,\vecb}_s\wedge\beta,\beta,k}_s:0\leq s\leq T)
%\end{align*}
is a sequence of \cadlag processes that converges $\Prob$-a.s.~pointwisely to the \cadlag process $(\tilde Y^{s,\beta\nu+s(\beta-\vecb)^ + ,\theta_s^{t,z,\veca,\vecb}\wedge T\beta,A^{t,\nu,\veca,\vecb}_s\wedge\beta,\beta}_s:0\leq s\leq T)$.
%\begin{align*}
%(\tilde Y^{s,\beta\nu+s(\beta-\vecb)^ + ,\theta_s^{t,z,\veca,\vecb}\wedge T\beta,A^{t,\nu,\veca,\vecb}_s\wedge\beta,\beta}_s:0\leq s\leq T).
%\end{align*}
We can thus use (\ref{Snell:lim}) of Theorem~\ref{thm:Snell} and find that
\begin{align*}
\tilde Y^{t,\nu,z,\veca,\vecb}_s=\esssup_{\tau \in \mcT_{s}} \E\bigg[&\int_s^{\tau\wedge T}\psi_{A^{t,\nu,\veca,\vecb}_r}(r,\theta_r^{t,\nu,z,\veca,\vecb})dr+\ett_{[\tau \geq T]}\Upsilon_{A^{t,\nu,\veca,\vecb}_r}(\theta_T^{t,\nu,z,\veca,\vecb})
\\
&+\ett_{[\tau < T]}\max_{\beta\in\mcI}\left\{-c_{\vecb,\beta}(\tau)+\tilde Y^{\tau,\beta\nu+\tau(\beta-\vecb)^ + ,\theta_\tau^{t,z,\veca,\vecb}\wedge T\beta,A^{t,\nu,\veca,\vecb}_\tau\wedge\beta,\beta}_\tau\right\}\Big| \mcG_s\bigg].
\end{align*}

\noindent\emph{\ref{vFAM:inS2e}) Limit in $\mcS^2_e$.} As $\tilde Y^{t,\nu,z,\veca,\vecb}_s+\int_0^s\psi_{A^{t,\nu,\veca,\vecb}_r}(r,\theta^{t,\nu,z,\veca,\vecb}_r)dr$ is the limit of an increasing sequence of \cadlag supermartingales it is also a \cadlag supermartingale (see \eg\cite{KarShreve1}). It remains to show that the limit is LCE. Rather than appealing to a uniform convergence argument as in the proof of Proposition~\ref{prop:Yklim} we give a direct, more intuitive proof, along the lines of \cite{BollanMSwitch1} and \cite{MartyrSigned}. We will look for a contradiction and let $(\gamma_j)_{j\geq 1}$ be a sequence of $\bbG$-stopping times such that $\gamma_j\nearrow\gamma\in\mcT$ and assume that
\begin{align*}
\lim_{j\to\infty} \E\left[\tilde Y^{t,\nu,z,\veca,\vecb}_{\gamma_j}\right]>\E\left[\tilde Y^{t,\nu,z,\veca,\vecb}_{\gamma}\right].
\end{align*}
Then, the previous step and the Doob-Meyer decomposition of the Snell envelope implies that
\begin{align*}
\lim_{j\to\infty}\tilde Y^{t,\nu,z,\veca,\vecb}_{\gamma_j}=\lim_{j\to\infty}\max_{\beta\in\mcI^{-b_n}}\left\{-c_{b_n,\beta}(\gamma_j)+\tilde Y^{\gamma_j,\nu\beta+\gamma_j(\beta-\vecb)^+,\theta^{t,z,\veca,\vecb}_{\gamma_j}\wedge T\beta,A^{t,\nu,\veca,\vecb}_{\gamma_j}\wedge\beta,\beta}_{\gamma_j}\right\}
\end{align*}
on some measurable set $M_1\subset\Omega$ with $\Prob(M_1)>0$ and
\begin{align}\nonumber
&\lim_{j\to\infty} \E\left[\ett_{M_1}\max_{\beta\in\mcI^{-b_n}}\left\{-c_{b_n,\beta}(\gamma_j)+\tilde Y^{\gamma_j,\nu\beta+\gamma_j(\beta-\vecb)^+,\theta^{t,z,\veca,\vecb}_{\gamma_j}\wedge T\beta,A^{t,\nu,\veca,\vecb}_{\gamma_j}\wedge\beta,\beta}_{\gamma_j}\right\}\right]
\\
&>\E\left[\ett_{M_1}\max_{\beta\in\mcI^{-b_n}}\left\{-c_{b_n,\beta}(\gamma)+\tilde Y^{\gamma,\nu\beta+\gamma(\beta-\vecb)^+,\theta^{t,z,\veca,\vecb}_{\gamma}\wedge T\beta,A^{t,\nu,\veca,\vecb}_{\gamma_j}\wedge\beta,\beta}_{\gamma}\right\}\right].\label{ekv:disk}
\end{align}
Since $c$ and $\tilde Y$ are both left limited there is a $\mcF_{\gamma}$-measurable random variable $\hat \beta_1\in \mcI^{-b_n}$ such that
\begin{align*}
&\lim_{j\to\infty}\max_{\beta\in\mcI^{-b_n}}\left\{-c_{b_n,\beta}(\gamma_j)+\tilde Y^{\gamma_j,\nu\beta+\gamma_j(\beta-\vecb)^+,\theta^{t,z,\veca,\vecb}_{\gamma_j}\wedge T\beta,A^{t,\nu,\veca,\vecb}_{\gamma_j}\wedge\beta,\beta}_{\gamma_j}\right\}
\\
&=\lim_{j\to\infty}\left\{-c_{b_n,\hat \beta_1}(\gamma_j)+\tilde Y^{\gamma_j,\nu\hat \beta_1+\gamma_j(\hat \beta_1-\vecb)^+,\theta^{t,z,\veca,\vecb}_{\gamma_j}\wedge T\hat \beta_1,A^{t,\nu,\veca,\vecb}_{\gamma_j}\wedge\hat \beta_1,\hat \beta_1}_{\gamma_j}\right\}.
\end{align*}
Now, by \eqref{ekv:disk} and continuity of the switching costs we must have
\begin{align*}
\lim_{j\to\infty} \E\left[\ett_{M_1}\tilde Y^{\gamma_j,\nu\hat \beta_1+\gamma_j(\hat \beta_1-\vecb)^+,\theta^{t,z,\veca,\vecb}_{\gamma_j}\wedge T\hat \beta_1,A^{t,\nu,\veca,\vecb}_{\gamma_j}\wedge\hat \beta_1,\hat \beta_1}_{\gamma_j}\right]>\E\left[\ett_{M_1}\tilde Y^{\gamma,\nu\hat \beta_1+\gamma(\hat \beta_1-\vecb)^+,\theta^{t,z,\veca,\vecb}_{\gamma}\wedge T\hat \beta_1,A^{t,\nu,\veca,\vecb}_{\gamma}\wedge\hat \beta_1,\hat \beta_1}_{\gamma}\right].
\end{align*}
Repeating this argument $l$ times we find that there is a sequence of measurable sets $M_1\supset M_2 \supset \cdots\supset M_l$ with $\Prob(M_l)>0$ and a sequence $\hat\beta_2,\hat\beta_3,\ldots,\hat\beta_l$ of $\mcF_{\gamma}$-measurable random variables such that $\hat\beta_{i+1}\in\mcI^{\hat\beta_i}$ for $i=1,\ldots,l-1$ and
\begin{align*}
&\lim_{j\to\infty}\tilde Y^{\gamma_j,\nu\underline{\hat \beta}_i+\gamma_j(\hat \beta_
i-\underline{\hat \beta}_i)^+,\theta^{t,z,\veca,\vecb}_{\gamma_j}\wedge T\underline{\hat \beta}_i,A^{t,\nu,\veca,\vecb}_{\gamma_j}\wedge \underline{\hat \beta}_i,\hat\beta_i}_{\gamma_j}
\\
&=\lim_{j\to\infty}\left\{-c_{\hat \beta_{i},\hat\beta_{i+1}}(\gamma_j)+\tilde Y^{\gamma_j,\nu\underline{\hat \beta}_{i+1}+\gamma_j(\hat \beta_
{i+1}-\underline{\hat \beta}_{i+1})^+,\theta^{t,z,\veca,\vecb}_{\gamma_j}\wedge T\underline{\hat \beta}_{i+1},A^{t,\nu,\veca,\vecb}_{\gamma_j}\wedge \underline{\hat \beta}_{i+1},\hat\beta_{i+1}}\right\},
\end{align*}
on $M_{i+1}$, with $\underline{\hat \beta}_{i}:=\hat\beta_1\wedge\hat\beta_2\wedge\cdots\wedge\hat\beta_i$. Now, since $\mcI$ is finite we can always find (possibly random) $l'$ and $l$ such that $l<l'$, $\beta_{l}=\beta_{l'}$ and $\underline{\hat \beta}_{l}=\underline{\hat \beta}_{l'}$ implying that on $M_{l'}$ we have
\begin{align*}
\lim_{j\to\infty}\left\{c_{\hat \beta_{l},\hat\beta_{l+1}}(\gamma_j)+\cdots+c_{\hat \beta_{l'-1},\hat\beta_{l'}}(\gamma_j)\right\}= 0
\end{align*}
contradicting the ``no free loop''-condition of Assumption~\ref{ass:prelim}.\ref{ass:onc}.\\

\noindent\emph{\ref{vFAM:cont}) Limit family continuous in $(\nu,z)$.}
In the proof of Proposition~\ref{prop:Hk} we showed that there is a $C>0$ such that for all $(p,q)\in \R^{2n}$
\begin{align*}
\E\bigg[\sup_{(t,\nu,z)\in\mcD_{(\veca,\vecb)}}|Y^{t,(\nu+p)^+\wedge T\vecb,(z+q)^+\wedge T,\veca,\vecb,k}_{t}-Y^{t,\nu,z,\veca,\vecb,k}_{t}|^2\bigg] \leq C(|p|+|q|),
\end{align*}
for all $k$. Letting $k\to\infty$ it follows that the limit family satisfies the same inequality.

This finishes the proof.\qed\\

\begin{rem}
Note that the results of this section naturally generalize to the case when $\bbF$ is a more general filtration, \eg generated by a Brownian motion and an independent Poisson random measure.
\end{rem}

We can now apply Theorem~\ref{thm:Snell}.\ref{Snell:DynP} to the recursion \eqref{ekv:Ydef} and get
\begin{align}\nonumber
\tilde Y^{t,\nu,z,\veca,\vecb}_t=\esssup_{\tau \in \mcT_{t}^{\bbF}} \E\bigg[&\int_t^{\tau\wedge T}e^{\int_t^r\lambda^{\nu,\vecb}_{\veca,\veca}(v)dv}\psi_{\veca}(r,\theta_r^{t,z,\veca})dr+\ett_{[\tau \geq T]}e^{\int_t^T\lambda^{\nu,\vecb}_{\veca,\veca}(v)dv}\Upsilon_{\veca}(\theta_T^{t,z,\veca})+
\\
&+\sum_{\veca'\in\mcA_{\veca,\vecb}^{-\veca}}\int_{t}^{\tau\wedge T}\lambda^{\nu,\vecb}_{\veca,\veca'}(r)e^{\int_t^r\lambda^{\nu,\vecb}_{\veca,\veca}(v)dv} Y^{r,\nu,\theta_r^{t,z,\veca}\wedge T\veca',\veca',\vecb,k}_r dr\nonumber
\\
&+\ett_{[\tau < T]}e^{\int_t^\tau\lambda^{\nu,\vecb}_{\veca,\veca}(v)dv}\max_{\beta\in\mcI}\left\{-c_{\vecb,\beta}(\tau)+\tilde Y^{\tau,\beta\nu+\tau(\beta-\vecb)^+,\theta_\tau^{t,z,\veca}\wedge T\beta,\veca\wedge\beta,\beta}_\tau\right\}\Big| \mcF_t\bigg].\label{ekv:YdefALT}
\end{align}

%%%%%%%%%%%%%%%%%%%%%%%%%%%%%%%%%%%%%%%%%%%%%%%%%%%%%%%%%%%%%%%%%%%%%%%%%%%%%%%%%%%%%%%%%%%%%%%%%%%%%%%%%%%%%%%%%%%%%%%%%%%%%%%%%%%%%%%%%%%%%%%%%%
%%%%%%%%%%%%%%%%%%%%%%%%%%%%%%%%%%%%%%%%%%%%%%%%%%%%%%%%%%%%%%%%%%%%%%%%%%%%%%%%%%%%%%%%%%%%%%%%%%%%%%%%%%%%%%%%%%%%%%%%%%%%%%%%%%%%%%%%%%%%%%%%%%
%%%%%%%%%%%%%%%%%%%%%%%%%%%%%%%%%%%%%%%%%%%%%%%%%%%%%%%%%%%%%%%%%%%%%%%%%%%%%%%%%%%%%%%%%%%%%%%%%%%%%%%%%%%%%%%%%%%%%%%%%%%%%%%%%%%%%%%%%%%%%%%%%%

\section{Value function representation\label{sec:ValFun}}
We now turn to the case when uncertainties in $\psi$ and $\Upsilon$ are modeled by a stochastic differential equation (SDE) as follows
\begin{align*}
dX_t&=a(t,X_t)dt+\sigma(t,X_t)dB_t,\quad t\in[0,T],\\
X_0&=x_0,
\end{align*}
where $a:[0,T]\times \R^m \to \R^m$ and $\sigma:[0,T]\times \R^m \to \R^{m\times m}$ are two deterministic, continuous functions that satisfy
\begin{equation*}
|a(t,x)|+|\sigma(t,x)|\leq C(1+|x|)
\end{equation*}
and
\begin{equation*}
|a(t,x)-a(t,x')|+|\sigma(t,x)-\sigma(t,x')|\leq C|x-x'|.
\end{equation*}
To be able to consider feedback-control formulations we will, for all $t\in [0,T]$ and $x\in\R^m$, define the process $(X_s^{t,x};0\leq s\leq T)$ as the strong solution to
\begin{align*}
dX_s^{t,x}&=a(s,X^{t,x}_s)ds+\sigma(s,X^{t,x}_s)dW_s,\quad \forall s\in(t,T],\\
X^{t,x}_s&=x,\quad \forall s\in[0,t].
\end{align*}
A standard result (see \eg Theorem 6.16 in Ch. 1 of \cite{XYZbook}) is that, for any $p\geq 1$, there exists a constant $C>0$ such that
\begin{equation}\label{ekv:Xbound1}
\E\bigg[\sup_{s\in [0,T]}|X^{t,x}_s|^p\bigg] \leq C (1+|x|^{p})
\end{equation}
and for all $t'\in [0,T]$ and all $x'\in\R$
\begin{align}
&\E\bigg[\sup_{s\in [0,T]}|X^{t,x}_s-X^{t',x'}_s|^p\bigg]\leq C (1+|x|^p)(|x-x'|^p+|t'-t|^{p/2}).\label{ekv:Xbound2}
\end{align}
We will consider the problem of finding a feedback strategy $u\in\mcU_t$ that, for each $(t,x,\nu,z,\veca,\vecb)\in\bar\mcD:=[0,T]\times\R^m\times\cup_{(\veca,\vecb)\in\mcJ}([0,T]^{\vecb}\times [0,T]^{\veca^+}\times (\veca,\vecb))$, maximizes
\begin{align}\nonumber
J_{(\veca,\vecb)}(t,x,\nu,z;u)&:=\E\bigg[\int_t^T\varphi_{\alpha^{t,x,\nu,\veca,\vecb,u}_s}(s,X_s^{t,x},\theta_s^{t,x,\nu,z,\veca,\vecb,u})ds
\\
&\quad + h_{\alpha^{t,x,\nu,\veca,\vecb,u}_T}(X^{t,x}_T,\theta_T^{t,x,\nu,z,\veca,\vecb,u}) -\sum_{j}c_{\beta_{j-1},\beta_j}(\tau_j)\bigg],\label{ekv:objfun}
\end{align}
where for each $\veca\in\mcA$ the deterministic function $\varphi_\veca:[0,T]\times\R^m\times[0,T]^{\veca^+}\to\R_+$ is locally Lipschitz and of polynomial growth in $x$, Lipschitz in $z$ and such that $\varphi_\veca(\cdot,x,z)$ is a \cadlag function and the deterministic function $h_{\veca}:\R^m \times[0,T]^{\veca^+}\to\R_+$ is locally Lipschitz in $x$ and Lipschitz in $z$.

Furthermore, we assume that the switching costs are deterministic and that, for each $(\veca,\vecb)\in\mcJ$ and $\nu\in[0,T]^{\vecb}$ and each $\veca'\in\mcA_{\veca,\vecb}$, the transition rates
\begin{equation}
\lambda^{\nu,\vecb}_{\veca,\veca'}(r):=\rho^{\nu,\vecb}_{\veca,\veca'}(r,X^{t,x}_r),\label{ekv:transINT}
\end{equation}
where $\rho^{\nu,\vecb}_{\veca,\veca'}:[0,T]\times \R^m\to\R$ is Lipschitz in $\nu$, bounded and locally Lipschitz in $x$ and $\rho^{\nu,\vecb}_{\veca,\veca'}(\cdot,x)$ is \cadlagp

We note that using the results in Sections~\ref{sec:VERthm} and \ref{sec:exist} there exists a control $u^*\in\mcU_t$ such that
\begin{equation*}
J_{(\veca,\vecb)}(t,x,\nu,z;u^*)\geq J_{(\veca,\vecb)}(t,x,\nu,z;u),
\end{equation*}
for all $u\in\mcU_t$ and more generally a unique family of processes $((Y^{t,x,\nu,z,\veca,\vecb}_s)_{0\leq s\leq T}:(t,x,\nu,z,\veca,\vecb)\in\bar\mcD)$ such that
\begin{align*}
Y^{t,x,\nu,z,\veca,\vecb}_s&=\esssup_{u\in \mcU_s} \E\bigg[\int_s^{T}\varphi_{\alpha_r^{t,x,\nu,\veca,\vecb,u}}(r,X^{t,x}_r,\theta^{t,x,\nu,z,\veca,\vecb,u}_r)dr
\\
&\quad +h_{\alpha_T^{t,x,\nu,\veca,\vecb,u}}(X^{t,x}_T,\theta^{t,x,\nu,z,\veca,\vecb,u}_T)
-\sum_{j=1}^N c_{\beta_{j-1},\beta_{j}}(\tau_j)\Big|\mcG_s\bigg],
\end{align*}
where $(\alpha_s^{t,x,\nu,\veca,\vecb,u}:0\leq s\leq T)$ is the process $(\alpha_s^{t,\nu,\veca,\vecb,u}:0\leq s\leq T)$ with transition intensities given by \eqref{ekv:transINT}. Furthermore, the following estimates hold:

\begin{prop}\label{prop:VFbnd}
There exists $K_Y>0$ such that, for each $(\veca,\vecb) \in \mcJ$, we have
\begin{equation*}
\E\left[\sup_{s\in [0,T]}|Y^{t,x,\nu,z,\veca,\vecb}_s|^2\right]\leq K_Y(1+|x|^{2p}), \quad\forall \,(t,x,\nu,z)\in \bar\mcD_{(\veca,\vecb)}.
\end{equation*}
Furthermore, $Y^{t,x,\nu,z,\veca,\vecb}_t$ is (indistinguishable from) a deterministic process and
\begin{equation*}
|Y^{t,x,\nu,z,\veca,\vecb}_t-Y^{t',x',\nu',z',\veca,\vecb}_{t'}|\to 0, \quad\text{as } (t',x',\nu',z')\to (t,x,\nu,z).
\end{equation*}
\end{prop}

\bigskip

\noindent\emph{Proof.} For the first part we note that, using Doob's maximal inequality and polynomial growth, there is a $p\geq 1$ such that
\begin{align*}
&\E\bigg[\sup_{s\in [0,T]}|Y^{t,x,\nu,z,\veca,\vecb}_s|^2\bigg]
\\
&\leq C\E\bigg[\int_0^T \max_{\veca\in\mcA}\max_{\zeta\in [0,T]^{\veca^+}}|\psi_{\veca}(v,X^{r,x}_v,\zeta)|^2 dv+\max_{\veca\in\mcA}\max_{\zeta\in [0,T]^{\veca^+}}|h_{\veca}(X^{r,x}_T,\zeta)|^2\bigg]
\\
&\leq C(1+|x|^{2p}).
\end{align*}
For the second part we argue as in the proof of Proposition~\ref{prop:Hk} and let $\hat u=(\hat \tau_1,\ldots,\hat \tau_{N};\hat \beta_1,\ldots,\hat \beta_N)\in \mcU_t^f$ be an optimal control corresponding to the reward $Y^{t,x,\nu,z,\veca,\vecb}_t$. Then,
\begin{align*}
Y^{t,x,\nu,z,\veca,\vecb}_t-Y^{t,x',\nu,z,\veca,\vecb}_{t}&\leq \E\bigg[\int_t^{T}\varphi_{\alpha_r^{t,x,\nu,\veca,\vecb,\hat u}}(r,X^{t,x}_r,\theta^{t,x,\nu,z,\veca,\vecb,\hat u}_r)-\varphi_{\alpha_r^{t,x',\nu,\veca,\vecb,\hat u}}(r,X^{t,x'}_r,\theta^{t,x',\nu,z,\veca,\vecb,\hat u}_r)dr
\\
&\quad +h_{\alpha_T^{t,x,\nu,\veca,\vecb,\hat u}}(X^{t,x}_T,\theta^{t,x,\nu,z,\veca,\vecb,\hat u}_T)
-h_{\alpha_T^{t,x',\nu,\veca,\vecb,\hat u}}(X^{t,x'}_T,\theta^{t,x',\nu,z,\veca,\vecb,\hat u}_T)\bigg],
\end{align*}
$\Prob$--a.s. Now, letting $\Psi_t^{x,x',\nu,z,\veca,\vecb,\hat u}$ denote the right hand side of the inequality we have
\begin{align*}
\Psi^{x,x',\nu,z,\veca,\vecb}_{t}&= \E\bigg[\int_t^{\hat\tau_1\wedge T}(e^{\int_t^{r}\rho^{\nu,\vecb}_{\veca,\veca}(v,X^{t,x}_v)dv}\varphi_{\veca}(r,X^{t,x}_r,\theta^{t,z,\veca}_r) -e^{\int_t^{r}\rho^{\nu,\vecb}_{\veca,\veca}(v,X^{t,x'}_v)dv}\varphi_{\veca}(r,X^{t,x'}_r,\theta^{t,z,\veca}_r))dr
\\
&\quad +\sum_{\veca'\in\mcA^{-\veca}_{\veca,\vecb}}\int_t^{\hat\tau_1\wedge T}(\rho^{\nu,\vecb}_{\veca,\veca'}(r,X^{t,x}_r)e^{\int_t^{r}\rho^{\nu,\vecb}_{\veca,\veca}(v,X^{t,x}_v)dv} Y^{r,X^{t,x}_r,\nu,\theta^{t,z,\veca}_r\wedge T\veca',\veca',\vecb}_r
\\
&\quad-\rho^{\nu,\vecb}_{\veca,\veca'}(r,X^{t,x'}_r)e^{\int_t^{r}\rho^{\nu,\vecb}_{\veca,\veca}(v,X^{t,x'}_v)dv} Y^{r,X^{t,x'}_r,\nu,\theta^{t,z,\veca}_r\wedge T\veca',\veca',\vecb}_r)dr
\\
&\quad + \ett_{[\hat\tau_1< T]}(e^{\int_t^{\hat\tau_1}\rho^{\nu,\vecb}_{\veca,\veca}(v,X^{t,x}_v)dv} Y^{{\hat\tau_1},X^{t,x}_{\hat\tau_1},\nu\hat\beta_1+{\hat\tau_1}(\hat\beta_1-\vecb)^+,\theta^{t,z,\veca}_{\hat\tau_1}\wedge T\hat\beta_1,\veca\wedge \hat\beta_1,\hat\beta_1}_{\hat\tau_1}
\\
&\quad-e^{\int_t^{\hat\tau_1}\rho^{\nu,\vecb}_{\veca,\veca}(v,X^{t,x'}_v)dv} Y^{{\hat\tau_1},X^{t,x'}_{\hat\tau_1},\nu\hat\beta_1+{\hat\tau_1}(\hat\beta_1-\vecb)^+,\theta^{t,z,\veca}_{\hat\tau_1}\wedge T\hat\beta_1,\veca\wedge \hat\beta_1,\hat\beta_1}_{\hat\tau_1})
\\
&\quad +\ett_{[\hat\tau_1\geq T]}(e^{\int_t^{T}\rho^{\nu,\vecb}_{\veca,\veca}(v,X^{t,x}_v)dv}h_{\veca}(X^{t,x}_T,\theta^{t,z,\veca}_T)-e^{\int_t^{T}\rho^{\nu,\vecb}_{\veca,\veca}(v,X^{t,x'}_v)dv}h_{\veca}(X^{t,x'}_T,\theta^{t,z,\veca}_T))\bigg].
\end{align*}
For $K>0$ we let $B:=\{\omega\in\Omega:\sup_{s\in[0,T]}(|X^{t,x}_s|+|X^{t,x'}_s|)\leq K\}$. Using the identity $ab-a'b'=1/2((a-a')(b+b')+(a+a')(b-b'))$ in combination with local Lipschitz continuity and polynomial growth gives the recursion
\begin{align*}
\Psi^{x,x',\nu,z,\veca,\vecb,\hat u}_{t}&\leq \E\bigg[C(K)\ett_{B}\Big(\int_t^{\hat\tau_1\wedge T}(1+|X^{t,x}_r|^p+|X^{t,x'}_r|^p)|X^{t,x}_r-X^{t,x'}_r|dr
\\
&\quad + \ett_{[\hat\tau_1\geq T]}(1+|X^{t,x}_T|^p+|X^{t,x'}_T|^p)|X^{t,x}_T-X^{t,x'}_T|
\\
&\quad +\sum_{\veca'\in\mcA^{-\veca}_{\veca,\vecb}}\int_t^{\hat\tau_1\wedge T}\Psi^{X^{t,x}_r,X^{t,x'}_r,\nu,\theta^{t,z,\veca}_{r}\wedge T\veca',\veca',\vecb,\hat u}_{r}dr\Big)
\\
&\quad +\ett_{B}\ett_{[\hat\tau_1< T]}\Psi^{X^{t,x}_{\hat\tau_1},X^{t,x'}_{\hat\tau_1},\nu\hat\beta_1+{\hat\tau_1}(\hat\beta_1-\vecb)^+,\theta^{t,z,\veca}_{\hat\tau_1}\wedge T\hat\beta_1,\wedge T\hat\beta_1,\hat\beta_1\hat u^{-1}}_{\hat\tau_1}\bigg]
\\
&\quad + C\E\bigg[\ett_{B^c}\Big(\int_t^{T}(1+|X^{t,x}_r|^p+|X^{t,x'}_r|^p)dr+1+|X^{t,x}_T|^p+|X^{t,x'}_T|^p\Big)\bigg],
\end{align*}
where $\hat u^{-1}:=(\hat \tau_2,\ldots,\hat \tau_{N};\hat \beta_2,\ldots,\hat \beta_N)$ and the constant $C(K)$ does not depend on the parameters $(t,x,x',z)$. Using an induction argument and finiteness of the strategy $\hat u$ together with H\"older's inequality we conclude that
\begin{align*}
Y^{t,x,\nu,z,\veca,\vecb}_t-Y^{t,x',\nu,z,\veca,\vecb}_{t}&\leq C(K)\E\bigg[\int_t^{T}|X^{t,x}_r-X^{t,x'}_r|dr+|X^{t,x}_T-X^{t,x'}_T|\bigg]
\\
&\quad+C\Prob(B^c)^{1/2}(1+|x|^{2p}+|x'|^{2p})^{1/2}.
\end{align*}
We note that, for all $K>0$ the first term approaches zero as $x'\to x$, by \eqref{ekv:Xbound2}. Furthermore, by \eqref{ekv:Xbound1}, $\Prob(B^c)\leq C\frac{1+|x|+|x'|}{K}$ and by symmetry we conclude that $\lim_{x'\to x}|Y^{t,x,\nu,z,\veca,\vecb}_t-Y^{t,x',\nu,z,\veca,\vecb}_{t}|\leq C\frac{1+|x|^{2p}+|x'|^{2p}}{K^{1/2}}$. Since $K>0$ was arbitrary, we conclude that $|Y^{t,x,\nu,z,\veca,\vecb}_t-Y^{t,x',\nu,z,\veca,\vecb}_{t}|\to 0$ as $x'\to x$. The second statement now follows by the $\Prob$-a.s.~continuity of $Y^{t,\nu,z,\veca,\vecb}_t$ in $(t,\nu,z)$ shown in Section~\ref{sec:exist}.\qed

%%%%%%%%%%%%%%%%%%%%%%%%%%%%%%%%%%%%%%%%%%%%%%%%%%%%%%%%%%%%%%%%%%%%%%%%%%%%%%%%%%%%%%%%%%%%%%%%%%%%%%%%%%%%%%%%%%%%%%%%%%%%%%%%%%%%%%%%%%%%%%%%%%

\bigskip

Now, \eqref{ekv:YdefALT} and repeated use of Theorem 8.5 in \cite{ElKaroui1} shows that for $k\geq 0$, there exist functions $(v_{(\veca,\vecb)}^{k})_{(\veca,\vecb)\in\mcJ}$ of polynomial growth, with $v_{(\veca,\vecb)}^{k}:[0,T]\times \R^m \times [0,T]^{\vecb}\times [0,T]^{\veca^+}\to \R$ such that
\begin{equation*}
Y^{t,x,\nu,z,\veca,\vecb,k}_t=v_{(\veca,\vecb)}^{k}(t,x,\nu,z).
\end{equation*}
Furthermore, by Proposition~\ref{prop:VFbnd} the functions $v_{(\veca,\vecb)}^{k}$ are continuous. Repeating the steps in the proof of Theorem~\ref{thm:limVT} we find that the sequence $(v_{(\veca,\vecb)}^{k})^{k\geq 0}_{(\veca,\vecb)\in\mcJ}$ converges pointwise to functions $v_{(\veca,\vecb)}:[0,T]\times \R^m \times [0,T]^{\vecb}\times [0,T]^{\veca^+}\to \R$ and that
\begin{equation*}
Y^{t,x,\nu,z,\veca,\vecb}_t=v_{(\veca,\vecb)}(t,x,\nu,z).
\end{equation*}
Now, by Proposition~\ref{prop:VFbnd} we find that the functions $v_{\vecb}$ are continuous and of polynomial growth. Finally, the verification theorem implies that the functions $v_{\vecb}$ are value functions for the stochastic control problem posed above and satisfy the following dynamic programming relation:
\begin{align}\nonumber
&v_{(\veca,\vecb)}(t,x,\nu,z)=\esssup_{\tau \in \mcT_{t}^{\bbF}} \E\bigg[\int_t^{\tau\wedge T}e^{\int_t^r\rho^{\nu,\vecb}_{\veca,\veca}(v,X^{t,x}_v)dv}\varphi_{\veca}(r,X^{t,x}_r,\theta_r^{t,z,\veca})dr
\\
&+\ett_{[\tau \geq T]}e^{\int_t^T\rho^{\nu,\vecb}_{\veca,\veca}(v,X^{t,x}_v)dv}h_{\veca}(X^{t,x}_T,\theta_T^{t,z,\veca}) + \nonumber
\\
&+\sum_{\veca'\in\mcA_{\veca,\vecb}^{-\veca}}\int_{t}^{\tau\wedge T}\rho^{\nu,\vecb}_{\veca,\veca'}(r,X^{t,x}_r)e^{\int_t^r\rho^{\nu,\vecb}_{\veca,\veca}(v,X^{t,x}_v)dv} v_{(\veca',\vecb)}(r,X^{t,x}_r,\nu,\theta_r^{t,z,\veca}\wedge T\veca') dr\nonumber
\\
&+\ett_{[\tau < T]}e^{\int_t^\tau\rho^{\nu,\vecb}_{\veca,\veca}(v,X^{t,x}_v)dv}\max_{\beta\in\mcI}\left\{-c_{\vecb,\beta}(\tau) + v_{(\veca\wedge\beta,\beta)}(\tau,X^{t,x}_\tau,\beta\nu +\tau(\beta-\vecb)^+,\theta_\tau^{t,z,\veca}\wedge T\beta)\right\}\Big| \mcF_t\bigg].\label{ekv:VFdef}
\end{align}
Using this formulation we can approximate the value function using either a Markov-chain approximation~\cite{numSCbok} or by reinforcement learning techniques~\cite{SuttonRL} (see \eg \cite{CarmLud,RAid,UUganget,ObactMMOR} for previous applications in optimal switching).

\bibliographystyle{plain}
\bibliography{RandDEL_ref}

\begin{thebibliography}{10}

\bibitem{RAid}
R.~A\"{i}d, L.~Campi, N.~Langren\'e, and H.~Pham.
\newblock A probabilistic numerical method for optimal multiple switching
  problems in high dimension.
\newblock {\em SIAM J. Financial Math.}, 5(1):191--231, 2014.

\bibitem{RAid2}
R.~A\"{i}d, S.~Federico, H.~Pham, and B.~Villeneuve.
\newblock Explicit investment rules with time-to-build and uncertainty.
\newblock {\em J. Econom. Dynam. Control}, 51:240--256, 2015.

\bibitem{BarIlan95}
A.~Bar-Ilan and A.~Sulem.
\newblock Explicit solution of inventory problems with delivery lags.
\newblock {\em Math. Oper. Res.}, 20(3), 1995.

\bibitem{BrekkeOksendal}
K.~A. Brekke and B.~{\O}ksendal.
\newblock Optimal switching in an economic activity under uncertainty.
\newblock {\em SIAM J. Control Optim.}, 32(4):1021--1036, 1994.

\bibitem{Brennan}
M.~J. Brennan and E.~S. Schwartz.
\newblock Evaluating natural resource investments.
\newblock {\em J. Bus.}, 58:135--157, 1985.

\bibitem{Bruder}
B.~Bruder and H.~Pham.
\newblock Impulse control problem on finite horizon with execution delay.
\newblock {\em Stochastic Process. Appl.}, 2009.

\bibitem{CarmLud}
R.~Carmona and M.~Ludkovski.
\newblock Pricing asset scheduling flexibility using optimal switching.
\newblock {\em Appl. Math. Finance}, 15:405--447, 2008.

\bibitem{ChassElieKharr}
J.~F. Chassagneux, R.~Elie, and I.~Kharroubi.
\newblock A note on existence and uniqueness for solutions of multidimensional
  reflected bsdes.
\newblock {\em Electron. Commun. Probab.}, 16:120--128, 2011.

\bibitem{CvitKar}
J.~Cvitanic and I.~Karatzas.
\newblock Backwards stochastic differential equations and {Dynkin} games.
\newblock {\em The Annals of Probability}, 24(4):2024--2056, 1996.

\bibitem{BollanMSwitch1}
B.~Djehiche, S.~Hamad\'ene, and A.~Popier.
\newblock A finite horizon optimal multiple switching problem.
\newblock {\em SIAM J. Control Optim.}, 47(4):2751--2770, 2009.

\bibitem{ElAsri}
B.~El~Asri and S.~Hamad\'ene.
\newblock The finite horizon optimal multi-modes switching problem: The
  viscosity solution approach.
\newblock {\em Appl. Math. Optim.}, 60:213--235, 2009.

\bibitem{ElKarouiLN}
N.~El~Karoui.
\newblock {\em Les aspects probabilistes du contr{\^o}le stochastique}.
\newblock Ecole d'Et\'e de SaintFlour IX 1979. Lecture Notes in Math. Springer,
  Berlin., 1981.

\bibitem{ElKaroui1}
N.~El-Karoui, C.~Kapoudjian, E.~Pardoux, S.~Peng, and M.~C. Quenez.
\newblock Reflected solutions of backward {SDEs} and related obstacle problems
  for {PDEs}.
\newblock {\em Ann. Probab.}, 25(2):702--737, 1997.

\bibitem{EliKharrCV}
Romuald Elie and Idris Kharroubi.
\newblock Bsde representations for optimal switching problems with controlled
  volatility.
\newblock {\em Stochastics and Dynamics}, 14(03), 2014.

\bibitem{HamRefBSDE}
S.~Hamad\`ene.
\newblock Reflected {BSDE's} with discontinuous barrier and application.
\newblock {\em Stochastics: An International Journal of Probability and
  Stochastic Processes}, 74(3-4):571--596, 2002.

\bibitem{HamJean}
S.~Hamad{\`e}ne and M.~Jeanblanc.
\newblock On the starting and stopping problem: application in reversible
  investments.
\newblock {\em Math. Oper. Res.}, 32(1):182--192, 2007.

\bibitem{HamZhang}
S.~Hamad\`ene and J.~Zhang.
\newblock Switching problem and related system of reflected backward {SDEs}.
\newblock {\em Stochastic Processes and their Applications}, 120(4):403--426,
  2010.

\bibitem{HuTang}
Y.~Hu and S.~Tang.
\newblock Multi-dimensional {BSDE} with oblique reflection and optimal
  switching.
\newblock {\em Prob. Theory and Related Fields}, 147(1-2):89--121, 2008.

\bibitem{KarShreve1}
I.~Karatzas and S.~E. Shreve.
\newblock {\em Brownian Motion and Stochastic Calculus}.
\newblock Springer-Verlag, New York, 1991.

\bibitem{KarShreve2}
I.~Karatzas and S.~E. Shreve.
\newblock {\em Methods of Mathematical Finance}.
\newblock Springer-Verlag, New York, 1998.

\bibitem{KharrSC}
I.~Kharroubi.
\newblock Optimal switching in finite horizon under state constraints.
\newblock {\em SIAM J. Control Optim.}, 54(4):2202--2233, 2016.

\bibitem{numSCbok}
H.~J. Kushner and P.~Dupuis.
\newblock {\em Numerical methods for stochastic control problems in continuous
  time}.
\newblock Springer, New York, 2nd edition, 2001.

\bibitem{Latifa}
I.~B. Latifa, J.~F. Bonnans, and M.~Mnif.
\newblock A general optimal multiple stopping problem with an application to
  swing options.
\newblock {\em Stochastic Analysis and Applications}, 33(4):715--739, 2015.

\bibitem{UUganget}
K.~Li, K.~Nystr\"om, and M.~Olofsson.
\newblock Optimal switching problems under partial information.
\newblock {\em Monte Carlo Methods and Applications}, 21(2):91--120, 2015.

\bibitem{LygerosSHS}
J.~Lygeros and M.~Prandini.
\newblock Stochastic hybrid systems: A powerful framework for complex, large
  scale applications.
\newblock {\em European Journal of Control}, 16(6):583--594, 2010.

\bibitem{MartyrSigned}
R.~Martyr.
\newblock Finite-horizon optimal multiple switching with signed switching
  costs.
\newblock {\em Math. Oper. Res.}, 41(4):1432--1447, 2016.

\bibitem{OksenImpulse}
B.~{\O}ksendal and A.~Sulem.
\newblock Optimal stochastic impulse control with delayed reaction.
\newblock {\em Appl. Math. Optim.}, 58:243--255, 2008.

\bibitem{LimFEED}
M.~Perninge.
\newblock A limited-feedback approximation scheme for optimal switching
  problems with execution delays.
\newblock {\em Math Meth Oper Res.}, 2017.

\bibitem{ObactMMOR}
M.~Perninge and L.~S\"oder.
\newblock Irreversible investments with delayed reaction: An application to
  generation re-dispatch in power system operation.
\newblock {\em Math. Meth. Oper. Res.}, 79:195--224, 2014.

\bibitem{Protter}
P.~Protter.
\newblock {\em Stochastic Integration and Differential Equations}.
\newblock Springer, Berlin, 2nd edition, 2004.

\bibitem{SuttonRL}
R.~S. Sutton and A.~G. Barto.
\newblock {\em Reinforcement Learning: An Introduction}.
\newblock The MIT Press Cambridge, Massachusetts, 2017.

\bibitem{XYZbook}
J.~Yong and X.~Y. Zhou.
\newblock {\em Stochastic Controls: Hamiltonian Systems and {HJB} equations}.
\newblock Springer, New York, 1999.

\end{thebibliography}
\end{document}